\documentclass{article}
\usepackage{amsmath,amssymb,graphicx,color,url,mathptmx}
\allowdisplaybreaks[3]
\newtheorem{thm}{Theorem}[section]

\newtheorem{prop}[thm]{Proposition}
\newtheorem{cor}[thm]{Corollary}
\newtheorem{lem}[thm]{Lemma}
\newtheorem{rem}[thm]{Remark}

\numberwithin{equation}{section}
\begin{document}

\title{Extremes of local times for simple random walks on symmetric trees}
\author{Yoshihiro Abe \footnote{Graduate School of Science, Kobe University,
Rokkodai-cho 1-1, Nada-ku, Kobe, 657-8501, Japan.
\newline e-mail: \texttt{yosihiro@math.kobe-u.ac.jp}}}
\date{}
\maketitle

\begin{abstract}
We consider local times of the simple random walk on the $b$-ary tree of depth $n$
and study a point process which encodes the location of the vertex with the maximal local time
and the properly centered maximum
over leaves of each subtree of depth $r_n$ rooted at the $(n-r_n)$ level,
where $(r_n)_{n \geq 1}$ satisfies $\lim_{n \to \infty} r_n = \infty$
and $\lim_{n \to \infty} r_n/n \in [0, 1)$.
We show that the point process weakly converges to a Cox process
with intensity measure $\alpha Z_{\infty} (dx) \otimes e^{-2\sqrt{\log b}~y}dy$,
where $\alpha > 0$ is a constant
and $Z_{\infty}$ is a random measure on $[0, 1]$
which has the same law as the limit of a critical random multiplicative cascade measure
up to a scale factor.
As a corollary, we establish convergence in law of the maximum of local times
over leaves to a randomly shifted Gumbel distribution.
\end{abstract}

\section{Introduction} \label{sec:intro}
Much efforts have been made in the study of the so-called log-correlated random field
such as the branching Brownian motion (BBM), the branching random walk (BRW),
and the two-dimensional discrete Gaussian free field (DGFF).
One of the remarkable features of these models is that
laws of their maxima share common properties:
each of the laws weakly converges to a randomly shifted Gumbel distribution
\cite{LS, Ai, BDZ}.   
It is believed that each of the limiting extremal processes of a wide class of log-correlated fields
converges to a so-called randomly shifted decorated Poisson point process \cite{SZ}
and it is established for the BBM \cite{ABBS, ABK3},
the BRW \cite{Ma},
and the two-dimensional DGFF \cite{BL3}.

It is well-known that local times of random walks on graphs have close relationships
with DGFFs thanks to ``the generalized second Ray-Knight theorem" \cite{EKMRS}
(this goes back to the Dynkin isomorphism \cite{Dy})
which has many applications, for example, to the cover time \cite{DLP, Di, Z}.
Since the occupation time field of the simple random walk
on the tree or on the two-dimensional lattice
is closely related to the BRW or two-dimensional DGFF respectively, 
it is natural to expect that their maxima and cover times belong to the universal class
mentioned above:
it is known that the cover times have subleading terms similarly to
other log-correlated fields \cite{DZ1, BK}
and that the cover time of the binary tree is tight \cite{BZ, BRZ},
but further details are still open.

In this paper, we consider local times
of the simple random walk on the $b$-ary tree of depth $n$
at time much larger than the maximal hitting time
and study convergence of a point process encoding 
extreme local maxima of the local times as $n \to \infty$.
 
To state our result, we begin with some notation.
We fix an arbitrary integer $b \ge 2$ throughout the paper.
We will write $T$ 
to denote the $b$-ary tree with root $\rho$:
this is a rooted tree whose vertices have exactly $b$ children.
Let $T_i$ be the $i$th generation of $T$.
Set $T_{\leq n} := \cup_{i=0}^n T_i.$ 
For $v \in T$, we will write $|v|$ to denote the depth of $v$.
For $u \in T$, let $T^u$ be the subtree of $T$ rooted at $u$,
and we define $T_i^u$ and $T_{\leq n}^u$ similarly.
For $v, u \in T$, let $v \wedge u$ be the most recent common ancestor of $v$ and $u$.
Let $X = \left(X_t,~t \geq 0,~P_v,~v \in T_{\leq n} \right)$ be the continuous-time 
simple random walk on $T_{\leq n}$
with exponential holding times of parameter $1$.
We define the local time of $X$ by
\begin{equation*}
L_t^n (v) := \frac{1}{\text{deg} (v)}\int_0^t 1_{\{X_s = v \}} ds,~~v \in T_{\leq n},~ t \ge 0,
\end{equation*}
where $\text{deg}(v)$ is the degree of $v$,
and the inverse local time by 
\begin{equation*}
\tau (t) := \inf \{s \ge 0 : L_s^n (\rho) > t \},~~~t \ge 0.
\end{equation*}
Let $E(T)$ be the set of all edges on $T$.
Let $(Y_e)_{e \in E(T)}$ be independent and identically distributed random variables 
whose common law is the normal distribution with mean $0$ and variance $1/2$.  
To each $v \in T$, we assign $h_v := \sum_{i=1}^{|v|} Y_{e_i^v}$,
where $e_1^v, \dotsc, e_{|v|}^v$ are the edges on the unique shortest path from $\rho$ to $v$.
We will call $(h_v)_{v \in T}$ a BRW on $T$.
It is well-known that the so-called derivative martingale
\begin{equation*} 
D_n := \sum_{v \in T_n} \left(\sqrt{\log b}~n - h_v \right) 
e^{- 2 \sqrt{\log b} \left(\sqrt{\log b}~n - h_v \right)}
\end{equation*}
converges almost surely as $n \to \infty$, and
the limit
\begin{equation} \label{eq:limit-derivative-martingale}
D_{\infty} := \lim_{n \to \infty} D_n
\end{equation}
is positive and finite almost surely
(see, for example, \cite[Theorem 5.1, 5.2]{BiKy} or \cite[Proposition A.3]{Ai}).
To each $v \in T$, we assign a distinct label $(\overline{v}_1, \dotsc, \overline{v}_{|v|})$ 
with $\overline{v}_i \in \{0, \dotsc, b-1 \}$ for all $1 \leq i \leq |v|$
so that the vertices with labels $(\overline{v}_1, \dotsc, \overline{v}_{|v|}, k)$, 
$k \in \{0, \dotsc, b-1 \}$ are children of $v$.
We define the location of $v \in T$ by
\begin{equation} \label{eq:location-point}
\sigma (v) := \sum_{i=1}^{|v|} \frac{\overline{v}_i}{b^i}.
\end{equation}
For each $n \in \mathbb{N}$ and 
$x \in [0, 1]$,
let $v(x)$ be the vertex in $T_n$ with
$x \in [\sigma (v(x)), \sigma (v (x)) + b^{- n}].$
(If we have two such vertices, we choose the one whose location is the largest.)
We define the random measure called a (critical) random multiplicative cascade measure by
\begin{equation} \label{eq:critical-cascade-measure}
Z_n (dx) := b^n \left(\sqrt{\log b}~ n - h_{v(x)} \right)
e^{- 2 \sqrt{\log b} \left(\sqrt{\log b}~ n - h_{v(x)} \right)} dx,
\end{equation}
where $dx$ is the Lebesgue measure on $[0, 1]$.
Barral, Rhodes, and Vargas \cite{BRV} observed that
\begin{equation} \label{eq:critical-mandelbrot-cascade}
\text{the weak limit}~Z_{\infty} := \lim_{n \to \infty} Z_n~\text{exists almost surely}.
\end{equation}
For each $v \in T_n$, set
$I_v := [\sigma (v), \sigma (v) + b^{-n}]$.
The random measure
$Z_{\infty}$ satisfies that
\begin{equation*}
\left(Z_{\infty} (I_v) \right)_{v \in T_n}
\stackrel{d}{=}
\left(e^{- 2 \sqrt{\log b} \left(\sqrt{\log b}~n - h_{v} \right)}
D_{\infty}^{(v)} \right)_{v \in T_n},
\end{equation*}
where 
``$\stackrel{d}{=}$'' means that the laws of the left and the right are the same
and $D_{\infty}^{(v)}, v \in T_n$ are independent copies of $D_{\infty}$
which are independent of $(h_v)_{v \in T_n}$. 
See \cite{BKNSW, BDK} for more details
on $Z_{\infty}$.
For each $(x, y) \in [0, 1] \times \mathbb{R}$, we write $\delta_{(x, y)}$
to denote the Dirac measure at $(x, y)$.
For each $0 \leq m \leq n$, we define the point process on $[0, 1] \times \mathbb{R}$ by
\begin{equation} \label{eq:point-process}
\Xi_{n, t}^{(m)} := 
\sum_{u \in T_{n - m} }
\delta_{\left(\sigma \left(\text{argmax}_{u}~L_{\tau (t)}^n \right),
~\max_{v \in T_{m}^u} \sqrt{L_{\tau (t)}^n (v)} - \sqrt{t} - a_n (t) \right)},
\end{equation}
where the centering sequence $a_n (t)$ is given by
\begin{equation} \label{eq:value-max}
a_n (t) := \sqrt{\log b}~ n - \frac{3}{4 \sqrt{\log b}} \log n - \frac{1}{4 \sqrt{\log b}} 
\log \left(\frac{\sqrt{t} + n}{\sqrt{t}} \right),
\end{equation} 
and for each $u \in T_{n - m}$, 
$\text{argmax}_{u}~L_{\tau (t)}^n$ is the vertex $v_*$ on 
$T_{m}^u \subset T_n$ with
$L_{\tau (t)}^n (v_*) = \max_{v \in T_{m}^u} L_{\tau (t)}^n (v)$.
(If two or more vertices on $T_{m}^u$ attain the maximum,
we take the one
whose location is the largest among such vertices.)
We regard $\Xi_{n, t}^{(m)}$ as an element of all Radon measures on Borel sets of 
$[0, 1] \times \mathbb{R}$
topologized with the vague topology.
Since this space is metrizable as a complete separable metric space,
we can consider convergence in law of sequences of random measures.
Given a random measure $\nu$ on $[0, 1] \times \mathbb{R}$,
we will write $\text{PPP}(\nu)$ to denote a point process on $[0, 1] \times \mathbb{R}$ 
which, conditioned on $\nu$, is a Poisson point process with intensity measure $\nu$
(that is $\text{PPP}(\nu)$ is a Cox process).
We now state the main result of this paper:
\begin{thm} \label{thm:convergence-point-process}
There exists $c_1 > 0$ such that 
for any sequence $(t_n)_{n \ge 1}$ with
$\lim_{n \to \infty} \frac{\sqrt{t_n}}{n} = \theta \in [0, \infty]$
and $t_n \ge c_1 n\log n$ for each $n \in \mathbb{N}$,
and any sequence $(r_n)_{n \geq 1}$ with $\lim_{n \to \infty} r_n = \infty$
and $\lim_{n \to \infty} r_n/n \in [0,~1)$,
the point process $\Xi_{n, t_n}^{(r_n)}$ converges in law to a Cox process
\begin{equation} \label{eq:limit-cox-process}
\text{PPP} \left( \frac{4}{\sqrt{\pi}} \beta_{*} \gamma_{*} Z_{\infty} (dx) 
\otimes 2 \sqrt{\log b}~ e^{- 2 \sqrt{\log b} ~y} dy \right)
\end{equation}
as $n \to \infty$,
where $Z_{\infty}$ is the random measure on $[0, 1]$ 
in (\ref{eq:critical-mandelbrot-cascade}),
\begin{equation} \label{eq:const-beta}
\beta_{*} := \begin{cases}
                           \sqrt{\frac{\theta + 1}{\theta + \sqrt{\log b}}} 
& \text{if} ~~\theta \in [0, \infty), \\ 
                             \\
                            1 & \text{if}~~\theta = \infty,
                            \end{cases}
\end{equation}
and
\begin{equation} \label{eq:const-gamma}
\gamma_{*} 
:= \lim_{\ell \to \infty} \int_{\ell^{2/5}}^{\ell} z e^{2 \sqrt{\log b}~z} 
\mathbb{P} \left(\max_{v \in T_{\ell}} h_v > \sqrt{\log b} ~\ell + z \right) dz.
\end{equation}
\end{thm}
\begin{rem} \label{rem:convergence-gamma}
The existence of the limit (\ref{eq:const-gamma}) is non-trivial.
It is proved in the proof of Proposition \ref{prop:limit-tail}.
Results similar to Theorem \ref{thm:convergence-point-process} are known for the BBM \cite{ABK2}
and the two-dimensional DGFF \cite{BL1, BL2}.
Our setting is inspired by \cite{BH}.
The convergence of the full extremal process has been established 
for the BBM \cite{ABBS, ABK3}, the BRW \cite{Ma},
and the two-dimensional DGFF \cite{BL3}.
Related convergence for the local times on the $b$-ary tree will be studied in a sequel paper.
\end{rem}
By Theorem \ref{thm:convergence-point-process}
and a tail estimate of the maximum of local times over leaves (Proposition \ref{prop:tail}(i) below), 
we have:
\begin{cor} \label{cor:convergence-in-law}
There exists $c_1 > 0$ such that for all $\lambda \in \mathbb{R}$ and any 
sequence $(t_n)_{n \ge 1}$ with
$\lim_{n \to \infty} \frac{\sqrt{t_n}}{n} = \theta \in [0, \infty]$
and $t_n \ge c_1 n\log n$ for each $n \in \mathbb{N}$, 
\begin{equation} \label{eq:convergence-in-law}
\lim_{n \to \infty} 
P_{\rho} \left(\max_{v \in T_{n}} \sqrt{L_{\tau (t_n)}^n (v)} \le \sqrt{t_n} + a_n (t_n) + \lambda \right)
= \mathbb{E} \left[ e^{- \frac{4}{\sqrt{\pi}} \beta_{*} \gamma_{*} D_{\infty} 
e^{- 2 \sqrt{\log b} ~\lambda} } \right],
\end{equation}
where $D_{\infty}$, $\beta_{*}$ and $\gamma_{*}$ are given 
by (\ref{eq:limit-derivative-martingale}), (\ref{eq:const-beta}) and (\ref{eq:const-gamma}), respectively.
\end{cor}
\begin{rem} \label{rem:relation-convergence-local-time-brw}
Let $(h_v)_{v \in T}$ be a BRW on $T$.
By Theorem \ref{thm:iso-thm} and Lemma \ref{lem:tail-max-brw-1}, one can show that
for all $\lambda \in \mathbb{R}$ and
any sequence $(t_n)_{n \geq 1}$ with $\lim_{n \to \infty} \sqrt{t_n}/n^2 = \infty$,
\begin{align} \label{eq:relation-convergence-local-time-brw}
&~~~~\lim_{n \to \infty} \mathbb{P} \left(\max_{v \in T_n} h_v \leq m_n + \lambda \right) 
\notag \\
&= \lim_{n \to \infty} 
P_{\rho} \left(\max_{v \in T_{n}} \sqrt{L_{\tau (t_n)}^n (v)} \le \sqrt{t_n} + m_n + \lambda \right)
= \mathbb{E} \left[ e^{- \frac{4}{\sqrt{\pi}} \gamma_{*} D_{\infty} 
e^{- 2 \sqrt{\log b} \lambda} } \right],
\end{align}
where $\gamma_*$ is given in (\ref{eq:const-gamma}) and the centering sequence $m_n$ is defined by
\begin{equation*}
m_n := \sqrt{\log b}~ n - \frac{3}{4 \sqrt{\log b}} \log n.
\end{equation*}
(Note that the convergence of the maximum of the BRW has 
already been established in \cite{Ba, Ai, BDZ2}.)
The centering sequence $a_n (t_n)$ in (\ref{eq:convergence-in-law})
is different from $m_n$ 
by the term $\frac{1}{4 \sqrt{\log b}} \log \left(\frac{\sqrt{t_n} + n}{\sqrt{t_n}} \right)$
which is non-negligible only when $\theta < \infty$.
\end{rem}
The organization of the paper is as follows.
Section \ref{sec:metric-tree} gives preliminary lemmas which we use repeatedly throughout the paper.
In Section \ref{sec:tail}, we obtain tail probabilities of the maximum of local times over leaves
which are essential to next sections. 
One can find that
for each leaf $v$,
the law of
the local time process along the path from $\rho$ to $v$
is the same as that of a zero-dimensional squared Bessel process
(see Lemma \ref{lem:local-time-bessel}).
By this and the Markov property of local time processes (see Lemma \ref{lem:markov-local-time}),
roughly speaking, one can regard the field of local times over the set of leaves
as a {\em branching Bessel process}.
This gives hints of how to estimate the tail of the maximum of local times over leaves:
we use the {\em constraint} first and second moment methods
developed in the BBM, BRW, and two-dimensional DGFF settings.
(See, for example, \cite{Br1, Ai, BDZ}.
We especially use techniques in \cite{DZ, BDZ}.)
Typical behavior of a vertex with extreme local time
is as follows:
the local time process along the path from the root to the vertex
stays below a curve and finally reach the maximal value at the vertex
(see the proof of Proposition \ref{prop:tail}).
In Section \ref{sec:geometry}, we show that
two leaves with local times near maxima are either very close or far away.
This suggests that local maximizers are distributed as a Poisson nature.
More technically, this implies that
$\Xi_{n, t_n}^{(n - q)}|_{[0, 1] \times [z,~\infty)} 
= \Xi_{n, t_n}^{(r_n)}|_{[0, 1] \times [z,~\infty)}$
with probability tending to $1$ as $n \to \infty$ and then $q \to \infty$,
which is one of the key steps in the proof of Theorem \ref{thm:convergence-point-process}.
In Section \ref{sec:limit-tail}, we obtain a limiting tail of the maximum of local times over leaves
which is crucial to study the Laplace functional of $\Xi_{n, t_n}^{(n-q)}$.
In the estimate, {\em entropic repulsion} (Lemma \ref{lem:density-barrier}(ii)) plays an important role:
this enables us to compute the tail of the maximum
by using the reflection principle of a Brownian motion.
In Section \ref{sec:pf-convergence-point-process}, 
we give the proof of Theorem \ref{thm:convergence-point-process} 
and Corollary \ref{cor:convergence-in-law}.

We should emphasize that it is more convenient to study ``continuous'' version of local times
rather than the original ``discrete'' ones especially when we estimate tail probabilities
of the maximum of the local times over leaves.
To take the advantage, motivated by \cite{Lu, Z}, we consider 
the local time process of the Brownian motion on the associated metric tree
as the ``continuous'' version.

We will write $c_1, c_2, \dotsc$ to denote positive universal constants
whose values are fixed within each argument.
We use $c_1 (M), c_2 (M) \dotsc$ for positive constants
which depend on $M$.
Given sequences $(c_n)_{n \ge 1}$ and $(c_n^{\prime})_{n \ge 1}$, we write 
$c_n^{\prime} = O(c_n)$
if there exists a universal constant $C$ such that $|c_n^{\prime}/c_n| \leq C$ for all $n \geq 1$.
We write $|S|$ to denote the cardinality of a set $S$.

\section{Preliminary lemmas} \label{sec:metric-tree}
In this section, we collect some lemmas which we use repeatedly throughout the paper.
We first recall the metric tree and the Brownian motions on it.
In the study of local times of random walks on graphs,
Lupu \cite{Lu} and Zhai \cite{Z} used
the corresponding metric graphs and Brownian motions.
We follow the approach and find its advantages
in obtaining precise tail probabilities of the maximum of local times over leaves on the $b$-ary tree.
Given a graph $G$, we will write $E (G)$ to denote the edge set of $G$. 
Let $T$ be the $b$-ary tree.
We regard
each $e \in E(T)$ as an interval of length $1/2$ by setting
$I_e := \{e \} \times \left(0, \frac{1}{2} \right)$.
Set $\bar{I}_e := I_e \cup \{e^{-}, e^{+} \}$, where $e^{-}, e^{+} \in T$ be the endpoints 
of the edge $e$. 
Let $\pi_e$ be the map from $I_e$ to $\left(0, \frac{1}{2} \right)$ defined by
$\pi_e ((e, x)) := x$.
We extend $\pi_e$ to the map from $\bar{I}_e$ to $\left[0, \frac{1}{2} \right]$ by setting 
$\pi_e (e^{-}) := 0, ~\pi_e (e^{+}) := \frac{1}{2}$.
We define a metric tree of depth $n$ by
\begin{equation*}
\widetilde{T}_{\leq n} := T_{\leq n} \cup \bigcup_{e \in E(T_{\leq n})} I_e.
\end{equation*}
For each $k \in \mathbb{N}$ and $v \in T$,
we will write $\widetilde{T}_{\leq k}^v$
to denote the metric tree corresponding to the subtree $T_{\leq k}^v$.
We define the metric $d(\cdot,~\cdot)$ on $\widetilde{T}_{\leq n}$ as follows:
for $x, y \in \widetilde{T}_{\leq n}$, 
let $e_x$ and $e_y$ be the edges
with $x \in I_{e_x}$ and $y \in I_{e_y}$, respectively.
In the case $I_{e_x} \neq I_{e_y}$, we define $d (x, y)$ by
\begin{equation*} 
\min \left\{
|\pi_{e_x} (x) - \pi_{e_x} (v)| + \frac{1}{2} d_g (v, u) + 
|\pi_{e_y} (u) - \pi_{e_y} (y)| : v \in \{e_x^{-}, e_x^{+} \}, 
u \in \{e_y^{-}, e_y^{+} \} \right \},
\end{equation*}
where $d_g$ is the graph distance on $T_{\leq n}$.
In the case $I_{e_x} = I_{e_y}$, we set
$d(x, y) := |\pi_{e_x} (x) - \pi_{e_y} (y)|$.
We define a measure $m$ on $\widetilde{T}_{\leq n}$ by
$$m(dx) := \sum_{e \in E(T_{\leq n})} 1_{I_e} (x) \nu_e (dx),$$
where $\nu_e := \nu \circ \pi_e$, and $\nu$ is the Lebesgue measure on $(0, 1/2)$.
We have a $m$-symmetric Hunt process on $\widetilde{T}_{\leq n}$ with continuous sample paths
such that on each $I_e$, it behaves like a standard Brownian motion on $(0, 1/2)$
until it hits $\{e^{-}, e^{+} \}$, and when it starts at a vertex $v$, it chooses one of
the edges incident to $v$ uniformly at random, and moves on it as described above.
See, for example, \cite{Fo, KPS, Lu} for the construction.
We write
$\widetilde{X} = (\widetilde{X}_t, t \ge 0, ~\widetilde{P}_x, x \in \widetilde{T}_{\leq n})$
to denote the process and call it a Brownian motion on $\widetilde{T}_{\leq n}$.
It is known that
$\widetilde{X}$ restricted to $T_{\leq n}$ behaves like a simple random walk on $T_{\leq n}$
in the following sense: for all $v \in T_{\leq n}$ and $1 \leq i \leq \text{deg} (v)$,
\begin{equation} \label{eq:rw-property}
\widetilde{P}_v (\widetilde{X}_{S_v} = v_i) = \frac{1}{\text{deg (v)}},
\end{equation}
where $v_1, \dotsc, v_{\text{deg} (v)}$ are vertices on $T_{\leq n}$ adjacent to $v$, and
$S_v$ is the hitting time of $\{v_1, \dotsc, v_{\text{deg} (v)} \}$ by $\widetilde{X}$.
See, for example, \cite[Theorem 2.1]{Fo} or \cite[Section 2]{Lu}.
By \cite[Section 2]{Lu}, $\widetilde{X}$ has a space-time continuous local time
$\{\widetilde{L}_t^n (x) : (t, x) \in [0, \infty) \times \widetilde{T}_{\leq n} \}$
and the following holds for each $v \in T_{\leq n}$ under $\widetilde{P}_v$:
\begin{equation} \label{eq:holding-time}
\widetilde{L}_{S_v}^n (v) \stackrel{d}{=}
~\text{Exp} \left(\frac{1}{\text{deg} (v)} \right),
\end{equation}
where $\text{Exp} (m)$ is an exponential random variable with mean $m$.
We define the inverse local time by
\begin{equation*}
\widetilde{\tau} (t) := \inf \{s \ge 0 : \widetilde{L}_s^n (\rho) > t\}, ~~~t > 0.
\end{equation*}
By (\ref{eq:rw-property}) and (\ref{eq:holding-time}), we have
\begin{equation} \label{eq:conti-disc-local-time}
\left(\widetilde{L}_{\widetilde{\tau} (t)}^n (v) \right)_{v \in T_{\leq n}} 
\text{under}~\widetilde{P}_{\rho}
\stackrel{d}{=}
\left(L_{\tau (t)}^n (v) \right)_{v \in T_{\leq n}} \text{under}~P_{\rho}.
\end{equation}
The following is the Markov property of local times of the Brownian motion on $\widetilde{T}_{\leq n}$.
The discrete version can be found in \cite[Lemma 2.6]{Di}.
\begin{lem} \label{lem:markov-local-time}
Fix $n \in \mathbb{N}, t > 0$, and $a \in T_{\leq n} \backslash T_{n}$.
Let $\mathcal{F}^{\uparrow}$ be the $\sigma$-field generated by
$\widetilde{L}_{\widetilde{\tau} (t)}^n (x),
~x \in \{a \} \cup \widetilde{T}_{\leq n} \backslash \widetilde{T}_{\leq n - |a|}^a$.
Then, the law of
$\left\{\widetilde{L}_{\widetilde{\tau} (t)}^n (x) : 
x \in \widetilde{T}_{\leq n - |a|}^a \backslash \{a \} \right\}$
under $\widetilde{P}_{\rho} (\cdot | \mathcal{F}^{\uparrow})$
is the same as that of
$\left\{\widetilde{L}_{\widetilde{\tau}^{\downarrow} 
\left(\widetilde{L}_{\widetilde{\tau} (t)}^n (a) \right)}^{\downarrow} (x) 
: x \in \widetilde{T}_{\leq n - |a|}^a \backslash \{a \} \right\}$
under $\widetilde{P}_a$,
where $\left\{\widetilde{L}_s^{\downarrow} (x) : (s, x) 
\in [0, \infty) \times \widetilde{T}_{\leq n - |a|}^a \right\}$
is a local time of a Brownian motion on $\widetilde{T}_{\leq n - |a|}^a$ and 
$\widetilde{\tau}^{\downarrow} (s) := \inf \{r: \widetilde{L}_r^{\downarrow} (a) > s \}$.
\end{lem}
The proof of Lemma \ref{lem:markov-local-time} is given in Section \ref{subsec:pf-markov-local-time}.

The generalized second Ray-Knight theorem connects local times
and BRWs:
\begin{thm} (\cite{Z}) \label{thm:iso-thm}
For all $t > 0$ and $n \in \mathbb{N}$, on the same probability space, one can construct a local time
$(L_{\tau (t)}^n (v))_{v \in T_{\leq n}}$ and two BRWs $(h_v)_{v \in T_{\leq n}}$,
$(h_v^{\prime})_{v \in T_{\leq n}}$ on $T_{\leq n}$ satisfying the following:  
\begin{equation} \label{eq:iso-thm-1}
(L_{\tau (t)}^n (v))_{v \in T_{\leq n}}~ \text{and}~ (h_v)_{v \in T_{\leq n}}~ \text{are independent},
\end{equation}
\begin{equation} \label{eq:iso-thm-2}
L_{\tau (t)}^n (v) + (h_v)^2 = \left(h_v^{\prime} + \sqrt{t} \right)^2,
~~\text{for each}~v \in T_{\leq n},~\text{almost surely}.
\end{equation}
\end{thm}
The construction of the coupling in Theorem \ref{thm:iso-thm} can be found in
the proof of Theorem 3.1 of \cite{Z}.
(Note that Zhai constructed the coupling in a more general setting and that in the context of \cite{Z}, 
the law of the DGFF on a $b$-ary tree is the same as 
that of our BRW scaled by $\sqrt{2}$.)
Let $C [0, ~\infty)$ be the space of real-valued continuous functions on $[0,~\infty)$
and $\mathcal{B}\left(C [0,~\infty)\right)$ be the $\sigma$-field generated by
cylinder sets in $C [0, ~\infty)$.
We have a nice connection between the local time and the $0$-dimensional squared Bessel process.
\begin{lem} (\cite[Lemma 7.7]{BK}) \label{lem:local-time-bessel}
For all $t > 0$ and $v \in T_n$,
\begin{equation*}
\left\{\widetilde{L}_{\widetilde{\tau} (t)}^n (v_s) : 0 \leq s \leq n \right\}~\text{under}~\widetilde{P}_{\rho}
\stackrel{d}{=}
\left\{\frac{1}{2} X_s : 0 \leq s \leq n \right\}~\text{under}~ \mathbb{Q}_{2t}^0,
\end{equation*}
where $\mathbb{Q}_{x}^d$ is a law 
on $\left(C [0, \infty),~\mathcal{B} \left(C [0, \infty) \right) \right)$
under which 
the coordinate process $\{X_s : s \geq 0 \}$ is 
a $d$-dimensional squared Bessel process started at $x$.
\end{lem}
Note that our setting is different from that of \cite[Lemma 7.7]{BK}.
Notwithstanding,
given Lemma \ref{lem:covariance-gff},
the proof of Lemma \ref{lem:local-time-bessel}
is almost the same as that of \cite[Lemma 7.7]{BK},
so we omit the proof of Lemma \ref{lem:local-time-bessel}.
It is known that the laws of $0$-dimensional and $1$-dimensional squared Bessel processes
are related to each other by the Radon-Nikodym derivative
\begin{equation} \label{eq:0-1-dim-bessel}
\frac{d\mathbb{Q}_x^0}{d\mathbb{Q}_x^1} \arrowvert_{\mathcal{F}_t \cap \{H_0 > t \}}
= \left(\frac{x}{X_t} \right)^{1/4} \exp \left( - \frac{3}{8} \int_0^t \frac{1}{X_s} ds \right),
\end{equation}
for all $t > 0$ and $x > 0$,
where $H_0 := \inf \{t \geq 0 : X_t = 0 \}$
and $\mathcal{F}_t$ is the $\sigma$-field generated by $\{X_s : s \leq t \}$. 
See, for example, \cite[(7.31)]{BK}. 
The transition semigroup $\{Q_t^0 : t \ge 0 \}$ of a $0$-dimensional squared Bessel process
is given by
\begin{equation} \label{eq:semigroup-bessel}
Q_t^0 (x, \cdot) = \exp \left(- \frac{x}{2t} \right) \delta_0 + \widetilde{Q}_t (x, \cdot),~~~x > 0,
\end{equation}
where $\delta_0$ is the Dirac measure at $0$.
$\widetilde{Q}_t (x, \cdot)$ in (\ref{eq:semigroup-bessel}) has the density
\begin{equation}
q_t^0 (x, y) = \frac{1}{2t} \sqrt{\frac{x}{y}} \exp \left(- \frac{x+y}{2t} \right) 
I_1 \left(\frac{\sqrt{xy}}{t} \right),~~~x, y \in (0, \infty),
\end{equation} 
where $I_1 (\cdot)$ is the modified Bessel function of the first kind
\begin{equation} \label{eq:modified-bessel-func}
I_1 (z) = \sum_{k=0}^{\infty} \left(\frac{z}{2} \right)^{2k+1} \frac{1}{k ! (k+1) !}.
\end{equation}
We will use the following asymptotic behavior of $I_1 (\cdot)$:
\begin{equation} \label{eq:asymp-modified-bessel-func}
I_1 (z) = \frac{e^z}{\sqrt{2 \pi z}} (1 + O(1/z)),~~\text{as}~z \to \infty.
\end{equation} 
See, for example, \cite[Chapter XI, \S1]{RY} or \cite[Section 2]{BM}
for the details on squared Bessel processes.

Let $(B_s, s \geq 0, P_x^B, x \in \mathbb{R})$
be a Brownian motion on $\mathbb{R}$ with variance $1/2$.
To estimate tail probabilities of the maximum of local times over leaves,
we frequently use the following. 
\begin{lem} (\cite[Lemma 3.6]{BDZ}) \label{lem:density-barrier}
Fix a constant $C > 0$. For $z > 0$ and $s > 0$, set
\begin{align} \label{eq:density-barrier-standard}
\mu_{s, z} (x) dx &:= 
P_0^B \left(B_s \in dx, 
B_r \leq z,~0 \leq \forall r \leq s \right) \notag \\
&~= \frac{1}{\sqrt{\pi s}} \left(e^{- \frac{x^2}{s}} - e^{- \frac{(2z-x)^2}{s}} \right) dx,~~x \leq z,
\end{align}
\begin{equation*}
\mu_{s, z}^* (x) dx :=
P_0^B \left(B_s \in dx,
B_r \leq z + z^{\frac{1}{20}} + C (r \wedge (s-r))^{\frac{1}{20}},~
0 \leq \forall r \leq s \right).
\end{equation*}
(i) There exists $c_1 > 0$ such that for any $z > 1$, $s > 0$, and $x \leq z + z^{\frac{1}{20}}$,
\begin{equation*}
\mu_{s, z}^* (x) \leq c_1 z (z + z^{\frac{1}{20}} - x) s^{- \frac{3}{2}} e^{- \frac{x^2}{s}}.
\end{equation*}
(ii) There exists $\delta_z$ with $\lim_{z \to \infty} \delta_z = 0$ such that
for all $z > 1$, $x \leq 0$, and $s \geq x^2 + z^2$,
\begin{equation*} 
\mu_{s, z}^* (x) \leq (1 + \delta_z) \mu_{s, z} (x).
\end{equation*}
\end{lem}

\section{Tail of maximum of local time over leaves} \label{sec:tail}
The aim of this section is to obtain the following tail estimates of the maximum of local times
of the simple random walk on the $b$-ary tree over leaves. 
Recall the definition of $a_n (t)$ from (\ref{eq:value-max}).
\begin{prop} \label{prop:tail}
(i) There exist $c_1, c_2 \in (0,~\infty)$ such that for all $t > 0, y \ge 0$, and $n \in \mathbb{N}$,
\begin{equation} \label{eq:prop-upper-tail-statement}
P_{\rho} \left(\max_{v \in T_n} \sqrt{L_{\tau (t)}^n (v)} \ge \sqrt{t} + a_n (t) + y \right)
\leq c_1 (1 + y) e^{- 2 \sqrt{\log b} ~y} e^{- c_2 \frac{y^2}{n}}.
\end{equation}
(ii) There exist $c_3 > 0$ and $n_0 \in \mathbb{N}$ such that 
for all $n \ge n_0$, $y \in [0, 2 \sqrt{n}]$, and $t \ge n$,
\begin{equation} \label{eq:prop-lower-tail-statement}
P_{\rho} \left(\max_{v \in T_n} \sqrt{L_{\tau (t)}^n (v)} \ge \sqrt{t} + a_n (t) + y \right)
\ge c_3 (1+y) e^{- 2 \sqrt{\log b} ~y}.
\end{equation}
\end{prop}

Given $v \in T$ and $s \in [0, |v|]$, let $v_s$ be the point on the unique path
from $\rho$ to $v$ with $d (\rho, v_s) = s/2$.
We first prove Proposition \ref{prop:tail}(i).
Fix $\kappa \in \left(\frac{1}{2\sqrt{\log b}}, \infty \right)$.
For $y > 0$ and $n \in \mathbb{N}$,
we define the event $G_y^n (t)$ by
\begin{equation} \label{eq:event-broken-barrier}
\Biggl \{\exists v \in T_n, \exists s \in [0, n]: 
\sqrt{\widetilde{L}_{\widetilde{\tau} (t)}^n (v_s)} 
\ge \sqrt{t} + \frac{a_n (t)}{n} s + \kappa (\log (s \wedge (n-s)))_+ + y + 1 \Biggr \},
\end{equation}
where $c_+ := \max \{c, 0 \}$.
We prove that $G_{y}^n (t)$ is a rare event,
that is, every local time process along the path from the root to a leaf stays below the curve
$s \mapsto \sqrt{t} + \frac{a_n (t)}{n} s + \kappa (\log (s \wedge (n-s)))_+ + y + 1$
with high probability:
\begin{lem} \label{lem:broken-barrier}
There exist $c_1, c_2 \in (0,~\infty)$ such that
for all $t > 0$, $n \in \mathbb{N}$, and $y \ge 0$,
\begin{equation} \label{eq:lem-broken-barrier-statement}
\widetilde{P}_{\rho} (G_y^n (t)) \leq c_1 (1 + y) e^{- 2\sqrt{\log b}~y} e^{- c_2 \frac{y^2}{n}}.
\end{equation}
\end{lem}
{\it Proof.} 
We first consider the case $y > M$, where $M$ is sufficiently large constant.
Set
\begin{equation} \label{eq:barrier-func}
g_{y, t, n} (s) := \sqrt{t} + \frac{a_n (t)}{n} s + \kappa (\log (s \wedge (n-s)))_+ + y + 1,
~~~0 \leq s \leq n,
\end{equation}
\begin{equation*}
m_{y, t, n} (j) := \sqrt{t} + \frac{a_n (t)}{n} j 
+ \kappa \min_{j \leq r \leq j+1} (\log (r \wedge (n-r)))_+ + y + 1,
~~~0 \leq j \leq n-1.
\end{equation*}
Recall the probability measure $\mathbb{Q}_x^d$ 
defined in Lemma \ref{lem:local-time-bessel} and set
\begin{equation*}
\tau := \inf \left\{s \geq 0 : \sqrt{X_s/2} \ge g_{y, t, n} (s) \right\},
\end{equation*}
where $X$ is a coordinate process.
Fix $\delta \in (0, 1)$.
By Lemma \ref{lem:local-time-bessel},
$\widetilde{P}_{\rho} (G_y^n (t))$ is bounded from above by
\begin{align} \label{eq:lem-broken-barrier-1}
&\widetilde{P}_{\rho} \Biggl(
\begin{minipage}{230pt} $0 \leq \exists j \leq n-1$, $\exists v \in T_{j+1}$,
$\exists s \in (j, j+1]$ : \\
$\sqrt{\widetilde{L}_{\widetilde{\tau} (t)}^n (v_r)} < g_{y, t, n} (r)$,
~$0 \leq  \forall r \leq j$, 
~$\sqrt{\widetilde{L}_{\widetilde{\tau} (t)}^n (v_s)} \ge g_{y, t, n} (s)$
\end{minipage}
\Biggr)  \notag \\
&\leq \sum_{j=0}^{n-1} b^{j+1}
\mathbb{Q}_{2t}^0 \Biggl(\tau \in (j, j+1], ~
\sqrt{X_{j+1}/2} \ge \delta m_{y, t, n} (j) \Biggr) \notag \\
&+ \sum_{j=0}^{n-1} b^{j+1}
\mathbb{Q}_{2t}^0 \Biggl(\tau \in (j, j+1],
~\sqrt{X_{j+1}/2} < \delta m_{y, t, n} (j) \Biggr) \notag \\
&=: \sum_{j=0}^{n-1} b^{j+1} I_{j}^{(1)}
+ \sum_{j=0}^{n-1} b^{j+1} I_{j}^{(2)}.
\end{align} 
Fix $0 \leq j \leq n-1$.
We first estimate $I_{j}^{(2)}$.
By the strong Markov property of a $0$-dimensional squared Bessel process
and (\ref{eq:semigroup-bessel}),
we have 
\begin{align} \label{eq:lem-broken-barrier-2}
&I_{j}^{(2)} = \mathbb{Q}_{2t}^0 \Biggl[ 1_{\{\tau \in (j, j+1) \}} 
\mathbb{Q}_{X_{\tau}}^0 
\left( \sqrt{X_{j+1-\tau}/2} < \delta m_{y, t, n} (j) \right) \Biggr] \notag \\
&= \mathbb{Q}_{2t}^0 \left[1_{\{\tau \in (j, j+1) \}} 
\exp \left(- \frac{X_{\tau}}{2(j+1-\tau)} \right) \right] \notag \\
&+ \mathbb{Q}_{2t}^0 \left[1_{\{\tau \in (j, j+1) \}} 
\int_0^{2 \delta^2 (m_{y, t, n} (j))^2} \frac{\sqrt{X_{\tau}/z}}{2(j+1-\tau)} 
e^{- \frac{X_{\tau} + z}{2(j+1-\tau)} }
I_1 \left(\frac{\sqrt{X_{\tau} z}}{j+1-\tau} \right) dz \right] \notag \\
&=: J_1 + J_2.
\end{align}
By the definition of $\tau$, we have
\begin{equation} \label{eq:lem-broken-barrier-3}
J_1 \leq \exp \{- (m_{y, t, n} (j))^2 \}.
\end{equation}
Assume that $\tau \in (j, j+1)$.
Recall the definition of $I_1$ from (\ref{eq:modified-bessel-func}). 
If $z \leq \frac{M (j+1-\tau)^2}{X_{\tau}}$, then we have
\begin{equation} \label{eq:lem-broken-barrier-4}
I_1 \left(\frac{\sqrt{X_{\tau}z}}{j+1-\tau} \right)
\leq \sum_{k=0}^{\infty} (\sqrt{M}/2)^{2k+1} \frac{1}{k ! (k+1)!} \leq c_1 (M).
\end{equation}
If $z > \frac{M(j+1-\tau)^2}{X_{\tau}}$, then 
by (\ref{eq:asymp-modified-bessel-func})
and the assumption that $M$ is sufficiently large,
we have
\begin{equation} \label{eq:lem-broken-barrier-5}
I_1 \left(\frac{\sqrt{X_{\tau}z}}{j+1-\tau} \right)
\leq c_2 \frac{e^{\frac{\sqrt{X_{\tau} z}}{j+1-\tau}}}{\sqrt{\frac{\sqrt{X_{\tau}z}}{j+1-\tau}}}.
\end{equation}
By (\ref{eq:lem-broken-barrier-4}) and (\ref{eq:lem-broken-barrier-5}),
we have
\begin{equation} \label{eq:lem-broken-barrier-6}
J_2 \leq c_3 (M) \left(\max_{j \leq r \leq j+1} g_{y, t, n} (r) \right)  
\exp \{ - c_4 (m_{y, t, n} (j))^2 \},
\end{equation}
where we have used the inequality
$(j+1-\tau)^{-1/2} e^{- \frac{(\sqrt{X_{\tau}} - \sqrt{z})^2}{2(j+1-\tau)}}
\leq e^{- c_4 (m_{y, t, n} (j))^2}$ for all $z \in [0, 2 \delta^2 (m_{y, t, n} (j))^2]$.
By (\ref{eq:lem-broken-barrier-3}) and (\ref{eq:lem-broken-barrier-6}), we have
\begin{equation} \label{eq:lem-broken-barrier-7}
\sum_{j=0}^{n-1} b^{j+1} I_{j}^{(2)}
\leq c_5(M) \sum_{j=0}^{n-1} b^{j+1} e^{- c_6 (\sqrt{t} + y + 1)^2} e^{- c_7 (\sqrt{t} + y + 1) j} 
\leq c_8 (M) e^{- c_6 (\sqrt{t} + y + 1)^2},
\end{equation}
where in the last inequality,
we have used the assumption that $y > M$ and $M$ is sufficiently large.

Next, we will estimate $I_{j}^{(1)}$.
Fix $1 \leq j \leq n-1$.
By (\ref{eq:0-1-dim-bessel}), $I_{j}^{(1)}$ is equal to
\begin{equation} \label{eq:lem-broken-barrier-8} 
\mathbb{Q}_{2t}^1 \left[ \left(\frac{2t}{X_{j+1}} \right)^{1/4} 
\exp \left(- \frac{3}{8} \int_0^{j+1} \frac{ds}{X_s} \right) 
1_{\left\{\tau \in (j, j+1], \sqrt{X_{j+1}/2} \geq \delta m_{y, t, n} (j) \right\} 
\cap \left\{H_0 > j+1 \right\}} \right].
\end{equation}
Recall that $(B_s, s \geq 0, P_x^B, x \in \mathbb{R})$
is a Brownian motion on $\mathbb{R}$ with variance $1/2$.
Since the law of a $1$-dimensional squared Bessel process
is the same as that of a square of a standard Brownian motion on $\mathbb{R}$,
(\ref{eq:lem-broken-barrier-8}) is bounded from above by $\sqrt{\frac{\sqrt{t}}{\delta m_{y, t, n} (j)}}$ times
\begin{align} \label{eq:lem-broken-barrier-8-2}
&P_{\sqrt{t}}^B \Biggl(B_r < g_{y, t, n} (r),~\forall r \in [0, j], 
~B_s \ge g_{y, t, n} (s),~\exists s \in [j, j+1] \Biggr) \notag \\
&\leq 
E_0^B \Biggl[1_{\left\{B_r
< g_{y, t, n} (r) - \sqrt{t},
~\forall r \in [0, j] \right\}} 
P_{B_j}^B \left(B_s
\geq g_{y, t, n} (j+s) - \sqrt{t},
~\exists s \in [0, 1] \right) \Biggr],
\end{align}
where we have used
the translation invariance and Markov property of $B$ in the last inequality.

Let $\widetilde{P}_j^B$ be the probability measure defined by
\begin{equation} \label{eq:change-measure}
\widetilde{P}_j^B (A) 
= E_0^B \left[1_A \exp \left \{\frac{2 a_n (t)}{n} B_j - \frac{(a_n (t))^2}{n^2} j \right \} \right],
~~~A \in \sigma(B_s : s \leq j).
\end{equation}
By the Girsanov theorem, under $\widetilde{P}_j^B$,
the process 
\begin{equation} \label{eq:bm-drift}
\left\{\widetilde{B}_s := B_s - \frac{a_n(t)}{n} s ~:~ 0 \leq s \leq j \right\}
\end{equation}
is a Brownian motion on $\mathbb{R}$ with variance $1/2$ started at $0$.
By the change of measure (\ref{eq:change-measure}),
the right of (\ref{eq:lem-broken-barrier-8-2}) is bounded from above by
\begin{align} \label{eq:lem-broken-barrier-9}
&\widetilde{E}_j^B \left[
\begin{minipage}{260pt} 
$e^{- \frac{2 a_n (t)}{n} \widetilde{B}_j
- \frac{(a_n (t))^2}{n^2} j } 
1_{\left\{\widetilde{B}_r < y + c_9 + \kappa \log (j \wedge (n-j))
+ c_{10} (r \wedge (j-r))^{1/20},~\forall r \in [0, j] \right\}}$
\\
$\times P_{\widetilde{B}_j}^B \left[\max_{0 \leq s \leq 1} B_s \geq y - c_{11}
+ \kappa \log (j \wedge (n-j)) \right]$
\end{minipage} \right] \notag \\
&\leq 
\widetilde{E}_j^B \left[
\begin{minipage}{260pt} 
$e^{- \frac{2 a_n (t)}{n} \widetilde{B}_j
- \frac{(a_n (t))^2}{n^2} j } 
1_{\left\{\widetilde{B}_r < y + c_9 + \kappa \log (j \wedge (n-j))
+ c_{10} (r \wedge (j-r))^{1/20},~\forall r \in [0, j] \right\}}$ \\
$\times 1_{\left\{\widetilde{B}_j < y - c_{11} + \kappa \log (j \wedge (n-j)) \right\}}$ \\
$\times P_{\widetilde{B}_j}^B \left[\max_{0 \leq s \leq 1} B_s \geq y - c_{11}
+ \kappa \log (j \wedge (n-j)) \right]$
\end{minipage} \right] \notag \\
&+
\widetilde{E}_j^B \left[
\begin{minipage}{260pt} 
$e^{- \frac{2 a_n (t)}{n} \widetilde{B}_j
- \frac{(a_n (t))^2}{n^2} j } 
1_{\left\{\widetilde{B}_r < y + c_9 + \kappa \log (j \wedge (n-j))
+ c_{10} (r \wedge (j-r))^{1/20},~\forall r \in [0, j] \right\}}$
\\
$\times 1_{\left\{\widetilde{B}_j \in [y - c_{11} + \kappa \log (j \wedge (n-j)),~y + c_9 
+ \kappa \log (j \wedge (n-j))] 
\right\}}$
\end{minipage} \right].
\end{align} 
To estimate the tail of $\max_{0 \leq s \leq 1} B_s$ in the first term of 
the right of (\ref{eq:lem-broken-barrier-9}),
we use the following:
\begin{equation} \label{eq:max-bm}
P_0^B \left(\max_{0 \leq s \leq 1} B_s \geq \lambda \right)
\leq e^{- \lambda^2},~~~\text{for each}~\lambda > 0
\end{equation}
(see, for example, \cite[Chapter 2, (8.4)]{KS}).
By Lemma \ref{lem:density-barrier}(i) and (\ref{eq:max-bm}),
the right of
(\ref{eq:lem-broken-barrier-9}) is bounded from above by
\begin{align} \label{eq:lem-broken-barrier-9-2}
&c_{12} e^{- \frac{2 a_n (t)}{n} y - \frac{(a_n (t))^2}{n^2} j}
\left(1 + \frac{y}{j+1} \right)(y + c_9 + \kappa \log (j \wedge (n-j)))  \notag \\
&\times
j^{- \frac{3}{2}}(j \wedge (n-j))^{- 2\sqrt{\log b}~\kappa} e^{- c_{13} \frac{y^2}{j+1}}.
\end{align}
Similarly, in the case $j=0$,
by (\ref{eq:0-1-dim-bessel}) and (\ref{eq:max-bm}), we have
\begin{equation} \label{eq:lem-broken-barrier-10}
I_{0}^{(1)} \leq c_{14} P_0^B \left(\max_{0 \leq s \leq 1} B_s \geq y+1 \right)
\leq c_{14} e^{- (y+1)^2}.
\end{equation}
Thus, by (\ref{eq:lem-broken-barrier-9-2}), (\ref{eq:lem-broken-barrier-10}),
and the condition $\kappa > 1/(2\sqrt{\log b})$,
we have
\begin{equation} \label{eq:lem-broken-barrier-11}
\sum_{j=0}^{n-1} b^{j+1} I_{j}^{(1)} 
\leq c_{15} (1+y) e^{- 2\sqrt{\log b}~y} e^{- c_{16} \frac{y^2}{n}}.
\end{equation}
Thus, by (\ref{eq:lem-broken-barrier-1}), (\ref{eq:lem-broken-barrier-7})
and (\ref{eq:lem-broken-barrier-11}), we have (\ref{eq:lem-broken-barrier-statement})
for $y > M$.
For $y \leq M$, (\ref{eq:lem-broken-barrier-statement}) holds
if we take $c_1$ in (\ref{eq:lem-broken-barrier-statement}) 
sufficiently large depending on $M$.~~~$\Box$ \\
We now prove Proposition \ref{prop:tail}(i). \\\\
{\em Proof of Proposition \ref{prop:tail}(i).}
Recall the definitions of the event $G_y^n (t)$ and the function $g_{y, t, n} (\cdot)$
from (\ref{eq:event-broken-barrier}) and (\ref{eq:barrier-func}). 
In view of Lemma \ref{lem:broken-barrier}, 
it is natural to impose the restriction that local time processes stay below the curve
$s \mapsto g_{y, t, n} (s)$:
we have
\begin{equation} \label{eq:prop-upper-tail-1}
P_{\rho} \left(\max_{v \in T_n} \sqrt{L_{\tau (t)}^n (v)} \ge \sqrt{t} + a_n (t) + y \right) 
\leq \widetilde{P}_{\rho} \left(\cup_{v \in T_n} E_v^n (t) \right) 
+ \widetilde{P}_{\rho} (G_y^n (t)),
\end{equation}
where for each $v \in T_n$, we define the event $E_v^n (t)$ by
\begin{equation*}
\left \{
\sqrt{\widetilde{L}_{\widetilde{\tau} (t)}^n (v_s)} \leq g_{y, t, n} (s),~\forall s \in [0, n], 
~\sqrt{\widetilde{L}_{\widetilde{\tau} (t)}^n (v)} 
\in \left[\sqrt{t} + a_n (t) + y, \sqrt{t} + a_n (t) + y + 1 \right] \right \}.
\end{equation*}
Fix any $v \in T_n$.
Recall the process $\widetilde{B}$ from (\ref{eq:bm-drift}).
By Lemma \ref{lem:local-time-bessel}, (\ref{eq:0-1-dim-bessel}), 
and the change of measure (\ref{eq:change-measure}) with $j=n$,
$\widetilde{P}_{\rho} (E_v^n (t))$ is bounded from above by
\begin{equation} \label{eq:prop-upper-tail-2}
\sqrt{\frac{\sqrt{t}}{\sqrt{t} + a_n (t) + y}}
\widetilde{E}_n^B \left[e^{- \frac{2 a_n (t)}{n} \widetilde{B}_n - \frac{(a_n (t))^2}{n}}
1_{\left\{\widetilde{B}_s \leq y + 1 + c_1 \kappa (s \wedge (n-s))^{\frac{1}{20}},
~\forall s \in [0, n], ~\widetilde{B}_n \in [y, y+1] \right\}} \right].
\end{equation}
By Lemma \ref{lem:density-barrier}(i), the right of (\ref{eq:prop-upper-tail-2})
is bounded from above by
\begin{equation} \label{eq:prop-upper-tail-2-2}
c_2 b^{-n} (1 + y) e^{- 2 \sqrt{\log b}~y} e^{- c_3 \frac{y^2}{n}}.
\end{equation}
Thus, by (\ref{eq:prop-upper-tail-1}), (\ref{eq:prop-upper-tail-2-2}),
and Lemma \ref{lem:broken-barrier}, we have (\ref{eq:prop-upper-tail-statement}).
~~~$\Box$ \\\\

Next, we prove Proposition \ref{prop:tail}(ii).
Fix $\delta \in (0, 1)$.
For $v \in T_n$, set the event
\begin{equation} \label{eq:event-optimal-barrier}
A_v^n (t) := \left \{
\begin{minipage}{210pt}
$\delta \sqrt{t} \leq \sqrt{\widetilde{L}_{\widetilde{\tau} (t)}^n (v_s)}
< \sqrt{t} + \frac{a_n (t)}{n} s + y + 1,~\forall s \in [0, n],$ \\
$\sqrt{\widetilde{L}_{\widetilde{\tau} (t)}^n (v)} \in [\sqrt{t} + a_n (t) + y, \sqrt{t} + a_n (t) + y + 1)$
\end{minipage}
\right \}.
\end{equation}
To obtain Proposition \ref{prop:tail}(ii),
we will apply the second moment method to $\sum_{v \in T_n} 1_{A_v^n (t)}$.
We first need the following:
\begin{lem} \label{lem:lower-tail-1}
There exist $c_1 > 0$ and $n_0 \in \mathbb{N}$ such that 
for all $n \geq n_0$, $y \in [0, 2\sqrt{n}]$, $t \geq n$, and $v \in T_n$,
\begin{equation} \label{eq:lem-lower-tail-1-statement}
\widetilde{P}_{\rho} (A_v^n (t)) \geq c_1 b^{-n} (1+y) e^{- 2 \sqrt{\log b}~y}.
\end{equation} 
\end{lem}
{\it Proof.}
Fix any $t \geq n$.
By Lemma \ref{lem:local-time-bessel} and
(\ref{eq:0-1-dim-bessel}),
$\widetilde{P}_{\rho} (A_v^n (t))$
is bounded from below by
\begin{align} \label{eq:lem-lower-tail-1-1}
&c_1 \sqrt{\frac{\sqrt{t}}{\sqrt{t} + a_n (t) + y + 1}} 
P_0^B \left( 
\begin{minipage}{180pt}
$-(1-\delta) \sqrt{t} \leq B_s < \frac{a_n (t)}{n} s + y + 1, 
~\forall s \in [0, n],$ \\
$B_n \in [a_n (t) + y, a_n (t) + y + 1)$
\end{minipage}
 \right) \notag \\
&\geq c_1 \sqrt{\frac{\sqrt{t}}{\sqrt{t} + a_n (t) + y + 1}} ~(J_1 - J_2),
\end{align}
where 
we set
\begin{equation*}
J_1 :=P_0^B \Biggl(B_s < \frac{a_n (t)}{n} s + y + 1,~
\forall s \in [0, n], ~
B_n \in [a_n (t) + y, a_n (t) + y + 1) \Biggr),
\end{equation*}
\begin{equation*}
J_2 := P_0^B \Biggl(B_s < - (1-\delta) \sqrt{t},~
\exists s \in [0, n], ~
B_n \in [a_n (t) + y, a_n (t) + y + 1) \Biggr).
\end{equation*}
We first obtain an upper bound of $J_2$.
By using the density
\begin{equation} \label{eq:joint-max-density}
P_0^B \left(B_s \in dx,~\max_{r \in [0, s]} B_r \in dz \right)
= \frac{4(2z - x)}{\sqrt{\pi s^3}} e^{- \frac{(2z-x)^2}{s}} dx dz,
~~s > 0, ~x \leq z,~z \geq 0,
\end{equation}
(see, for example, \cite[Chapter 2, Proposition 8.1]{KS}),
for all $n \in \mathbb{N}$ and $y \in [0, 2\sqrt{n}]$, we have
\begin{align} \label{eq:lem-lower-tail-1-2}
J_2 &= P_0^B \left(\max_{0 \leq s \leq n} B_s > (1-\delta) \sqrt{t},
~B_n \in (- a_n (t) - y - 1, - a_n (t) - y] \right) \notag \\
&\leq c_2 b^{- n} n e^{- c_3 \sqrt{t}} \sqrt{\frac{\sqrt{t} + n}{\sqrt{t}}} e^{- 2 \sqrt{\log b}~y},
\end{align}
where we have used the symmetry of $B$ in the first equality.

Next, we obtain a lower bound of $J_1$.
Recall the process $\widetilde{B}$ from (\ref{eq:bm-drift}).
By the change of measure (\ref{eq:change-measure}) with $j=n$, we have
\begin{equation} \label{eq:lem-lower-tail-1-3}
J_1 = \widetilde{E}_n^B 
\left[ e^{- \frac{2 a_n (t)}{n} \widetilde{B}_n - \frac{(a_n (t))^2}{n}}
1_{\{\widetilde{B}_s < y+1,
~\forall s \in [0, n],~\widetilde{B}_n \in [y, ~y+1) \}} \right].
\end{equation}
By the reflection principle (\ref{eq:density-barrier-standard}),
for all $n \geq n_0$ ($n_0$ is sufficiently large) and $y \in [0, 2\sqrt{n}]$,
(\ref{eq:lem-lower-tail-1-3}) is bounded from below by
\begin{align} \label{eq:lem-lower-tail-1-3-2}
&c_4 b^{- n} \sqrt{\frac{\sqrt{t} + n}{\sqrt{t}}} (y+1) e^{- 2 \sqrt{\log b}~y},
\end{align}
where we used
$e^{- \frac{(y+1-z)^2}{n}} - e^{- \frac{(y+1 + z)^2}{n}} \geq c_5 \frac{y+1}{n}$
for all $z \in [1/2, 1]$.
Thus, by (\ref{eq:lem-lower-tail-1-1}), (\ref{eq:lem-lower-tail-1-2}), and (\ref{eq:lem-lower-tail-1-3-2}),
we have (\ref{eq:lem-lower-tail-1-statement}).~~~$\Box$ \\\\
To obtain upper bounds of $\widetilde{P}_{\rho} \left(A_u^n (t) \cap A_v^n (t) \right), u, v \in T_n$,
we need the following:
\begin{lem} \label{lem:lower-tail-2}
(1) There exists $c_1 > 0$ such that
for all $n \in \mathbb{N}$, $t > 0$, $v \in T_n$, $0 \leq \ell \leq n-1$,
$s < (\sqrt{t} + \frac{a_n (t)}{n} \ell + y + 1)^2$, and $y \geq 0$,
\begin{align} \label{eq:lem-lower-tail-2-statement}
&\widetilde{P}_{v_{\ell}}
\left(
\begin{minipage}{200pt}
$\sqrt{\widetilde{L}_{\widetilde{\tau}^{\downarrow} (s)}^{\downarrow} (v_r)}
< \sqrt{t} + \frac{a_n (t)}{n} r + y + 1, ~\forall r \in [\ell, ~n],$ \\
$\sqrt{\widetilde{L}_{\widetilde{\tau}^{\downarrow} (s)}^{\downarrow} (v)}
\in [\sqrt{t} + a_n (t) + y,~\sqrt{t} + a_n (t) + y + 1)$
\end{minipage}
\right) \notag \\
&\leq c_1 (n-\ell)^{- 3/2} \sqrt{\frac{\sqrt{s}}{\sqrt{t} + a_n (t) + y}}
\left(\sqrt{t} - \sqrt{s} + \frac{a_n (t)}{n} \ell + y + 1 \right)
e^{- \frac{(\sqrt{t} - \sqrt{s} + a_n (t) + y)^2}{n-\ell}},
\end{align}
where
$\left\{\widetilde{L}_r^{\downarrow} (x) : 
(r, x) \in [0, \infty) \times \widetilde{T}_{\leq n - \ell}^{v_{\ell}} \right\}$
is a local time of a Brownian motion 
on $\widetilde{T}_{\leq n - \ell}^{v_{\ell}}$
and
$\widetilde{\tau}^{\downarrow} (s) 
:= \inf \{r \geq 0 : \widetilde{L}_r^{\downarrow} (v_{\ell}) > s \}$. \\
(2) 
There exists $c_3 > 0$ such that for all $n \in \mathbb{N}$, $t > 0$, $v \in T_n$, and $y \geq 0$,
\begin{equation} \label{eq:lem-lower-tail-2-statement-2}
\widetilde{P}_{\rho} (A_v^n (t)) \leq c_3 b^{- n} (y+1) e^{- 2 \sqrt{\log b}~y}.
\end{equation}
\end{lem}
{\it Proof.} We first prove (1).
Recall the process 
$\widetilde{B}$ from (\ref{eq:bm-drift}).
By Lemma \ref{lem:local-time-bessel},
(\ref{eq:0-1-dim-bessel}), and the change of measure (\ref{eq:change-measure}),
the left of (\ref{eq:lem-lower-tail-2-statement}) is bounded from above by
\begin{equation} \label{eq:lem-lower-tail-2-1}
\sqrt{\frac{\sqrt{s}}{\sqrt{t} + a_n (t) + y}} \widetilde{E}_{n-\ell}^B
\left[
\begin{minipage}{205pt}
$e^{- \frac{2 a_n (t)}{n} \widetilde{B}_{n-\ell} - \frac{(a_n (t))^2}{n^2} (n-\ell) } 
1_{\left\{\widetilde{B}_{n-\ell}
\in \sqrt{t} - \sqrt{s} + \frac{a_n (t)}{n} \ell + y + [0, 1) \right\}}$ \\
$\times 1_{\left\{\widetilde{B}_r < \sqrt{t} - \sqrt{s} + \frac{a_n (t)}{n} \ell + y + 1,
~\forall r \in [0, n - \ell] \right \}}$
\end{minipage}
\right].
\end{equation}
By the reflection principle (\ref{eq:density-barrier-standard}),
(\ref{eq:lem-lower-tail-2-1}) is bounded from above by
\begin{equation*}
c_1 \sqrt{\frac{\sqrt{s}}{\sqrt{t} + a_n (t) + y}} 
(n-\ell)^{- 3/2} \left(\sqrt{t} - \sqrt{s} + \frac{a_n (t)}{n} \ell + y + 1 \right)
e^{- \frac{(\sqrt{t} - \sqrt{s} + a_n (t) + y)^2}{n-\ell}},
\end{equation*}
where we have used 
the inequality $1 - e^{-x} \leq x$ for each $x \geq 0$.
Thus, we have obtained (\ref{eq:lem-lower-tail-2-statement}).
The inequality (\ref{eq:lem-lower-tail-2-statement-2})
immediately follows from (1) with $s = t$ and $\ell = 0$.
~~~$\Box$ \\\\
{\it Proof of Proposition \ref{prop:tail}(ii).}
Fix any $n \geq n_0$, $t \geq n$ and $y \in [0, 2\sqrt{n}]$,
where we take $n_0 \in \mathbb{N}$ large enough.
Set
\begin{equation*}
Z := \sum_{v \in T_n} 1_{A_v^n (t)}.
\end{equation*}
We have
\begin{equation} \label{eq:prop-lower-tail-1} 
P_{\rho} \left(\max_{v \in T_n} \sqrt{L_{\tau (t)}^n (v)} \geq \sqrt{t} + a_n (t) + y \right)
\geq \widetilde{P}_{\rho} (Z \geq 1) 
\geq \frac{\left(\widetilde{E}_{\rho}[Z] \right)^2}{\widetilde{E}_{\rho} [Z^2]}.
\end{equation} 
By Lemma \ref{lem:lower-tail-1}, we have
\begin{equation} \label{eq:prop-lower-tail-2}
\widetilde{E}_{\rho}[Z] \geq c_1 (1+y) e^{- 2 \sqrt{\log b}~y}.
\end{equation}
The rest of the proof focuses on obtaining an upper bound of $\widetilde{E}_{\rho} [Z^2]$.
We have
\begin{equation} \label{eq:prop-lower-tail-3}
\widetilde{E}_{\rho} [Z^2]
= \widetilde{E}_{\rho}[Z]
+ \sum_{\ell=0}^{n-1} \sum_{\begin{subarray}{c} v, u \in T_n, \\ |v \wedge u| = \ell \end{subarray}}
\widetilde{P}_{\rho} (A_v^n (t) \cap A_u^n (t) ).
\end{equation}
By Lemma \ref{lem:lower-tail-2} (2), we have
\begin{equation} \label{eq:1st-moment}
\widetilde{E}_{\rho} [Z] \leq c_2 (y+1) e^{- 2 \sqrt{\log b}~y}.
\end{equation}
Fix $1 \leq \ell \leq n-1$ and $v, u \in T_n$ with $|v \wedge u| = \ell$.
Let $\left\{\widetilde{L}_s^{\downarrow} (x) : 
(s, x) \in [0, \infty) \times \widetilde{T}_{\leq n - \ell}^{v_{\ell}} \right \}$
be a local time of a Brownian motion on $\widetilde{T}_{\leq n - \ell}^{v_{\ell}}$. 
Set $\widetilde{\tau}^{\downarrow} (s) := \inf \{r \geq 0 : \widetilde{L}_r^{\downarrow} (v_{\ell}) > s \}.$
For $w \in \{v, u \}$ and $s \geq 0$, we define the event $C_w^{\downarrow} (s)$ by
\begin{equation*}
\left \{
\begin{minipage}{220pt}
$\delta \sqrt{t} \leq 
\sqrt{\widetilde{L}_{\widetilde{\tau}^{\downarrow} (s)}^{\downarrow} (w_r)}
< \sqrt{t} + \frac{a_n (t)}{n} r + y + 1,~\forall r \in [\ell, ~n],$ \\
$\sqrt{\widetilde{L}_{\widetilde{\tau}^{\downarrow} (s)}^{\downarrow} (w)}
\in [\sqrt{t} + a_n (t) + y,~\sqrt{t} + a_n (t) + y + 1)$
\end{minipage}
\right \}.
\end{equation*}
By Lemma \ref{lem:markov-local-time}, we have
\begin{align} \label{eq:prop-lower-tail-4}
&~~~\widetilde{P}_{\rho} (A_v^n (t) \cap A_u^n (t)) \notag \\
&= \widetilde{E}_{\rho} \left[
1_{\left\{\delta \sqrt{t} \leq \sqrt{\widetilde{L}_{\widetilde{\tau} (t)}^n (v_s)} 
< \sqrt{t} + \frac{a_n (t)}{n} s + y + 1,~\forall s \in [0,~ \ell] \right\}} 
\widetilde{P}_{v_{\ell}} 
\left( \bigcap_{w \in \{v, u \}}
C_w^{\downarrow} \left(\widetilde{L}_{\widetilde{\tau} (t)}^n (v_{\ell}) \right)
 \right) \right] \notag \\
&= \widetilde{E}_{\rho} \left[
1_{\left\{\delta \sqrt{t} \leq \sqrt{\widetilde{L}_{\widetilde{\tau} (t)}^n (v_s)} 
< \sqrt{t} + \frac{a_n (t)}{n} s + y + 1,~\forall s \in [0,~ \ell] \right\}} 
\prod_{w \in \{v, u \}} \widetilde{P}_{v_{\ell}} 
\left(C_w^{\downarrow} \left(\widetilde{L}_{\widetilde{\tau} (t)}^n (v_{\ell}) \right) \right) \right] 
\notag \\
&\leq \sum_{i=0}^{\lceil (1-\delta) \sqrt{t} + \frac{a_n (t)}{n} \ell + y \rceil}
 \widetilde{E}_{\rho} \left[
\begin{minipage}{160pt}
$1_{\left\{\delta \sqrt{t} \leq \sqrt{\widetilde{L}_{\widetilde{\tau} (t)}^n (v_s)} 
< \sqrt{t} + \frac{a_n (t)}{n} s + y + 1,~\forall s \in [0,~\ell] \right \}}$ \\
$\times 1_{\left\{\sqrt{\widetilde{L}_{\widetilde{\tau} (t)}^n (v_{\ell})}
\in \sqrt{t} + \frac{a_n (t)}{n} \ell + y + 1 + [- i - 1,~-i) \right\}}$ \\
$\times \prod_{w \in \{v, u \}}\widetilde{P}_{v_{\ell}} 
\left(C_w^{\downarrow} \left(\widetilde{L}_{\widetilde{\tau} (t)}^n (v_{\ell}) \right) \right)$
\end{minipage}
\right],
\end{align}
where we have used the independence of $C_v^{\downarrow} (s)$ and $C_u^{\downarrow} (s)$
for each $s \geq 0$.
(The independence follows from that of two types of excursions of a Brownian motion
around $v_{\ell}$ on 
$\widetilde{T}_{\leq n - \ell - 1}^{v_{\ell + 1}} \cup I_{\{v_{\ell}, ~v_{\ell + 1} \}}$
or on
$\widetilde{T}_{\leq n - \ell - 1}^{u_{\ell + 1}} \cup I_{\{v_{\ell}, ~u_{\ell + 1} \}}$.)
Fix $i \leq \left \lceil (1-\delta) \sqrt{t} + \frac{a_n (t)}{n} \ell + y \right \rceil$.
Lemma \ref{lem:lower-tail-2} implies that
under 
$\left\{\sqrt{\widetilde{L}_{\widetilde{\tau} (t)}^n (v_{\ell})}
\in \sqrt{t} + \frac{a_n (t)}{n} \ell + y + 1 + [- i - 1,~-i) \right\}$,
for all $w \in \{v, u \}$,
$\widetilde{P}_{v_{\ell}} 
[C_w^{\downarrow} (\widetilde{L}_{\widetilde{\tau} (t)}^n (v_{\ell})]$
is bounded from above by
\begin{equation} \label{eq:prop-lower-tail-5}
c_3 (n-\ell)^{- \frac{3}{2}} (i+1) 
\sqrt{\frac{\sqrt{t} + \frac{a_n (t)}{n} \ell + y + 1 - i}{\sqrt{t} + a_n (t) + y}}
~e^{- \frac{\left(\frac{a_n (t)}{n} (n-\ell) + i - 1 \right)^2}{n-\ell}}.
\end{equation} 
By almost the same argument as the proof of Lemma \ref{lem:lower-tail-2} (1),
we have
\begin{align} \label{eq:prop-lower-tail-6}
&\widetilde{P}_{\rho} \left(
\begin{minipage}{200pt}
$\sqrt{\widetilde{L}_{\widetilde{\tau} (t)}^n (v_s)} 
< \sqrt{t} + \frac{a_n (t)}{n} s + y + 1,~\forall s \in [0, \ell],$ \\
$\sqrt{\widetilde{L}_{\widetilde{\tau} (t)}^n (v_{\ell})}
\in \sqrt{t} + \frac{a_n (t)}{n} \ell + y + 1 + [- i - 1,~-i)$
\end{minipage}
\right) \notag \\
&\leq c_4 \sqrt{\frac{\sqrt{t}}{\sqrt{t} + \frac{a_n (t)}{n}\ell + y - i}}
b^{-\ell} \ell^{- \frac{3}{2}} (i+1) (y+1) e^{- 2 \sqrt{\log b}~y}
e^{\frac{2 a_n (t)}{n} i} e^{\frac{3 \log n}{2n} \ell} 
e^{\frac{\log \left(\frac{\sqrt{t} + n}{\sqrt{t}} \right)}{2n} \ell}.
\end{align}
By (\ref{eq:prop-lower-tail-5}) and (\ref{eq:prop-lower-tail-6}),
the right of (\ref{eq:prop-lower-tail-4}) is bounded from above by
\begin{equation} \label{eq:prop-lower-tail-7}
c_5 b^{- 2n + \ell} \ell^{- \frac{3}{2}} (n-\ell)^{-3} n^3 (y+1) e^{- 2 \sqrt{\log b}~y}
e^{- \frac{3 \log n}{2n} \ell} 
e^{- \frac{\log \left(\frac{\sqrt{t} + n}{\sqrt{t}} \right)}{2n} \ell}
\sqrt{\frac{\sqrt{t} + \frac{a_n (t)}{n}\ell + y}{\sqrt{t}}}.
\end{equation}
By (\ref{eq:prop-lower-tail-7}), we have
\begin{equation} \label{eq:prop-lower-tail-8}
\sum_{\ell=1}^{n-1} \sum_{\begin{subarray}{c} v, u \in T_n, \\ |v \wedge u| = \ell \end{subarray}}
\widetilde{P}_{\rho} (A_v^n (t) \cap A_u^n (t))
\leq c_6 (y+1) e^{- 2 \sqrt{\log b}~y}.
\end{equation}
In the case $\ell = 0$, by Lemma \ref{lem:lower-tail-2} (2), we have
\begin{equation} \label{eq:prop-lower-tail-9}
\sum_{\begin{subarray}{c} v, u \in T_n, \\ |v \wedge u| = 0 \end{subarray}}
\widetilde{P}_{\rho} (A_v^n (t) \cap A_u^n (t))
= \sum_{\begin{subarray}{c} v, u \in T_n, \\ |v \wedge u| = 0 \end{subarray}}
\widetilde{P}_{\rho} (A_v^n (t)) \widetilde{P}_{\rho} (A_u^n (t)) 
\leq c_7 (y+1) e^{- 2 \sqrt{\log b}~y},
\end{equation}
where we have used 
the independence of $A_v^n (t)$ and $A_u^n (t)$ for each $v, u \in T_n$ with $|v \wedge u| = 0$
in the first equality. (The independence follows from that of two types of excursions
of a Brownian motion on $\widetilde{T}_{\leq n}$
around $\rho$ restricted to $\widetilde{T}_{\leq n-1}^{v_1} \cup I_{\{\rho, v_1\}}$
or to $\widetilde{T}_{\leq n-1}^{u_1} \cup I_{\{\rho, u_1\}}$).
Thus, by (\ref{eq:prop-lower-tail-1})-(\ref{eq:1st-moment}) 
and (\ref{eq:prop-lower-tail-8})-(\ref{eq:prop-lower-tail-9}), 
we have (\ref{eq:prop-lower-tail-statement}).~~~$\Box$

\section{Geometry of near maxima} \label{sec:geometry}
In this section, we will prove that two leaves with local times near maxima
are either very close or far away.
More specifically, the following is the aim of this section.
\begin{prop} \label{prop:geometry}
There exist $c_1, c_2 \in (0, \infty)$, $n_0, r_0 \in \mathbb{N}$, and $t_0 > 0$ such that
for all $n \geq n_0$, $t \geq t_0$, and $r_0 \leq r \leq n/4$,
\begin{equation} \label{eq:geometry-statement}
P_{\rho} \left(
\begin{minipage}{180pt}
$\exists v, u \in T_n
~\text{with}~r \leq |v \wedge u| \leq n-r :$ \\
$\sqrt{L_{\tau (t)}^n (v)},
~\sqrt{L_{\tau (t)}^n (u)} \geq \sqrt{t} + a_n (t) - c_1 \log r$
\end{minipage}
\right) 
\leq c_2 r^{- 1/8}.
\end{equation}
\end{prop}
\begin{rem} \label{rem:geometry}
Results similar to Proposition \ref{prop:geometry} are known 
for the BBM \cite{ABK1} and 
the two-dimensional DGFF \cite{DZ}.
\end{rem}  
For $n \in \mathbb{N}$, $t > 0$, and $k \in \mathbb{Z}$, set
\begin{equation*} 
\Gamma_k^n (t) := \left\{v \in T_n : \sqrt{L_{\tau (t)}^n (v)} 
\in [\sqrt{t} + a_n (t) - k - 1, \sqrt{t} + a_n (t) - k] \right \}.
\end{equation*}
\begin{rem} \label{rem:geometry-2}
Fix $n^{\prime} > n \geq 1$, $t > 0$, and $k \in \mathbb{Z}$.
Set
\begin{equation*}
\Gamma_k^{n^{\prime}, n} (t) := \left\{v \in T_n : \sqrt{L_{\tau (t)}^{n^{\prime}} (v)} 
\in [\sqrt{t} + a_n (t) - k - 1, \sqrt{t} + a_n (t) - k] \right \}.
\end{equation*}
Since the law of the simple random walk on $T_{\leq n^{\prime}}$
watched only on $T_{\leq n}$ is the same as 
that of the simple random walk on $T_{\leq n}$,
we have 
$|\Gamma_k^{n^{\prime}, n} (t)| \stackrel{d}{=} |\Gamma_k^n (t)|$. 
\end{rem}
In the proof of Proposition \ref{prop:geometry}, we will use the following repeatedly.
\begin{lem} \label{lem:level-set-tail}
(i) There exist $c_1 > 0$ and $t_0 > 0$ such that 
for all $n \in \mathbb{N}$, $t \geq t_0$, $k \leq -1$,
and $\lambda \in \mathbb{R}$ with $k + \lambda \geq 0$,
\begin{align} \label{eq:lem-level-set-tail-2-statement}
&~~~P_{\rho} \left(|\Gamma_k^n (t)| \geq e^{2 \sqrt{\log b}~(k+\lambda)} \right) \notag\\
&\leq c_1 (\lambda + 1)
e^{- 2 \sqrt{\log b}~\left(\lambda - \frac{3}{4 \sqrt{\log b}} \log (\lceil (k+\lambda)^2 \rceil \vee 1)
- \frac{1}{4 \sqrt{\log b}} \log
\left(\frac{\sqrt{t} + a_n (t) - k - 1 +\lceil (k+\lambda)^2 \rceil \vee 1 }{\sqrt{t} + a_n (t) - k - 1} \right)
 \right)}.
\end{align}
(ii) There exist $c_2 > 0$ and $t_0 > 0$ such that 
for all $n \in \mathbb{N}$, $t \geq t_0$, $\lambda > 0$,
and $k \geq 0$ with $\sqrt{t} + a_n (t) - k - 1 \geq c_2$,
(\ref{eq:lem-level-set-tail-2-statement}) holds.
\end{lem}
{\it Proof.}
(i) Fix $n \in \mathbb{N}$, $t > 0$, $k \leq -1$, and $\lambda \in \mathbb{R}$ 
with $k + \lambda \geq 0$.
Set $r := \lceil (k + \lambda)^2 \rceil \vee 1$ and
$y := \lambda - \frac{3}{4 \sqrt{\log b}} \log r 
- \frac{1}{4 \sqrt{\log b}} \log \left(\frac{\sqrt{t} + a_n (t) - k - 1 + r}{\sqrt{t} + a_n (t) - k - 1} \right)$.
If $y < 0$, then it is clear that (\ref{eq:lem-level-set-tail-2-statement}) holds
because
$P_{\rho} \left(|\Gamma_k^n (t)| \geq e^{2 \sqrt{\log b}~(k+\lambda)} \right)
\leq 1 \leq (\lambda + 1) e^{- 2 \sqrt{\log b}~y}$.
So, we may assume that $y \geq 0$.
Fix any $K > 0$. We have
\begin{align} \label{eq:lem-level-set-tail-2-1}
&~~~P_{\rho} \left(\max_{v \in T_{n + r}} \sqrt{L_{\tau (t)}^{n+r} (v)} \geq \sqrt{t} + a_{n+r} (t) + y \right) 
\notag \\
&\geq P_{\rho} \left(|\Gamma_k^{n+r, n} (t)| \geq K,~
\max_{v \in T_{n + r}} \sqrt{L_{\tau (t)}^{n+r} (v)} \geq \sqrt{t} + a_{n+r} (t) + y \right) \notag \\
&= P_{\rho} \left(|\Gamma_k^{n+r, n} (t)| \geq K \right)
- P_{\rho} \left(|\Gamma_k^{n+r, n} (t)| \geq K,
\max_{v \in T_{n + r}} \sqrt{L_{\tau (t)}^{n+r} (v)} < \sqrt{t} + a_{n+r} (t) + y \right).
\end{align}
We estimate the second term on the right-hand side of (\ref{eq:lem-level-set-tail-2-1}).
By Lemma \ref{lem:markov-local-time}, we have
\begin{align} \label{eq:lem-level-set-tail-2-2}
&~~~P_{\rho} \left(|\Gamma_k^{n+r, n} (t)| \geq K,~
\max_{v \in T_{n + r}} \sqrt{L_{\tau (t)}^{n+r} (v)} < \sqrt{t} + a_{n+r} (t) + y \right) 
\notag \\
&\leq \sum_{\begin{subarray}{c} S \subset T_n, \\ |S| \geq K \end{subarray}}
P_{\rho} \left(\Gamma_k^{n+r, n} (t) = S,
~\max_{v \in T_r^u} \sqrt{L_{\tau (t)}^{n+r} (v)} < \sqrt{t} + a_{n+r} (t) + y,
~\forall u \in S \right) \notag \\
&= \sum_{\begin{subarray}{c} S \subset T_n, \\ |S| \geq K \end{subarray}}
E_{\rho} \Biggl[ 1_{\left\{\Gamma_k^{n+r, n} (t) = S \right\}}
\prod_{u \in S} \widetilde{P}_u \left(\max_{v \in T_r^u} 
\sqrt{\widetilde{L}_{\widetilde{\tau}^{\downarrow} 
\left(L_{\tau (t)}^{n+r} (u)\right)}^{\downarrow} (v)} 
< \sqrt{t} + a_{n+r} (t) + y \right)
\Biggr],
\end{align}
where for each $u \in T_n$, 
$\left\{\widetilde{L}_s^{\downarrow} (x) : (s, x) \in [0,~\infty) \times \widetilde{T}_{\leq r}^u \right\}$
is a local time of a Brownian motion on $\widetilde{T}_{\leq r}^u$
and $\widetilde{\tau}^{\downarrow} (q) := \inf \{s \geq 0: \widetilde{L}_s^{\downarrow} (u) > q \}$.
We omit the subscript $u$ in $\widetilde{L}_s^{\downarrow} (x)$
and $\widetilde{\tau}^{\downarrow}$ to simplify the notation.

We estimate each probability on the right-hand side of (\ref{eq:lem-level-set-tail-2-2}).
Fix $S \subset T_n$ with $|S| \geq K$ and $u \in S$.
Note that under the event that $\Gamma_k^{n+r, n} (t) = S$, we have
\begin{equation*}
\sqrt{t} + a_{n+r} (t) + y 
\leq \sqrt{L_{\tau (t)}^{n+r} (u)} + a_r \left(L_{\tau (t)}^{n+r} (u) \right)
+ k + \lambda + 1.
\end{equation*}
By this and
Proposition \ref{prop:tail-max-revisit}
for $t \geq t_0$,
where $t_0$ is sufficiently large, we have
\begin{align} \label{eq:lem-level-set-tail-2-3}
&~~~\widetilde{P}_u \left(\max_{v \in T_r^u} 
\sqrt{\widetilde{L}_{\widetilde{\tau}^{\downarrow} \left(L_{\tau (t)}^{n+r} (u)\right)}^{\downarrow} (v)} 
< \sqrt{t} + a_{n+r} (t) + y \right) \notag \\
&\leq 1 - \widetilde{P}_u \left(\max_{v \in T_r^u} 
\sqrt{\widetilde{L}_{\widetilde{\tau}^{\downarrow} \left(L_{\tau (t)}^{n+r} (u)\right)}^{\downarrow} (v)} 
\geq \sqrt{L_{\tau (t)}^{n+r} (u)} + a_r \left(L_{\tau (t)}^{n+r} (u) \right) + k + \lambda + 1 \right) 
\notag \\
&\leq 1 - c_1 e^{- 2 \sqrt{\log b} ~(k + \lambda + 1)}.
\end{align}
By (\ref{eq:lem-level-set-tail-2-1})-(\ref{eq:lem-level-set-tail-2-3}), we have
\begin{align} \label{eq:lem-level-set-tail-2-4}
&~~~~P_{\rho} \left(\max_{v \in T_{n + r}} \sqrt{L_{\tau (t)}^{n+r} (v)} \geq \sqrt{t} + a_{n+r} (t) + y \right) 
\notag \\
&\geq \left(1 - \exp \left \{- c_1 K e^{- 2 \sqrt{\log b} (k + \lambda + 1)} \right\} \right)
P_{\rho} \left(|\Gamma_k^{n+r, n} (t) | \geq K \right).
\end{align}
By Remark \ref{rem:geometry-2},
(\ref{eq:lem-level-set-tail-2-4}) with $K := e^{2 \sqrt{\log b} (k + \lambda)}$,
and Proposition \ref{prop:tail}(i), 
we have (\ref{eq:lem-level-set-tail-2-statement}).

(ii) The proof of (ii) is almost the same as that of (i), so we omit the detail.~~~$\Box$ \\

For the rest of this section, we focus on proving the following.
\begin{lem} \label{lem:geometry-key} 
Fix $0 < \underline{c} < \overline{c} < \frac{3}{4 \sqrt{\log b}}$. 
There exist $c_1 > 0$, $n_0, s_0 \in \mathbb{N}$, and $t_0 > 0$ such that 
for all $n \geq n_0$, $t \geq t_0$, and $s_0 \leq s \leq n - s_0$,
\begin{align} \label{eq:lem-geometry-key-statement}
&~~~P_{\rho} \left[
\begin{minipage}{240pt}
$\exists v, u \in T_n~\text{with}~|v \wedge u| = s :$\\
$\sqrt{L_{\tau (t)}^n (v)}, \sqrt{L_{\tau (t)}^n (u)}
\geq \sqrt{t} + a_n (t) - (\overline{c} - \underline{c}) \log (s \wedge (n-s))$
\end{minipage}
\right] \notag \\
&\leq c_1 (\log (s \wedge (n-s)))^8
(s \wedge (n-s))^{-3 + 4 \overline{c} \sqrt{\log b} - 2 \underline{c} \sqrt{\log b}}
\notag \\
&~~~+ c_1 (\log (s \wedge (n-s)))^6 (s \wedge (n-s))^{- 2 \underline{c} \sqrt{\log b}}.
\end{align}
\end{lem}
Before we prove this, let us show that 
Lemma \ref{lem:geometry-key} implies Proposition \ref{prop:geometry}.\\\\
{\it Proof of Proposition \ref{prop:geometry} via Lemma \ref{lem:geometry-key}.}
Fix any $n \geq n_0$, $t \geq t_0$, and
$r_0 \leq r \leq n/4$, where we take $n_0, r_0 \in \mathbb{N}$
and $t_0 > 0$ sufficiently large. 
By Lemma \ref{lem:geometry-key} with
$\underline{c} = \frac{5}{8 \sqrt{\log b}}$ and $\overline{c} = \frac{11}{16 \sqrt{\log b}}$,
the left of (\ref{eq:geometry-statement}) is bounded from above by
\begin{equation*} 
c_{1}
\sum_{s=r}^{n-r}
\left\{ (\log (s \wedge (n-s)))^8 (s \wedge (n-s))^{- 3/2} 
+ (\log (s \wedge (n-s)))^6 (s \wedge (n-s))^{- 5/4} \right \}.
\end{equation*}
This is bounded from above by $c_2 r^{- 1/8}$.~~$\Box$
\\\\
{\it Proof of Lemma \ref{lem:geometry-key}.}
Fix any $n \geq n_0$, $t \geq t_0$, $s_0 \leq s \leq n - s_0$,
where we take $n_0, s_0 \in \mathbb{N}$
and $t_0 > 0$ sufficiently large.
Set $z := \underline{ c} \log (s \wedge (n-s))$.
The left of (\ref{eq:lem-geometry-key-statement})
is bounded from above by
\begin{equation} \label{eq:lem-geometry-key-1}
P_{\rho} \left[
\begin{minipage}{240pt}
$\exists k \in \mathbb{Z}~\text{with}~\sqrt{t} + a_{s} (t) - k > 0,
~\exists j \geq -k$~\text{s.t.} \\
$|\Gamma_k^{n, s} (t)| \in \left[e^{2 \sqrt{\log b} (k+j)},~e^{2 \sqrt{\log b} (k+j+1)} \right),$ \\
$\exists w \in \Gamma_k^{n, s} (t),
~\exists w_1, w_2 \in T_1^w~\text{with}~w_1 \neq w_2~\text{s.t.}~ \forall i \in \{1, 2 \}$, \\
$\max_{v \in T_{n-s-1}^{w_i}} \sqrt{L_{\tau (t)}^n (v)} \geq \sqrt{t} + a_n (t) 
- \overline{c} \log (s \wedge (n-s) )
+ z$
\end{minipage}
\right].
\end{equation}
For $k \in \mathbb{Z}$, we set
$j^* (k, z) := \lceil \max \{|k|,~z \} \rceil$.
Fix sufficiently large constant $t_* > 0$. 
We will decompose (\ref{eq:lem-geometry-key-1}) into three terms
with respect to $k$ and $j$:
(i) $k \leq \sqrt{t} + a_s (t) - t_* - 1$ and $-k \leq j \leq j^{*} (k, z)$;
(ii) $k \leq \sqrt{t} + a_s (t) - t_* - 1$ and $j > j^{*} (k, z)$;
(iii) $k \geq \sqrt{t} + a_s (t) - t_* - 1$.
Then, (\ref{eq:lem-geometry-key-1}) is bounded from above by
\begin{align} \label{eq:lem-geometry-key-1-2}
&\sum_{\begin{subarray}{c} k \in \mathbb{Z}, \\ \sqrt{t} + a_s (t) - k - 1 \geq t_* \end{subarray}}
\sum_{j= -k}^{j^* (k, z)}
P_{\rho} \left[
\begin{minipage}{230pt}
$|\Gamma_k^{n, s} (t)| \in \left[e^{2 \sqrt{\log b} (k+j)},~e^{2 \sqrt{\log b} (k+j+1)} \right),$ \\
$\exists w \in \Gamma_k^{n, s} (t),
\exists w_1, w_2 \in T_1^w~\text{with}~w_1 \neq w_2~\text{s.t.} \forall i \in \{1, 2 \}$, \\
$\max_{v \in T_{n-s-1}^{w_i}} \sqrt{L_{\tau (t)}^n (v)} \geq \sqrt{t} + a_n (t) 
- \overline{c} \log (s \wedge (n-s) )
+ z$
\end{minipage}
\right] \notag \\
&~~~
+ \sum_{\begin{subarray}{c} k \in \mathbb{Z}, \\ \sqrt{t} + a_s (t) - k - 1 \geq t_* \end{subarray}}
P_{\rho} \left(|\Gamma_k^{n, s} (t)| \geq e^{2 \sqrt{\log b} (k+j^* (k, z))} \right) \notag \\
&~~~
+ \sum_{\begin{subarray}{c} k \in \mathbb{Z}, 
\\ 0 < \sqrt{t} + a_s (t) - k  \leq t_* + 1 \end{subarray}}
P_{\rho} \left[
\begin{minipage}{230pt}
$\exists w \in T_s~\text{with}~\sqrt{L_{\tau (t)}^n (w)} \leq t_* + 1$, \\
$\exists w_1, w_2 \in T_1^w~\text{with}~w_1 \neq w_2~\text{s.t.}~\forall i \in \{1, 2 \}$, \\
$\max_{v \in T_{n-s-1}^{w_i}} \sqrt{L_{\tau (t)}^n (v)} \geq \sqrt{t} + a_n (t) 
- \overline{c} \log (s \wedge (n-s)) + z$
\end{minipage}
\right] \notag \\
&=: \sum_{\begin{subarray}{c} k \in \mathbb{Z}, \\ \sqrt{t} + a_s (t) - k - 1 \geq t_* \end{subarray}}
\sum_{j= -k}^{j^* (k, z)} J_1 (k, j) 
+ J_2 
+ \sum_{\begin{subarray}{c} k \in \mathbb{Z}, 
\\ 0 < \sqrt{t} + a_s (t) - k  \leq t_* + 1 \end{subarray}} J_3 (k).
\end{align} 
For each $w \in T_s$,
let $\widetilde{L}^{\downarrow}$ be a local time of a Brownian motion on 
$\widetilde{T}_{\leq n-s}^w$ and set
$\widetilde{\tau}^{\downarrow} (q) := \inf \{p \geq 0 : \widetilde{L}_p^{\downarrow} (w) > q \}$.
We omit the subscript $w$ in $\widetilde{L}^{\downarrow}$ and $\widetilde{\tau}^{\downarrow} (q)$.
Fix $k \leq \sqrt{t} + a_s (t) - t_* - 1$ and $- k \leq j \leq j^{*} (k, z)$.
Set $I_{k, j} := [e^{2 \sqrt{\log b} (k+j)},~e^{2 \sqrt{\log b} (k+j+1)} )$.
$J_1 (k, j)$ is bounded from above by
\begin{align} \label{eq:lem-geometry-key-2}
&\sum_{\begin{subarray}{c} S \subset T_s, \\ 
|S| \in I_{k, j} \end{subarray}}
\sum_{w \in S}
\sum_{\begin{subarray}{c} w_1, w_2 \in T_1^w, \\ w_1 \neq w_2 \end{subarray}}
J_1^{k, j} (S, w, w_1, w_2),
\end{align}
where $J_1^{k, j} (S, w, w_1, w_2)$ is given by
\begin{equation} \label{eq:lem-geometry-key-2-2}
P_{\rho} \left[
\Gamma_k^{n, s} (t) = S, 
\max_{v \in T_{n-s-1}^{w_i}} \sqrt{L_{\tau (t)}^n (v)} \geq \sqrt{t} + a_n (t) 
- \overline{c} \log (s \wedge (n-s) )
+ z,~\forall i \in \{1, 2 \}
\right].
\end{equation}
Fix $S \subset T_s$ with $|S| \in I_{k, j}$,
$w \in S$, and $w_1, w_2 \in T_1^w$ with $w_1 \neq w_2$.
By Lemma \ref{lem:markov-local-time}
and the independence of two types of
excursions restricted to $\widetilde{T}_{\leq n - s - 1}^{w_1} \cup I_{\{w, w_1 \}}$
or to $\widetilde{T}_{\leq n - s - 1}^{w_2} \cup I_{\{w, w_2 \}}$, 
(\ref{eq:lem-geometry-key-2-2}) is equal to
\begin{equation} \label{eq:lem-geometry-key-2-3}
E_{\rho} \left[
\begin{minipage}{307pt}
$1_{\{\Gamma_k^{n, s} (t) = S \}}$ \\
$\times \displaystyle \prod_{i \in \{1,~2 \}} \widetilde{P}_w \left(\max_{v \in T_{n-s-1}^{w_i}} 
\sqrt{\widetilde{L}_{\widetilde{\tau}^{\downarrow} 
\left(L_{\tau (t)}^n (w) \right)}^{\downarrow} (v)} \geq \sqrt{t} + a_n (t) 
- \overline{c} \log (s \wedge (n-s) )+ z \right)$
\end{minipage}
\right].
\end{equation}
By the symmetry of the $b$-ary tree and (\ref{eq:conti-disc-local-time}),
(\ref{eq:lem-geometry-key-2-3})
is bounded from above by
\begin{align} \label{eq:lem-geometry-key-2-4}
&P_{\rho} \left[\Gamma_k^{n, s} (t) = S \right] \notag \\
&\times
P_{\rho} \left(\max_{v \in T_{n-s}} 
\sqrt{L_{\tau \left((\sqrt{t} + a_s (t) - k)^2 \right)}^{n-s} (v)} 
\geq \sqrt{t} + a_n (t) 
- \overline{c} \log (s \wedge (n-s) )+ z \right)^2.
\end{align}
By (\ref{eq:lem-geometry-key-2-4}), 
(\ref{eq:lem-geometry-key-2}) is bounded from above by
\begin{align} \label{eq:lem-geometry-key-2-5}
&c_1
e^{2 \sqrt{\log b} (k+j)}
P_{\rho} \left(|\Gamma_k^{n, s} (t) | \geq e^{2 \sqrt{\log b} (k+j)} \right) \notag \\
&\times P_{\rho} \left(\max_{v \in T_{n-s}} 
\sqrt{L_{\tau \left((\sqrt{t} + a_s (t) - k)^2 \right)}^{n-s} (v)} 
\geq \sqrt{t} + a_n (t) 
- \overline{c} \log (s \wedge (n-s) )+ z \right)^2.
\end{align}
We estimate (\ref{eq:lem-geometry-key-2-5})
in different ways according to three cases:
(a) $k \geq - z$ and $j \geq 0$; (b) $k < -z$; (c) $j < 0$.
In the case (a), we use Proposition \ref{prop:tail}(i)
, Lemma \ref{lem:level-set-tail}, and Remark \ref{rem:geometry-2}.
In the case (b), we only use Lemma \ref{lem:level-set-tail} and Remark \ref{rem:geometry-2}
and estimate the square of the probability in (\ref{eq:lem-geometry-key-2-5}) just by $1$.
In the case (c), we only use Proposition \ref{prop:tail}(i) 
and estimate the probability in the first display of (\ref{eq:lem-geometry-key-2-5}) just by $1$.
Note that for $k \in \mathbb{Z}$ with $k\geq -z$, we have 
\begin{align*}
\sqrt{t} + a_n (t) - \overline{c} \log (s \wedge (n-s)) + z
&\geq (\sqrt{t} + a_s (t) - k) + a_{n-s} ((\sqrt{t} + a_s (t) - k)^2) \notag \\
&~~~+ \left(\frac{3}{4 \sqrt{\log b}} - \overline{c} \right) \log (s \wedge (n-s)) + k + z - c_2.
\end{align*}
Recall that $z = \underline{c} \log (s \wedge (n-s))$.
By these observations,
for sufficiently large $t_*$,
the first term of (\ref{eq:lem-geometry-key-1-2}) is bounded from above by
\begin{align} \label{eq:lem-geometry-key-4}
&c_3 (\log (s \wedge (n-s)))^8
(s \wedge (n-s))^{-3 + 4\overline{c} \sqrt{\log b} - 2 \underline{c} \sqrt{\log b}} \notag \\
&+c_3 (\log (s \wedge (n-s)))(s \wedge (n-s))^{- 2 \underline{c} \sqrt{\log b}}.
\end{align}

Next, we estimate the second term of (\ref{eq:lem-geometry-key-1-2}).
By Lemma \ref{lem:level-set-tail} and Remark \ref{rem:geometry-2}, we have
\begin{equation} \label{eq:lem-geometry-key-5}
J_2 
\leq c_4 (\log (s \wedge (n-s)))^6 (s \wedge (n-s))^{- 2 \underline{c} \sqrt{\log b}}.
\end{equation}

Finally, we estimate the third term of (\ref{eq:lem-geometry-key-1-2}).
Fix $k \geq \sqrt{t} + a_s (t) - t_* - 1$.
By Lemma \ref{lem:markov-local-time}, $J_3 (k)$ is bounded from above by
\begin{equation} \label{eq:lem-geometry-key-6}
\sum_{w \in T_s}
\sum_{\begin{subarray}{c} w_1, w_2 \in T_1^w, \\ w_1 \neq w_2 \end{subarray}} 
J_3^k (w, w_1, w_2),
\end{equation}
where $J_3^k (w, w_1, w_2)$ is given by
\begin{align} \label{eq:lem-geometry-key-6-1}
E_{\rho} \left[
\begin{minipage}{300pt}
$1_{\left\{\sqrt{L_{\tau (t)}^n (w)} \leq t_* + 1 \right\}}$ \\
$\times \displaystyle \prod_{i \in \{1, 2 \}}
\widetilde{P}_w \left[\max_{v \in T_{n - s - 1}^{w_i}} 
\sqrt{\widetilde{L}_{\widetilde{\tau}^{\downarrow} \left(L_{\tau (t)}^n (w) \right)}^{\downarrow} (v)} 
\geq \sqrt{t} + a_n (t) - \overline{c} \log (s \wedge (n - s)) + z \right]$
\end{minipage}
\right].
\end{align}
By (\ref{eq:lem-geometry-key-6-1}) 
together with the symmetry of the $b$-ary tree,
(\ref{eq:lem-geometry-key-6}) is bounded from above by
\begin{equation} \label{eq:lem-geometry-key-6-2}
c_5
b^s 
P_{\rho} \biggl(\max_{v \in T_{n-s}} 
\sqrt{L_{\tau \left((t_* + 1)^2 \right)}^{n-s} (v)}
\geq \sqrt{t} + a_n (t) 
- \overline{c} \log (s \wedge (n-s)) + z \biggr)^2.
\end{equation}
Since $t_0$ and $n_0$ are sufficiently large, we have
\begin{align*}
\sqrt{t} + a_n (t) - \overline{c} \log (s \wedge (n-s)) + z
&\geq (t_* + 1) + a_{n-s} ((t_* + 1)^2) \notag \\
&~~~+ \frac{1}{2} \sqrt{\log b}~s - \overline{c} \log (s \wedge (n-s)) + z.
\end{align*}
By Proposition \ref{prop:tail}(i) and (\ref{eq:lem-geometry-key-6-2}), 
the third term of (\ref{eq:lem-geometry-key-1-2}) is bounded from above by
\begin{equation} \label{eq:lem-geometry-key-7}
c_6 s^2 b^{- s} (s \wedge (n-s))^{4 (\overline{c} - \underline{c}) \sqrt{\log b}}.
\end{equation}
Thus, by (\ref{eq:lem-geometry-key-1-2}), (\ref{eq:lem-geometry-key-4}),
(\ref{eq:lem-geometry-key-5}), and (\ref{eq:lem-geometry-key-7}), we have
(\ref{eq:lem-geometry-key-statement}).~~~$\Box$

\section{Limiting tail of the maximum of local times} \label{sec:limit-tail}
The aim of this section is to prove the exact asymptotics of the tail
of the maximum of local times.
Recall the constants $\beta_*$ and $\gamma_*$
from (\ref{eq:const-beta}) and (\ref{eq:const-gamma}).
\begin{prop} \label{prop:limit-tail}
Fix positive sequences $(y_j^{+})_{j \geq 1}, (y_j^{-})_{j \geq 1}$
with $y_j^{-} \leq y_j^{+}$ for each $j \geq 1$ and $\lim_{j \to \infty} y_j^{-} = \infty$.
For each $j \geq 1$, fix sequences $(t_n^{+} (j))_{n \geq 1}, (t_n^{-} (j))_{n \geq 1}$
with $\lim_{n \to \infty} \sqrt{t_n^{+} (j)}/n = \lim_{n \to \infty} \sqrt{t_n^{-} (j)}/n = \theta \in [0, \infty]$
and
$t_n^{-} (j) \leq t_n^{+} (j)$ for each $n \geq 1$.
For all $\varepsilon > 0$, there exists $j_0 \in \mathbb{N}$ such that
the following holds for each $j \geq j_0$:
there exists $n_0 (j) \in \mathbb{N}$ such that for all $n \geq n_0 (j)$,
\begin{equation} \label{eq:prop-limit-tail-statement}
\frac{P_{\rho} \left(\max_{v \in T_n} \sqrt{L_{\tau (t)}^n (v)} > \sqrt{t} + a_n (t) + y_j \right)}
{y_j e^{- 2 \sqrt{\log b}~y_j}}
\in \left[\frac{4}{\sqrt{\pi}} \beta_* \gamma_* - \varepsilon,
~\frac{4}{\sqrt{\pi}} \beta_* \gamma_* + \varepsilon \right],
\end{equation}
uniformly in $y_j$ and $t$ with
\begin{align} \label{eq:prop-limit-tail-assumption-uniformity}
y_j^{-} \leq y_j \leq y_j^{+},~~~
t \geq c_* n \log n~
\text{and}~t_n^{-} (j) \leq t \leq t_n^{+} (j),
\end{align} 
for some constant $c_* > 0$ not depending on $\varepsilon$, $j$, $n$.
\end{prop}  
Fix $\delta \in (0,~1)$ and $\kappa > \frac{1}{2 \sqrt{\log b}}$.
Fix  $t > 0$, $y > 0$, $\ell (y) \in \mathbb{N}$, $n > \ell (y)$.
We will approximate the tail of the maximum of local times by
$\widetilde{E}_{\rho} [\Lambda_{y, \ell (y)}^n (t)]$, where
$$\Lambda_{y, \ell (y)}^n (t) := \sum_{v \in T_{n - \ell (y)}} 1_{F_{v, y, \ell(y)}^n (t)},$$
and for each $v \in T_{n - \ell (y)}$,
\begin{equation*}
F_{v, y, \ell(y)}^n (t) := \left\{
\begin{minipage}{220pt}
$\delta \sqrt{t} \leq 
\sqrt{\widetilde{L}_{\widetilde{\tau} (t)}^n (v_s)}
\leq \sqrt{t} + \frac{a_n (t)}{n} s + y,~\forall s \in [0,~n - \ell(y)]$, \\
$\max_{u \in T_{\ell(y)}^v} \sqrt{\widetilde{L}_{\widetilde{\tau} (t)}^n (u)} > \sqrt{t} + a_n (t) + y$
\end{minipage}
\right \}.
\end{equation*}
To do so, we need an intermediate approximation
$\widetilde{E}_{\rho} [\widetilde{\Lambda}_{y, \ell (y)}^{n} (t)]$,
where
$$\widetilde{\Lambda}_{y, \ell (y)}^{n} (t) := \sum_{v \in T_{n - \ell (y)}} 
1_{\widetilde{F}_{v, y, \ell(y)}^{n} (t)},$$
and for each $v \in T_{n - \ell (y)}$,
\begin{align*}
\widetilde{F}_{v, y, \ell(y)}^{n} (t) := 
\left\{ 
\begin{minipage}{280pt}
$\delta \sqrt{t} \leq
\sqrt{\widetilde{L}_{\widetilde{\tau} (t)}^n (v_s)} \leq
\sqrt{t} + \frac{a_n (t)}{n} s + y + y^{\frac{1}{20}} + \kappa (\log (s \wedge (n-\ell (y)-s)))_+$, \\
$\forall s \in [0,~n - \ell(y)], ~
\max_{u \in T_{\ell(y)}^v} \sqrt{\widetilde{L}_{\widetilde{\tau} (t)}^n (u)} > \sqrt{t} + a_n (t) + y$
\end{minipage}
\right \}.
\end{align*}
In the approximation, the entropic repulsion (Lemma \ref{lem:density-barrier} (ii)) plays an important role.
In Lemma \ref{lem:limit-tail-3},
we show that
$P_{\rho} (\max_{v \in T_n} \sqrt{L_{\tau (t)}^n (v)} \geq \sqrt{t} + a_n (t) + y)
/\widetilde{E}_{\rho} [\Lambda_{y, \ell(y)}^n (t)]$ is close to $1$.
To obtain an upper bound of this, we need the following:
\begin{lem} \label{lem:limit-tail-1}
There exist $c_1, c_2 \in (0, \infty)$, $y_0 > 0$, and $\{\delta_{y^{\prime}} : y^{\prime} > 0 \}$ with
$\lim_{y^{\prime} \to \infty} \delta_{y^{\prime}} = 0$ such that the following holds:
for all $y \geq y_0$ and $\ell (y) > e^{\frac{8\sqrt{\log b}}{3} y^{1/20}}$, 
there exists $n_0 = n_0 (y, \ell (y)) \in \mathbb{N}$ such that
for all $n \geq n_0$ and $t \geq c_1 n \log n$,
\begin{equation} \label{eq:lem-limit-tail-1-statement}
\frac{\widetilde{E}_{\rho} \left[\Lambda_{y, \ell(y)}^n (t) \right]}
{\widetilde{E}_{\rho} \left[\widetilde{\Lambda}_{y, \ell(y)}^n (t) \right]} 
\geq (1-\delta_y)
\left(1 - c_2 (\ell(y))^{- 1/2} - \delta_y \right).
\end{equation}
\end{lem}
{\it Proof.}
Fix any $y \geq y_0$ and $\ell (y) > e^{\frac{8 \sqrt{\log b}}{3} y^{1/20}}$,
where we take $y_0 > 0$ large enough.
Throughout the proof, given $n \in \mathbb{N}$,
we assume that $t \geq c_* n \log n$, where $c_*$ is a sufficiently large positive constant. 
Fix $v \in T_{n - \ell(y)}$.
Let $\widetilde{L}^{\downarrow}$ be a local time of a Brownian motion on
$\widetilde{T}_{\leq \ell(y)}^v$ and set
$\widetilde{\tau}^{\downarrow} (s) := \inf \{r \geq 0 : \widetilde{L}_r^{\downarrow} (v) > s \}$.
Recall the definitions of
$\mu_{n - \ell (y), y}$, $\mu_{n - \ell(y), y}^*$, and $\delta_y$ from Lemma \ref{lem:density-barrier}.
By Lemma \ref{lem:markov-local-time}, we have
\begin{align} \label{eq:lem-limit-tail-1-1}
&~~~\widetilde{P}_{\rho} \left(\widetilde{F}_{v, y, \ell(y)}^n (t) \right)
- \widetilde{P}_{\rho} \left(F_{v, y, \ell(y)}^n (t) \right) \notag \\
&= \widetilde{E}_{\rho} \left[
\begin{minipage}{265pt}
$1_{\left\{\delta \sqrt{t} \leq \sqrt{\widetilde{L}_{\widetilde{\tau} (t)}^n (v_s)} 
\leq \sqrt{t} + \frac{a_n (t)}{n} s + y + y^{\frac{1}{20}} + \kappa (\log (s \wedge (n-\ell (y)-s)))_+,~
\forall s \in [0,~ n - \ell(y)] \right\}}$ \\
$\times \widetilde{P}_v \left(\max_{u \in T_{\ell (y)}^v} 
\sqrt{\widetilde{L}_{\widetilde{\tau}^{\downarrow} \left(L_{\tau (t)}^n (v) \right)}^{\downarrow} (u)}
> \sqrt{t} + a_n (t) + y \right)$
\end{minipage}
\right] \notag \\
&- \widetilde{E}_{\rho} \left[
\begin{minipage}{210pt}
$1_{\left\{\delta \sqrt{t} \leq \sqrt{\widetilde{L}_{\widetilde{\tau} (t)}^n (v_s)} 
\leq \sqrt{t} + \frac{a_n (t)}{n} s + y,~
\forall s \in [0,~n - \ell(y)] \right\}}$  \\
$\times \widetilde{P}_v \left(\max_{u \in T_{\ell (y)}^v} 
\sqrt{\widetilde{L}_{\widetilde{\tau}^{\downarrow} \left(L_{\tau (t)}^n (v) \right)}^{\downarrow} (u)}
> \sqrt{t} + a_n (t) + y \right)$
\end{minipage}
\right].
\end{align}
Recall the process $\widetilde{B}$ from (\ref{eq:bm-drift}).
By Lemma \ref{lem:local-time-bessel}, (\ref{eq:0-1-dim-bessel}), and the change of measure
(\ref{eq:change-measure}),
(\ref{eq:lem-limit-tail-1-1}) is bounded from above by
\begin{align} \label{eq:lem-limit-tail-1-1-2}
&\widetilde{E}_{n - \ell (y)}^B \left[
1_{\left\{- (1 - \delta) \sqrt{t} - \frac{a_n (t)}{n} s
\leq \widetilde{B}_s \leq y + y^{\frac{1}{20}} 
+ \kappa (\log (s \wedge (n-\ell (y)-s)))_+,~
\forall s \in [0,~ n - \ell(y)] \right\}}
\psi (\widetilde{B}_{n - \ell (y)})
\right] \notag \\
&- \widetilde{E}_{n - \ell (y)}^B \left[
1_{\left\{- (1 - \delta) \sqrt{t} - \frac{a_n (t)}{n} s
\leq \widetilde{B}_s \leq y,~
\forall s \in [0,~ n - \ell(y)] \right\}} \psi (\widetilde{B}_{n - \ell (y)}) \right],
\end{align}
where
\begin{align} \label{eq:lem-limit-tail-1-1-2-psi}
\psi (x) := &\sqrt{\frac{\sqrt{t}}{\sqrt{t} + \frac{a_n (t)}{n} (n - \ell (y)) + x}} 
\cdot e^{- \frac{2 a_n (t)}{n} x - \frac{(a_n (t))^2}{n^2} (n - \ell (y))}
\notag \\
&\times \widetilde{P}_v \left(\max_{u \in T_{\ell (y)}^v} 
\sqrt{\widetilde{L}_{\widetilde{\tau}^{\downarrow} 
\left((\sqrt{t} + \frac{a_n (t)}{n} (n - \ell (y)) + x )^2 \right)}^{\downarrow} (u)}
> \sqrt{t} + a_n (t) + y \right).
\end{align}
By Lemma \ref{lem:density-barrier}, (\ref{eq:lem-limit-tail-1-1-2}) is bounded from above by
\begin{align} \label{eq:lem-limit-tail-1-1-3}
&\int_{[0,~y + y^{\frac{1}{20}}] \cup 
[- (1-\delta) \sqrt{t} - \frac{a_n (t)}{n} (n - \ell (y)),~- \ell(y)] }
 \mu_{n - \ell (y),~y}^* (x) \psi (x) dx \notag \\
&~~~+ \int_{- \ell(y)}^0
\delta_y \mu_{n - \ell (y),~y} (x) \psi (x) dx
\notag \\
& =: J_1 + J_2.
\end{align}
By Proposition \ref{prop:tail}(i)
and the assumption $\ell (y) > e^{\frac{8\sqrt{\log b}}{3} y^{1/20}}$, 
taking $n_0 (y, \ell(y)) \in \mathbb{N}$
large enough, we have
for all $n \geq n_0 (y, \ell(y))$
\begin{equation} \label{eq:lem-limit-tail-1-2}
J_1 \leq c_1 b^{- (n - \ell (y))} y e^{- 2 \sqrt{\log b}~y} (\ell (y))^{- 1/2}.
\end{equation}

By the change of measure (\ref{eq:change-measure}),
$J_2$ is bounded from above by
\begin{align} \label{eq:lem-limit-tail-1-3}
&\delta_y E_0^B \left[ 
\begin{minipage}{230pt}
$1_{\{\frac{a_n (t)}{n}(n - \ell (y)) - \ell (y) \leq B_{n - \ell (y)} \leq \frac{a_n (t)}{n} (n - \ell (y)) \}}$
$\times 1_{\left\{- (1 - \delta) \sqrt{t} \leq B_s 
\leq \frac{a_n (t)}{n} s + y,~\forall s \in [0,~n - \ell (y)] \right\}} 
\sqrt{\frac{\sqrt{t}}{\sqrt{t} + B_{n - \ell (y)}}}$ \\
$\times \widetilde{P}_{v} \left( \max_{u \in T_{\ell (y)}^v}
\sqrt{\widetilde{L}_{\widetilde{\tau}^{\downarrow} \left((\sqrt{t}+B_{n - \ell (y)})^2 \right)}
^{\downarrow} (u)} > \sqrt{t} + a_n (t) + y \right)$
\end{minipage}
\right] \notag \\
&+ \delta_y E_0^B \left[ 
\begin{minipage}{280pt}
$1_{\left\{\frac{a_n (t)}{n} (n - \ell (y)) - \ell (y) \leq B_{n - \ell (y)}
\leq \frac{a_n (t)}{n} (n - \ell(y)),
~\displaystyle \min_{s \in [0,~ n - \ell (y)]} B_s < - (1-\delta) \sqrt{t} \right\}}$ \\
$\times
\sqrt{\frac{\sqrt{t}}{\sqrt{t} + B_{n - \ell (y)}}}
\widetilde{P}_{v} \left(\displaystyle \max_{u \in T_{\ell (y)}^v}
\sqrt{\widetilde{L}_{\widetilde{\tau}^{\downarrow} \left((\sqrt{t}+B_{n - \ell (y)})^2 \right)}
^{\downarrow} (u)} > \sqrt{t} + a_n (t) + y \right)$
\end{minipage}
\right] \notag \\
&=: J_{2, 1} + J_{2, 2}.
\end{align}

By Lemma \ref{lem:local-time-bessel} and (\ref{eq:0-1-dim-bessel}), we have for each $n \geq 1$
\begin{equation} \label{eq:lem-limit-tail-1-4}
J_{2, 1} \leq c_2 \delta_y \widetilde{P}_{\rho} (F_{v, y, \ell(y)}^n (t)),
\end{equation}
where we have used the fact that under the event that
$\sqrt{X_s/2} \geq \delta \sqrt{t}$ for all $0 \leq s \leq n - \ell (y)$, 
we have $\exp \left(\frac{3}{8} \int_0^{n - \ell(y)} \frac{ds}{X_s} \right) \leq c_2$
under the assumption $t \geq c_* n \log n$.  

By the symmetry of $B$, (\ref{eq:joint-max-density}), and Proposition \ref{prop:tail}(i), 
taking $n_0 = n_0 (y, \ell(y)) \in \mathbb{N}$ large enough,
we have for all $n \geq n_0$
\begin{equation} \label{eq:lem-limit-tail-1-5}
J_{2, 2} 
\leq c_3 \delta_y b^{- (n - \ell (y))} e^{- 2 \sqrt{\log b}~y},
\end{equation}
where we have used the assumption $t \geq c_* n \log n$, $c_*$ is large enough.
Thus, by (\ref{eq:lem-limit-tail-1-1})-(\ref{eq:lem-limit-tail-1-5}), we have
\begin{equation} \label{eq:lem-limit-tail-1-6}
(1 + c_2 \delta_y) \widetilde{E}_{\rho} (\Lambda_{y, \ell(y)}^n (t))
\geq \widetilde{E}_{\rho} (\widetilde{\Lambda}_{y, \ell(y)}^n (t))
- c_1 y e^{- 2 \sqrt{\log b}~y} (\ell (y))^{- 1/2} 
- c_3 \delta_y e^{- 2 \sqrt{\log b}~y}.
\end{equation}

In the remainder of the proof, we obtain a lower bound of 
$\widetilde{E}_{\rho} (\widetilde{\Lambda}_{y, \ell(y)}^n (t))$.
Recall the definition of the event $G_{y + y^{1/20} - 2}^{n - \ell (y)} (t)$
from (\ref{eq:event-broken-barrier}).
Let $\widetilde{G}$ be the slightly modified version of this event given by
\begin{equation*}
\left\{ \begin{minipage}{260pt}
$\exists v \in T_{n - \ell (y)}, \exists s \in [0, n - \ell(y)]:$ \\
$\sqrt{\widetilde{L}_{\widetilde{\tau} (t)}^n (v_s)} 
> \sqrt{t} + \frac{a_n (t)}{n} s + \kappa (\log (s \wedge (n - \ell (y) -s)))_+ + y + y^{1/20}$
\end{minipage} \right \}.
\end{equation*}
$\widetilde{E}_{\rho} (\widetilde{\Lambda}_{y, \ell(y)}^n (t))$ is bounded from below by
\begin{align} \label{eq:lem-limit-tail-1-7}
&\widetilde{P}_{\rho} \left[
\left \{\exists v \in T_{n - \ell (y)} :~
\begin{minipage}{150pt}
$\min_{s \in [0,~ n - \ell (y)]}
\sqrt{\widetilde{L}_{\widetilde{\tau} (t)}^n (v_s)} \geq \delta \sqrt{t}$, \\
$\max_{u \in T_{\ell (y)}^v} \sqrt{\widetilde{L}_{\widetilde{\tau} (t)}^n (u)} > \sqrt{t} + a_n (t) + y$
\end{minipage}
\right \} \cap \left(\widetilde{G} \right)^c
 \right] \notag \\
&\geq P_{\rho} \left[\max_{u \in T_n} \sqrt{L_{\tau (t)}^n (u)} > \sqrt{t} + a_n (t) + y \right]
- \widetilde{P}_{\rho} \left[G_{y + y^{1/20} - 2}^{n - \ell (y)} (t) \right] \notag \\
&- \widetilde{P}_{\rho} \left[
\exists v \in T_{n - \ell (y)} :~
\begin{minipage}{150pt} 
$\min_{s \in [0,~ n - \ell (y)]} \sqrt{\widetilde{L}_{\widetilde{\tau} (t)}^n (v_s)} < \delta \sqrt{t}$, \\
$\max_{u \in T_{\ell (y)}^v} \sqrt{\widetilde{L}_{\widetilde{\tau} (t)}^n (u)} > \sqrt{t} + a_n (t) + y$
\end{minipage}
\right],
\end{align}
where we have used the inequality
$\frac{a_n (t)}{n} s \geq \frac{a_{n - \ell (y)} (t)}{n - \ell (y)} s - 1, ~s \in [0, n - \ell (y)]$,
which implies that
\begin{equation} \label{eq:barrier-variation}
\widetilde{P}_{\rho} [\widetilde{G}] 
\leq \widetilde{P}_{\rho} \left[G_{y + y^{1/20} - 2}^{n - \ell (y)} (t) \right].
\end{equation}
(Note that we have also used the fact that
the law of $\{\widetilde{L}_{\widetilde{\tau} (t)}^n (x) : x \in \widetilde{T}_{\leq n - \ell (y)} \}$
is the same as that of 
$\{\widetilde{L}_{\widetilde{\tau} (t)}^{n - \ell (y)} (x) : x \in \widetilde{T}_{\leq n - \ell (y)} \}$.)

Fix $v \in T_{n - \ell (y)}$. 
Recall the definitions of $\widetilde{L}^{\downarrow}$ and $\widetilde{\tau}^{\downarrow} (\cdot)$
from the beginning of the proof.
By Lemma \ref{lem:markov-local-time},
the third term on the right-hand side of (\ref{eq:lem-limit-tail-1-7})
is bounded from above by
\begin{equation} \label{eq:lem-limit-tail-1-10-1}
b^{n - \ell (y)} \widetilde{E}_{\rho} 
\left[
\begin{minipage}{220pt}
$1_{\left\{\min_{s \in [0,~ n - \ell (y)]} \sqrt{\widetilde{L}_{\widetilde{\tau} (t)}^n (v_s)} < \delta \sqrt{t},
~\sqrt{\widetilde{L}_{\widetilde{\tau} (t)}^n (v)} > 0 \right \}}$ \\
$\times \widetilde{P}_{v} \left(\max_{u \in T_{\ell (y)}^{v}}
\sqrt{\widetilde{L}_{\widetilde{\tau}^{\downarrow} 
\left(\widetilde{L}_{\widetilde{\tau} (t)}^n (v) \right)}^{\downarrow} (u)}
> \sqrt{t} + a_n (t) + y \right)$
\end{minipage}
\right].
\end{equation}
By Lemma \ref{lem:local-time-bessel}, (\ref{eq:0-1-dim-bessel}),
and the symmetry of the Brownian motion, 
(\ref{eq:lem-limit-tail-1-10-1}) is bounded from above by
\begin{equation} \label{eq:lem-limit-tail-1-10-2}
b^{n - \ell (y)} E_0^B \left[
\begin{minipage}{240pt}
$1_{\left\{\max_{s \in [0,~ n - \ell (y)]} B_s >  (1 - \delta) \sqrt{t},~
B_{n - \ell (y)} <  \sqrt{t} \right \}} 
\sqrt{\frac{\sqrt{t}}{\sqrt{t} - B_{n - \ell (y)}}}$ \\
$\times \widetilde{P}_{v} \left(\max_{u \in T_{\ell (y)}^{v}}
\sqrt{\widetilde{L}_{\widetilde{\tau}^{\downarrow} 
\left(\left(\sqrt{t} - B_{n - \ell (y)} \right)^2 \right) }^{\downarrow} (u)}
> \sqrt{t} + a_n (t) + y \right)$
\end{minipage}
\right].
\end{equation}
We estimate the indicator function in the expectation in (\ref{eq:lem-limit-tail-1-10-2}) 
from above
by $\sum_{i=1}^3 1_{E_i}$, where
\begin{align*}
&E_1 := \{(1 - \delta) \sqrt{t} \leq B_{n - \ell (y)} \leq \sqrt{t} \}, \\
&E_2 := \{- a_n (t) + \sqrt{\log b}~ \ell (y) \leq B_{n - \ell (y)} < (1 - \delta) \sqrt{t},~~
\max_{s \in [0, n - \ell (y)]} B_s > (1 - \delta) \sqrt{t} \}, \\
&E_3 := \{B_{n - \ell (y)} < - a_n (t) + \sqrt{\log b}~ \ell (y),~~
\max_{s \in [0, n - \ell (y)]} B_s > (1 - \delta) \sqrt{t} \}.
\end{align*}
Let $H_i, i \in \{1, 2, 3 \}$ be the expectation obtained 
from the one in (\ref{eq:lem-limit-tail-1-10-2})
by replacing the indicator function in it with $1_{E_i}$.
In particular, the right of (\ref{eq:lem-limit-tail-1-10-2}) is bounded from above by
$b^{n - \ell (y)} (H_1 + H_2 + H_3)$.

To estimate $H_1$, we use Proposition \ref{prop:tail}(i)
and the density of $B_{n - \ell (y)}$.
To estimate $H_2$,
we use Proposition \ref{prop:tail}(i) and (\ref{eq:joint-max-density}).
To estimate $H_3$, we use (\ref{eq:joint-max-density})
and bound the probability in $H_3$ from above just by $1$. 
Taking $n_0 = n_0(y, \ell (y)) \in \mathbb{N}$ large enough,
for all $n \geq n_0$, 
(\ref{eq:lem-limit-tail-1-10-2}) is bounded from above by
\begin{equation} \label{eq:lem-limit-tail-1-10-3}
c_4 \varepsilon_n y e^{- 2 \sqrt{\log b}~y},
\end{equation}
where $\varepsilon_n, n \geq 0$ is a sequence with $\varepsilon_n \to 0$ as $n \to \infty$.
By (\ref{eq:lem-limit-tail-1-7}) and (\ref{eq:lem-limit-tail-1-10-3})
together with Proposition \ref{prop:tail}(ii) and Lemma \ref{lem:broken-barrier},
taking $n_0 = n_0 (y, \ell(y)) \in \mathbb{N}$ large enough,
we have
for all $n \geq n_0$
\begin{equation} \label{eq:lem-limit-tail-1-11}
\widetilde{E}_{\rho} (\widetilde{\Lambda}_{y, \ell(y)}^n (t))
\geq c_5 y e^{- 2 \sqrt{\log b}~y}.
\end{equation}
By (\ref{eq:lem-limit-tail-1-6}) and (\ref{eq:lem-limit-tail-1-11}), we have
(\ref{eq:lem-limit-tail-1-statement}).~~~$\Box$\\
To obtain a lower bound of
$P_{\rho} (\max_{v \in T_n} \sqrt{L_{\tau (t)}^n (v)} \geq \sqrt{t} + a_n (t) + y)
/\widetilde{E}_{\rho} [\Lambda_{y, \ell(y)}^n (t)]$,
we need the following:
\begin{lem} \label{lem:limit-tail-2}
There exist $c_1, c_2 \in (0, ~\infty)$, $y_0 > 0$ such that the following holds:
for all $y \geq y_0$
and $\ell(y) > e^{\frac{8 \sqrt{\log b}}{3} y^{1/20}}$,
there exists $n_0 = n_0 (y, \ell(y)) \in \mathbb{N}$ such that
for all $n \geq n_0$ and $t \geq c_1 n \log n$,
\begin{equation} \label{eq:lem-limit-tail-2-statement}
\frac{\widetilde{E}_{\rho} 
\left[\left(\Lambda_{y, \ell(y)}^n (t) \right)^2 \right]}
{\widetilde{E}_{\rho} \left[\Lambda_{y, \ell(y)}^n (t) \right]}
\leq 1 + c_2 y^{- 1/2}.
\end{equation}
\end{lem} 
{\it Proof.}
Fix any $y \geq y_0$ and $\ell (y) > e^{\frac{8 \sqrt{\log b}}{3} y^{1/20}}$,
where we take $y_0 > 0$ large enough.
Throughout the proof, given $n \in \mathbb{N}$, we assume 
$t \geq c_* n \log n$ for some sufficiently large $c_* > 0$.
We have
\begin{equation} \label{eq:lem-limit-tail-2-1}
\widetilde{E}_{\rho} \left[\left(\Lambda_{y, \ell(y)}^n (t) \right)^2 \right]
= \widetilde{E}_{\rho} [\Lambda_{y, \ell(y)}^n (t)]
+ \sum_{k = 1}^{n - \ell (y)} 
\sum_{\begin{subarray}{c} v, w \in T_{n - \ell (y)}, \\ |v \wedge w| = n - \ell (y) - k \end{subarray}}
\widetilde{P}_{\rho} \left(F_{v, y, \ell(y)}^n (t) \cap F_{w, y, \ell(y)}^n (t) \right).
\end{equation}
Fix $1 \leq k \leq n - \ell (y) - 1$ 
and $v, w \in T_{n - \ell (y)}$ with $|v \wedge w| = n - \ell (y) - k$.
Let $\widetilde{L}^{\downarrow}$ be a local time of a Brownian motion on
$\widetilde{T}_{\leq \ell(y) + k}^{v \wedge w}$ and set
$\widetilde{\tau}^{\downarrow} (s) := \inf \{r \geq 0 : \widetilde{L}_r ^{\downarrow} (v \wedge w) > s \}$.
By Lemma \ref{lem:markov-local-time},
$\widetilde{P}_{\rho} (F_{v, y, \ell(y)}^n (t) \cap F_{w, y, \ell(y)}^n (t))$ is equal to
\begin{equation} \label{eq:lem-limit-tail-2-2}
\widetilde{E}_{\rho} \Biggl[
1_{\left\{\delta \sqrt{t} \leq \sqrt{\widetilde{L}_{\widetilde{\tau} (t)}^n (v_s)}
\leq \sqrt{t} + \frac{a_n (t)}{n} s + y,~
\forall s \in [0,~n - \ell (y) - k] \right \}} \prod_{x \in \{v,~w \}} P^x \Biggr],
\end{equation}
where for each $x \in \{v, w \}$,
\begin{equation*}
P^x := \widetilde{P}_{v \wedge w} \left(
\begin{minipage}{190pt}
$\forall s \in [n - \ell (y) - k, ~n - \ell (y)]$, \\
$\delta \sqrt{t} \leq 
\sqrt{\widetilde{L}_{\widetilde{\tau}^{\downarrow} \left(\widetilde{L}_{\widetilde{\tau} (t)}^n (v \wedge w) \right)}
^{\downarrow} (x_s)} \leq \sqrt{t} + \frac{a_n (t)}{n} s + y$, \\ 
$\max_{u \in T_{\ell (y)}^x} 
\sqrt{\widetilde{L}_{\widetilde{\tau}^{\downarrow} 
\left(\widetilde{L}_{\widetilde{\tau} (t)}^n (v \wedge w) \right)}^{\downarrow} (u)} > \sqrt{t} + a_n (t) + y$
\end{minipage}
\right).
\end{equation*}

Fix $0 \leq i \leq \lfloor (1 - \delta) \sqrt{t} + \frac{a_n (t)}{n} (n - \ell (y) - k) + y \rfloor$.
Assume that 
\begin{equation} \label{eq:lem-limit-tail-2-assumption}
\sqrt{\widetilde{L}_{\widetilde{\tau} (t)}^n (v \wedge w)} 
\in \sqrt{t} + \frac{a_n (t)}{n} (n - \ell (y) - k) + y - i + \left(- 1, ~0 \right].
\end{equation}
Under the assumption (\ref{eq:lem-limit-tail-2-assumption}), 
we estimate the probabilities in (\ref{eq:lem-limit-tail-2-2}).
Fix $x \in \{v, w \}$.
Let $\widetilde{L}^{\downarrow \downarrow}$ be a local time of a Brownian motion
on $\widetilde{T}_{\leq \ell(y)}^x$.
We define the inverse local time by
$\widetilde{\tau}^{\downarrow \downarrow} (s) 
:= \inf \{r \geq 0 : \widetilde{L}_r^{\downarrow \downarrow} (x) > s \}$.
By Lemma \ref{lem:markov-local-time},
$P^x$ is equal to
\begin{equation} \label{eq:lem-limit-tail-2-3}
\widetilde{E}_{v \wedge w} \left[
\begin{minipage}{260pt}
$1_{\left \{\delta \sqrt{t} \leq 
\sqrt{\widetilde{L}_{\widetilde{\tau}^{\downarrow} \left(\widetilde{L}_{\widetilde{\tau} (t)}^n (v \wedge w) \right)}
^{\downarrow} (x_s)} \leq \sqrt{t} + \frac{a_n (t)}{n} s + y,
~\forall s \in [n - \ell (y) - k, ~n - \ell (y)] \right \}}$ \\
$\times \widetilde{P}_x \left(\max_{u \in T_{\ell (y)}^x}
\sqrt{\widetilde{L}_{\widetilde{\tau}^{\downarrow \downarrow} 
\left(\widetilde{L}_{\widetilde{\tau}^{\downarrow} \left(\widetilde{L}_{\widetilde{\tau} (t)}^n (v \wedge w) \right)}
^{\downarrow} (x)  \right)}
^{\downarrow \downarrow} (u)} > \sqrt{t} + a_n (t) + y \right)$
\end{minipage}
\right].
\end{equation}
By Lemma \ref{lem:local-time-bessel}, (\ref{eq:0-1-dim-bessel}), 
and the change of measure (\ref{eq:change-measure})-(\ref{eq:bm-drift}),
(\ref{eq:lem-limit-tail-2-3}) is bounded from above by
\begin{equation} \label{eq:lem-limit-tail-2-3-2}
\widetilde{E}_k^B \left[
\begin{minipage}{270pt}
$e^{- \frac{2 a_n (t)}{n} \widetilde{B}_k - \frac{(a_n (t))^2}{n^2} k}
\sqrt{\frac{\sqrt{\widetilde{L}_{\widetilde{\tau} (t)}^n (v \wedge w)}}
{\sqrt{\widetilde{L}_{\widetilde{\tau} (t)}^n (v \wedge w)} + \frac{a_n (t)}{n} k + \widetilde{B}_k}}$ \\ 
$\times 1_{\left \{\delta \sqrt{t} - \frac{a_n (t)}{n}s \leq 
\sqrt{\widetilde{L}_{\widetilde{\tau} (t)}^n (v \wedge w)} + \widetilde{B}_s
\leq \sqrt{t} + \frac{a_n (t)}{n} (n - \ell (y) - k) + y,
~\forall s \in [0, k] \right \}}$ \\
$\times \widetilde{P}_x \left(\displaystyle \max_{u \in T_{\ell (y)}^x}
\sqrt{\widetilde{L}_{\widetilde{\tau}^{\downarrow \downarrow} 
\left(\left(\sqrt{\widetilde{L}_{\widetilde{\tau} (t)}^n (v \wedge w)} 
+ \frac{a_n (t)}{n} k + \widetilde{B}_k \right)^2 \right)}
^{\downarrow \downarrow} (u)} > \sqrt{t} + a_n (t) + y \right)$
\end{minipage}
\right].
\end{equation}
By (\ref{eq:lem-limit-tail-2-assumption})
and (\ref{eq:density-barrier-standard}),
(\ref{eq:lem-limit-tail-2-3-2}) is bounded from above by
\begin{align} \label{eq:lem-limit-tail-2-3-3}
&\int_0^{(1 -\delta) \sqrt{t} + \frac{a_n (t)}{n} (n - \ell (y)) + y + 1}
\frac{1}{\sqrt{\pi k}} \left(e^{- \frac{(i+1-z)^2}{k}} - e^{- \frac{(i+1+z)^2}{k}} \right) \notag \\
&~~~~~~~~~~~~~~\times
e^{\frac{2 a_n (t)}{n} z - \frac{2 a_n (t)}{n} (i+1) - \frac{(a_n (t))^2}{n^2} k}
\sqrt{\frac{\sqrt{t} + \frac{a_n (t)}{n} (n - \ell (y) - k) + y - i}
{\sqrt{t} + \frac{a_n (t)}{n} (n - \ell (y)) + y - z}} \notag \\
&~~~~~~~~~~~~~~~~~~\times \widetilde{P}_x \left(\max_{u \in T_{\ell (y)}^x}
\sqrt{\widetilde{L}_{\widetilde{\tau}^{\downarrow \downarrow} 
\left(\left(\sqrt{t} + \frac{a_n (t)}{n} (n - \ell (y)) + y + 1 - z \right)^2 \right)}
^{\downarrow \downarrow} (u)} > \sqrt{t} + a_n (t) + y \right) dz.
\end{align}
We use the following estimate in the integrand of (\ref{eq:lem-limit-tail-2-3-3}):
$e^{- \frac{(i+1-z)^2}{k}} - e^{- \frac{(i+1+z)^2}{k}} \leq 1$
for each $1 \leq k \leq \lfloor y \rfloor$, and 
$e^{- \frac{(i+1-z)^2}{k}} - e^{- \frac{(i+1+z)^2}{k}} \leq \frac{4(i+1)z}{k}$ 
for each $\lfloor y \rfloor \leq k \leq n - \ell (y) - 1$.
By this and Proposition \ref{prop:tail}(i), 
the right-hand side of (\ref{eq:lem-limit-tail-2-3-3}) is bounded from above by
\begin{equation} \label{eq:lem-limit-tail-2-6}
c_1 A b^{- k}  e^{- \frac{2 a_n (t)}{n} i}
e^{\frac{3 k \log n}{2n}} e^{\frac{k \log \left(\frac{\sqrt{t} + n}{\sqrt{t}} \right)}{2n}}
\sqrt{\frac{\sqrt{t} + \frac{a_n (t)}{n} (n - \ell (y) - k) + y - i}{\sqrt{t} + n}},
\end{equation}
where
$A := k^{- 1/2} (\ell (y))^{- 1/2}$
if $1 \leq k \leq \lfloor y \rfloor$, and
$A := k^{- 3/2} (i+1)$
if $\lfloor y \rfloor \leq k \leq n - \ell (y) - 1$.
Recall the events in the indicator function in (\ref{eq:lem-limit-tail-2-2})
and in (\ref{eq:lem-limit-tail-2-assumption}).
We estimate the probability of the intersection of these events.
Using Lemma \ref{lem:local-time-bessel}, (\ref{eq:0-1-dim-bessel}), and
and the change of measure (\ref{eq:change-measure})-(\ref{eq:bm-drift})
for $n \geq n_0$ ($n_0 = n_0(y, \ell (y)) \in \mathbb{N}$ large enough),
we have
\begin{align} \label{eq:lem-limit-tail-2-7}
&\widetilde{P}_{\rho} \left(
\begin{minipage}{240pt}
$\delta \sqrt{t} \leq \sqrt{\widetilde{L}_{\widetilde{\tau} (t)}^n (v_s)}
\leq \sqrt{t} + \frac{a_n (t)}{n} s + y,~
\forall s \in [0,~ n - \ell (y) - k]$, \\
$\sqrt{\widetilde{L}_{\widetilde{\tau} (t)}^n (v \wedge w)} 
\in \sqrt{t} + \frac{a_n (t)}{n} (n - \ell (y) - k) + y - i + \left(- 1, ~0 \right]$
\end{minipage}
\right) \notag \\
&\leq \widetilde{E}_{n - \ell (y) - k}^B
\left[
\begin{minipage}{250pt}
$1_{\{\widetilde{B}_s \leq y,~\forall s \in [0,~n - \ell (y) - k],
~\widetilde{B}_{n - \ell (y) - k} - y \in (- i - 1,~- i] \}}$ \\
$\times e^{- \frac{2 a_n (t)}{n} \widetilde{B}_{n - \ell (y) - k} - \frac{(a_n (t))^2}{n^2} (n - \ell (y) - k)}
\sqrt{\frac{\sqrt{t}}{\sqrt{t} + \frac{a_n (t)}{n} (n - \ell (y) - k) + \widetilde{B}_{n - \ell (y) - k}}}$
\end{minipage}
\right].
\end{align}
By (\ref{eq:density-barrier-standard}), 
the right of (\ref{eq:lem-limit-tail-2-7}) is bounded from above by
\begin{align} \label{eq:lem-limit-tail-2-7-1}
&\int_i^{i+1} \frac{1}{\sqrt{\pi (n - \ell (y) -k)}} 
\left(e^{- \frac{(y - z)^2}{n - \ell (y) - k}} - e^{- \frac{(y + z)^2}{n - \ell (y) - k}} \right) \notag \\
&~~~~~~~~\times e^{\frac{2 a_n (t)}{n} z - \frac{2 a_n (t)}{n} y - \frac{(a_n (t))^2}{n^2} (n - \ell (y) - k)}
\sqrt{\frac{\sqrt{t}}{\sqrt{t} + \frac{a_n (t)}{n} (n - \ell (y) - k) + y - z}}~~~ dz.
\end{align}
We will use the following in the integrand of (\ref{eq:lem-limit-tail-2-7-1}):
$e^{- \frac{(y - z)^2}{n - \ell (y) - k}} - e^{- \frac{(y + z)^2}{n - \ell (y) - k}} \leq
\frac{4yz}{n - \ell (y) - k}$
for each $1 \leq k \leq n - \ell (y) - \lfloor y \rfloor$, and
$e^{- \frac{(y - z)^2}{n - \ell (y) - k}} - e^{- \frac{(y + z)^2}{n - \ell (y) - k}} \leq 1$
for each $n - \ell (y) - \lfloor y \rfloor \leq k \leq n - \ell (y) - 1$.
By this, the right-hand side of (\ref{eq:lem-limit-tail-2-7-1}) is bounded from above by
\begin{equation} \label{eq:lem-limit-tail-2-7-2}
\begin{minipage}{280pt}
$c_2 A^{\prime} b^{- (n - \ell (y) - k)}  e^{\frac{2 a_n (t)}{n} i} e^{- 2 \sqrt{\log b}~y}$ \\
$\times e^{\frac{3 (n - \ell (y) - k)}{2n} \log n}
e^{\frac{(n - \ell (y) - k) \log \left(\frac{\sqrt{t} + n}{\sqrt{t}} \right)}{2n}}
\sqrt{\frac{\sqrt{t}}{\sqrt{t} + \frac{a_n (t)}{n} (n - \ell (y) - k) + y - i - 1}}$,
\end{minipage}
\end{equation}
where
$A^{\prime} := (n - \ell (y) - k)^{- 3/2} (i+1) y$
if $1 \leq k \leq n - \ell (y) - \lfloor y \rfloor$, and
$A^{\prime} := (n - \ell (y) - k)^{- 1/2}$
if $n - \ell (y) - \lfloor y \rfloor \leq k \leq n - \ell (y) - 1$.

We divide the sum over $1 \leq k < n - \ell (y)$
in (\ref{eq:lem-limit-tail-2-1})
into sums over the following:
(a) $1 \leq k \leq \lfloor y \rfloor$;
(b) $\lfloor y \rfloor < k \leq \lfloor (n - \ell (y))/2 \rfloor$;
(c) $\lfloor (n - \ell (y))/2 \rfloor < k \leq n - \ell (y) - \lfloor y \rfloor$;
(d) $n - \ell (y) - \lfloor y \rfloor < k \leq n - \ell (y) - 1$.
We make remarks on how to estimate the sums
over difficult regimes, (b) and (c):
In the regime (b), we use the fact that
$k^{- \frac{3}{2}} e^{\frac{3k \log n}{2n}} e^{\frac{k \log (\frac{\sqrt{t} + n}{\sqrt{t}})}{2n}}$
is bounded from above by a universal constant.
In the regime of (c), we use the estimate
$$\frac{1}{(n - \ell (y) - k)^{\frac{3}{2}}} \sqrt{\frac{\sqrt{t} + \frac{a_n (t)}{n} (n - \ell (y) - k) + y}{\sqrt{t}}}
\leq \frac{c}{(n - \ell (y) - k)^{\frac{3}{2}}} + \frac{c}{t^{\frac{1}{4}} (n - \ell (y)  - k) },$$
for some universal constant $c$.

By (\ref{eq:lem-limit-tail-2-2}), (\ref{eq:lem-limit-tail-2-6}), 
and (\ref{eq:lem-limit-tail-2-7-2}), we have
\begin{equation} \label{eq:lem-limit-tail-2-8}
\sum_{k=1}^{n - \ell (y) - 1}
\sum_{\begin{subarray}{c} v, w \in T_{n - \ell (y)}, \\ |v \wedge w| = n - \ell (y) - k \end{subarray}} 
\widetilde{P}_{\rho} \left(F_{v, y, \ell(y)}^n (t) \cap F_{w, y, \ell(y)}^n (t) \right) 
\leq c_3 y^{1/2} e^{- 2 \sqrt{\log b}~y}.
\end{equation}
In the case $k = n - \ell (y)$, by the independence of 
excursions of a Brownian motion around $\rho$, we have
\begin{align} \label{eq:lem-limit-tail-2-9}
\sum_{\begin{subarray}{c} v, w \in T_{n - \ell (y)}, \\ |v \wedge w| = 0 \end{subarray}} 
\widetilde{P}_{\rho} \left(F_{v, y, \ell(y)}^n (t) \cap F_{w, y, \ell(y)}^n (t) \right) 
&= \sum_{\begin{subarray}{c} v, w \in T_{n - \ell (y)}, \\ |v \wedge w| = 0 \end{subarray}} 
\widetilde{P}_{\rho} \left(F_{v, y, \ell(y)}^n (t) \right) 
\widetilde{P}_{\rho} \left(F_{w, y, \ell(y)}^n (t) \right)  \notag \\
&\leq \left(\widetilde{E}_{\rho} [\Lambda_{y, \ell(y)}^n (t)] \right)^2.
\end{align}
Thus, by (\ref{eq:lem-limit-tail-2-1}), (\ref{eq:lem-limit-tail-2-8}), and (\ref{eq:lem-limit-tail-2-9}),
we have
\begin{equation} \label{eq:lem-limit-tail-2-10}
\widetilde{E}_{\rho} \left[\left(\Lambda_{y, \ell(y)}^n (t) \right)^2 \right]
\leq \widetilde{E}_{\rho} [\Lambda_{y, \ell(y)}^n (t)] + c_3 y^{1/2} e^{- 2 \sqrt{\log b}~y}
+ \left(\widetilde{E}_{\rho} [\Lambda_{y, \ell(y)}^n (t)] \right)^2.
\end{equation}
\\ \\
We will obtain upper and lower bounds of $\widetilde{E}_{\rho} [\Lambda_{y, \ell(y)}^n (t)]$.
By (\ref{eq:lem-limit-tail-1-6}) and (\ref{eq:lem-limit-tail-1-11}),
taking $n_0 = n_0 (y, \ell(y))$ large enough,
we have for all $n \geq n_0$
\begin{equation} \label{eq:lem-limit-tail-2-11}
\widetilde{E}_{\rho} [\Lambda_{y, \ell(y)}^n (t)]
\geq c_4 y e^{- 2 \sqrt{\log b}~y}.
\end{equation}
By the arguments in (\ref{eq:lem-limit-tail-1-1}) and (\ref{eq:lem-limit-tail-1-1-2}),
$\widetilde{E}_{\rho} [\Lambda_{y, \ell(y)}^n (t)]$ is bounded from above
by $b^{n - \ell (y)}$ times the second term of (\ref{eq:lem-limit-tail-1-1-2}).
By this together with (\ref{eq:density-barrier-standard}) and Proposition \ref{prop:tail}(i),
we have
\begin{equation} \label{eq:lem-limit-tail-2-12}
\widetilde{E}_{\rho} [\Lambda_{y, \ell(y)}^n (t)] \leq c_5 y e^{- 2 \sqrt{\log b}~y}.
\end{equation}
Thus, by (\ref{eq:lem-limit-tail-2-10})-(\ref{eq:lem-limit-tail-2-12}),
we have (\ref{eq:lem-limit-tail-2-statement}).~~~$\Box$\\
Using Lemma \ref{lem:limit-tail-1} and \ref{lem:limit-tail-2}, we prove the following:
\begin{lem} \label{lem:limit-tail-3}
There exist $c_1, c_2, c_3 \in (0, ~\infty)$, $y_0 > 0$, and
$\{\delta_{y^{\prime}} : y^{\prime} > 0 \}$ with $\lim_{y^{\prime} \to \infty} \delta_{y^{\prime}} = 0$
such that the following holds:
for all $y \geq y_0$ and $\ell(y) > e^{\frac{8 \sqrt{\log b}}{3} y^{1/20}}$, 
there exists $n_0 = n_0(y, \ell(y)) \in \mathbb{N}$ such that
for all $n \geq n_0$ and $t \geq c_1 n \log n$,
\begin{equation} \label{eq:lem-limit-tail-3-statement-1}
\frac{P_{\rho} \left(\max_{v \in T_n} \sqrt{L_{\tau (t)}^n (v)} > \sqrt{t} + a_n (t) + y \right) }
{\widetilde{E}_{\rho} \left[\Lambda_{y, \ell(y)}^n (t) \right]}  \\
\leq 
(1 + \delta_y) \left(1 + e^{- c_2 ~y^{1/20}} + \delta_y \right),
\end{equation}
\begin{equation} \label{eq:lem-limit-tail-3-statement-2}
\frac{P_{\rho} \left(\max_{v \in T_n} \sqrt{L_{\tau (t)}^n (v)} > \sqrt{t} + a_n (t) + y \right) }
{\widetilde{E}_{\rho} \left[\Lambda_{y, \ell(y)}^n (t) \right]}
\geq 1 - c_3 y^{- 1/2}.
\end{equation}
\end{lem}
{\it Proof.}
Fix any $y \geq y_0$, $\ell(y) > e^{\frac{8 \sqrt{\log b}}{3} y^{1/20}}$,
$n \geq n_0$, and $t \geq c_* n \log n$,
where we take $y_0 > 0$, $n_0 = n_0 (y, \ell(y)) \in \mathbb{N}$, $c_* > 0$ large enough.
We first obtain the upper bound.
Recall the event $G_{y + y^{1/20} - 2}^{n - \ell (y)} (t)$ from (\ref{eq:event-broken-barrier})
and $\varepsilon_n$ from (\ref{eq:lem-limit-tail-1-10-3}).
Recall the inequality in (\ref{eq:barrier-variation}).
We have
\begin{align} \label{eq:lem-limit-tail-3-1}
&~~~~P_{\rho} \left(\max_{v \in T_n} \sqrt{L_{\tau (t)}^n (v)} > \sqrt{t} + a_n (t) + y \right) \notag \\
&= P_{\rho} \left(
\exists v \in T_{n - \ell (y)} :
\max_{u \in T_{\ell (y)}^v} \sqrt{L_{\tau (t)}^n (u)} > \sqrt{t} + a_n (t) + y \right) \notag \\
&\leq \widetilde{E}_{\rho} [\widetilde{\Lambda}_{y, \ell(y)}^n (t)]
+ \widetilde{P}_{\rho} \left(G_{y + y^{1/20} - 2}^{n - \ell (y)} (t)\right) \notag \\
&~~~+\widetilde{P}_{\rho} \left(\exists v \in T_{n - \ell (y)} :~
\begin{minipage}{160pt}
$\max_{u \in T_{\ell (y)}^v} \sqrt{\widetilde{L}_{\widetilde{\tau} (t)}^n (u)} > \sqrt{t} + a_n (t) + y$, \\
$\min_{s \in [0,~ n - \ell (y)]} \sqrt{\widetilde{L}_{\widetilde{\tau} (t)}^n (v_s)} < \delta \sqrt{t}$
\end{minipage}
\right) \notag \\
&\leq \left(1 + c_1 e^{- 2 \sqrt{\log b} ~y^{1/20}} 
+ \varepsilon_n \right)
\widetilde{E}_{\rho} [\widetilde{\Lambda}_{y, \ell(y)}^n (t)],
\end{align}
where we have used Lemma \ref{lem:broken-barrier},
(\ref{eq:lem-limit-tail-1-10-3}), and (\ref{eq:lem-limit-tail-1-11}) in the last inequality.
By (\ref{eq:lem-limit-tail-3-1}) and Lemma \ref{lem:limit-tail-1}, we have
(\ref{eq:lem-limit-tail-3-statement-1}).
By Lemma \ref{lem:limit-tail-2}, we have
\begin{align*} 
&~~~~P_{\rho} \left(\max_{v \in T_n} \sqrt{L_{\tau (t)}^n (v)} > \sqrt{t} + a_n (t) + y \right) 
\geq \widetilde{P}_{\rho} \left(\Lambda_{y, \ell(y)}^n (t) \geq 1 \right) 
\geq \frac{\left(\widetilde{E}_{\rho} [\Lambda_{y, \ell(y)}^n (t)] \right)^2}
{\widetilde{E}_{\rho} \left[\left(\Lambda_{y, \ell(y)}^n (t) \right)^2 \right]} \notag \\
&\geq \left(1 - c_2 y^{- 1/2} \right) 
\widetilde{E}_{\rho} [\Lambda_{y, \ell(y)}^n (t)],
\end{align*}
which proves (\ref{eq:lem-limit-tail-3-statement-2}).~~~$\Box$
\\ \\
For each interval $I \subset \mathbb{R}$, set
\begin{equation} \label{eq:limit-tail-lambda-interval}
\Lambda_{y, \ell(y), I}^n (t) 
:= \sum_{v \in T_{n - \ell (y)}} 1_{F_{v, y, \ell(y)}^n (t) \cap 
\left\{\sqrt{\widetilde{L}_{\widetilde{\tau} (t)}^n (v)} \in \sqrt{t} + \frac{a_n (t)}{n} (n - \ell (y)) + I 
\right \}}.
\end{equation}
Set the interval
\begin{equation} \label{eq:limit-tail-interval}
J_{\ell (y)} := \left(- \ell (y) + y, - (\ell (y))^{2/5} + y \right].
\end{equation}
Then the following holds:
\begin{lem} \label{lem:limit-tail-4}
There exist $c_1, c_2 \in (0,~\infty)$ and $y_0 > 0$ such that the following holds:
for all $y \geq y_0$ and $\ell(y) > e^{\frac{8 \sqrt{\log b}}{3} y^{1/20}}$, 
there exists $n_0 = n_0 (y, \ell(y)) \in \mathbb{N}$ such that
for all $n \geq n_0$ and $t \geq c_1 n \log n$,
\begin{equation} \label{eq:lem-limit-tail-4-statement}
\frac{\widetilde{E}_{\rho} \left[\Lambda_{y, \ell(y), J_{\ell (y)}}^n (t) \right]}
{\widetilde{E}_{\rho} \left[\Lambda_{y, \ell(y)}^n (t) \right]}
\geq 1 - c_2 (\ell (y))^{- 1/5}.
\end{equation}
\end{lem}
{\it Proof.}
Fix any $y \geq y_0$ and $\ell(y) > e^{\frac{8 \sqrt{\log b}}{3} y^{1/20}}$,
where we take $y_0 > 0$ large enough.
Throughout the proof, given $n \in \mathbb{N}$,
we assume that $t \geq c_* n \log n$ for some sufficiently large $c_* > 0$.
Fix $v \in T_{n - \ell (y)}$ and any interval $I \subset \mathbb{R}$.
Recall the definition of $\psi$ from (\ref{eq:lem-limit-tail-1-1-2-psi}).
By similar arguments to those in (\ref{eq:lem-limit-tail-1-1}) and (\ref{eq:lem-limit-tail-1-1-2}),
taking $n_0 = n_0(y, \ell(y)) \in \mathbb{N}$ large enough, we have for all $n \geq n_0$, 
\begin{align} \label{eq:lem-limit-tail-4-1}
&\widetilde{P}_{\rho} \left(F_{v, y, \ell(y)}^n (t) \cap 
\left\{\sqrt{\widetilde{L}_{\widetilde{\tau} (t)}^n (v)} \in \sqrt{t} + \frac{a_n (t)}{n} (n - \ell (y)) + I 
\right \} \right) \notag \\
&\leq \widetilde{E}_{n - \ell (y)}^B \left[
1_{\left\{- (1 - \delta) \sqrt{t} - \frac{a_n (t)}{n} s \leq \widetilde{B}_{s} \leq y,
~\forall s \in [0,~n - \ell (y)],~\widetilde{B}_{n - \ell (y)} \in I \right\}}
\psi (\widetilde{B}_{n - \ell (y)})
\right].
\end{align}
By Proposition \ref{prop:tail}(i) and (\ref{eq:density-barrier-standard}),
(\ref{eq:lem-limit-tail-4-1}) is bounded from above by
\begin{align} \label{eq:lem-limit-tail-4-1-2}
&c_1 b^{- (n - \ell (y))} (\ell (y))^{- 3/2} y e^{- 2 \sqrt{\log b}~y} \notag \\
&\times \int_{(y - I) \cap I_{n, t, y}}
z (\log \ell (y) + z) 
e^{- c_2 \frac{z^2}{\ell (y)}} 
\sqrt{\frac{\sqrt{t} + n}{\sqrt{t} + \frac{a_n (t)}{n}(n - \ell (y)) + y - z + \ell (y) }} ~~dz,
\end{align}
where $I_{n, t, y} := [0,~(1 - \delta) \sqrt{t} + \frac{a_n (t)}{n} (n - \ell (y)) + y]$
and we have used the inequality
$e^{- \frac{(y - z)^2}{n - \ell (y)}} - e^{- \frac{(y + z)^2}{n - \ell (y)}} \leq \frac{4yz}{n - \ell (y)}$
for $z \geq 0$.
In the cases $I = (- \infty,~y - \ell (y)]$ and $I = \left(y- (\ell (y))^{\frac{2}{5}},~\infty \right)$,
the right-hand side of (\ref{eq:lem-limit-tail-4-1-2}) is bounded from above by
\begin{equation*}
c_3 b^{- (n - \ell (y))} (\ell (y))^{- 1/5} y e^{- 2 \sqrt{\log b}~y}
\end{equation*}
for all $n \geq n_0$, where $n_0 = n_0 (y, \ell(y)) \in \mathbb{N}$ large enough.
By this and (\ref{eq:lem-limit-tail-2-11}), we have
\begin{align*} 
&~~~~\widetilde{E}_{\rho} \left[\Lambda_{y, \ell(y)}^n (t) \right] \notag \\
&= \widetilde{E}_{\rho} \left[\Lambda_{y, \ell(y), J_{\ell (y)}}^n (t) \right]
+ \widetilde{E}_{\rho} \left[\Lambda_{y, \ell(y), (- \infty,~y - \ell (y)]}^n (t) \right]
+ \widetilde{E}_{\rho} \left[\Lambda_{y, \ell(y), \left(y - (\ell (y))^{2/5},~\infty \right)}^n (t) \right] 
\notag \\
&\leq \widetilde{E}_{\rho} \left[\Lambda_{y, \ell(y), J_{\ell (y)}}^n (t) \right]
+ c_4 (\ell (y))^{- 1/5} \widetilde{E}_{\rho} \left[\Lambda_{y, \ell(y)}^n (t) \right],
\end{align*}
which implies (\ref{eq:lem-limit-tail-4-statement}).~~~$\Box$
\\ \\
{\it Proof of Proposition \ref{prop:limit-tail}.}
Let $(h_v)_{v \in T}$ be a BRW on $T$ defined in Section \ref{sec:intro}.
By Lemma \ref{lem:tail-max-brw-1}, one can show that the sequence
\begin{equation} \label{eq:sequence-gamma}
\left(\int_{\ell^{2/5}}^{\ell}
z e^{2 \sqrt{\log b}~z}
\mathbb{P} \left[\max_{u \in T_{\ell}} h_u
> \sqrt{\log b}~\ell + z \right] dz \right)_{\ell \geq 1}
\end{equation}
is bounded from above and away from $0$.
Fix a nondecreasing sequence $(\ell_0 (y_k^+))_{k \geq 1}$ with
$\ell_0 (y_k^+) > e^{\frac{8 \sqrt{\log b}}{3} (y_k^+)^{1/20}}$ for each $k \geq 1$.
By the boundedness of the sequence (\ref{eq:sequence-gamma}),
there exists a subsequence $(\ell_0 (y_{k_j}^+))_{j \geq 1}$ of $(\ell_0 (y_k^+))_{k \geq 1}$
such that the limit
\begin{equation} \label{eq:const-gamma-candidate}
\widetilde{\gamma_*} := \lim_{j \to \infty}
\int_{(\ell_0 (y_{k_j}^+))^{2/5}}^{\ell_0 (y_{k_j}^+)}
z e^{2 \sqrt{\log b}~z}
\mathbb{P} \left[\max_{u \in T_{\ell_0 (y_{k_j}^+)}} h_u
> \sqrt{\log b}~\ell_0 (y_{k_j}^+) + z \right] dz \in (0, \infty)
\end{equation}
exists.
We set
\begin{equation*}
\ell_j := \ell_0 (y_{k_j}^+),~~j \geq 1.
\end{equation*}
Note that by the definition of $\ell_0 (y_j^+)$,
for any $y_j$ with $y_j \leq y_j^+$, we have 
\begin{equation*}
\ell_j \geq \ell_0 (y_j^+) > e^{\frac{8 \sqrt{\log b}}{3} (y_j^+)^{1/20}}
\geq e^{\frac{8 \sqrt{\log b}}{3} (y_j)^{1/20}}.
\end{equation*}
Fix $j \geq 1$ and $y_j$ with $y_j^{-} \leq y_j \leq y_j^+$.
Recall definitions $J_{\ell_j}$ and $\Lambda_{y_j, \ell_j, J_{\ell_j}}^n (t)$
from (\ref{eq:limit-tail-interval}) and (\ref{eq:limit-tail-lambda-interval}).
Fix $v^* \in T_{n - \ell_j}$.
Let $\widetilde{L}^{\downarrow}$ be a local time of a Brownian motion on
$\widetilde{T}_{\leq \ell_j}^{v^*}$. 
We set
$\widetilde{\tau}^{\downarrow} (s) := \inf \{r \geq 0 : \widetilde{L}_r^{\downarrow} (v^*) > s \}$.
Recall $\widetilde{P}_{n - \ell_j}^B$ and $\widetilde{B}$ from
(\ref{eq:change-measure}) and (\ref{eq:bm-drift}).
We define $\phi (x)$ 
by replacing $y$, $\ell (y)$, and $v$ 
in the definition of $\psi (x)$ in (\ref{eq:lem-limit-tail-1-1-2-psi})
by $y_j$, $\ell_j$, and $v^{*}$, respectively.
By similar arguments to those in (\ref{eq:lem-limit-tail-1-1}) and (\ref{eq:lem-limit-tail-1-1-2}),
we have for $t \geq n \log n$,
\begin{align} \label{eq:prop-limit-tail-1}
&\widetilde{P}_{\rho} \left[
F_{v^*, y_j, \ell_j}^n (t) \cap \left\{\sqrt{\widetilde{L}_{\widetilde{\tau} (t)}^n (v^*)}
\in \sqrt{t} + \frac{a_n (t)}{n} (n - \ell_j) + J_{\ell_j} \right \} \right] \notag \\
&= \left(1 + O\left(\frac{1}{\log n} \right) \right)
\widetilde{E}_{n - \ell_j}^B \left[
1_{\left\{\widetilde{B}_s \leq y_j,
~\forall s \in [0,~ n - \ell_j],~
\widetilde{B}_{n - \ell_j} \in J_{\ell_j} \right \}} 
\phi \left(\widetilde{B}_{n - \ell_j} \right) \right] \notag \\
&~~~~- \left(1 + O\left(\frac{1}{\log n} \right) \right)
\widetilde{E}_{n - \ell_j}^B \left[
\begin{minipage}{180pt}
$1_{\left\{\widetilde{B}_s \leq y_j,
~\forall s \in [0,~ n - \ell_j],~
\widetilde{B}_{n - \ell_j} \in J_{\ell_j} \right \}}$ \\
$\times 1_{\left\{\widetilde{B}_s
< - (1 - \delta) \sqrt{t} - \frac{a_n (t)}{n} s,~\exists s \in [0,~ n - \ell_j] \right \}}
\phi \left(\widetilde{B}_{n - \ell_j} \right)$
\end{minipage}
\right] \notag \\
&=: K_1 - K_2,
\end{align}
where
we have used the fact that for $t \geq n \log n$,~
$\exp \left(- \frac{3}{8} \int_0^{n - \ell_j} \frac{ds}{X_s} \right) = 1 + O((\log n)^{-1})$
under the event that 
$\sqrt{X_s/2} \geq \delta \sqrt{t}$ for all $0 \leq s \leq n - \ell_j$.

We first estimate $K_2$. 
By (\ref{eq:joint-max-density}) and Proposition \ref{prop:tail}(i),
for all $\varepsilon > 0$ and $j \geq 1$, 
there exists $n_0(j) = n_0 (y_j^{-}, y_j^{+}, \ell_j) \in \mathbb{N}$ 
such that for all $n \geq n_0(j)$,
we have
\begin{equation} \label{eq:prop-limit-tail-2}
K_2
\leq c_1 b^{- (n - \ell_j)} \sqrt{\ell_j} n^{1 - (1 - \delta)^2 c_*} y_j e^{- 2 \sqrt{\log b}~y_j} 
\leq \varepsilon b^{- (n - \ell_j)} y_j e^{- 2 \sqrt{\log b}~y_j},
\end{equation}
uniformly in $y_j$ and $t$ satisfying (\ref{eq:prop-limit-tail-assumption-uniformity})
(we take $c_*$ large enough).
Next, we estimate $K_1$.
By the density (\ref{eq:density-barrier-standard}),
$K_1$ is equal to
\begin{align} \label{eq:prop-limit-tail-3}
&\int_{(\ell_j)^{2/5}}^{\ell_j} 
\frac{1 + O(1/\log n)}{\sqrt{\pi (n - \ell_j)}} \left(e^{- \frac{(z - y_j)^2}{n - \ell_j}} - e^{- \frac{(z + y_j)^2}{n - \ell_j}}
\right) \notag \\
&\times e^{\frac{2 a_n (t)}{n} z - \frac{2 a_n (t)}{n} y_j - \frac{(a_n (t))^2}{n^2} (n - \ell_j)}
\sqrt{\frac{\sqrt{t}}{\sqrt{t} + \frac{a_n (t)}{n} (n - \ell_j) + y_j - z}} 
~~\widetilde{P}_{v^{*}} [A_{y_j, \ell_j}^{n, t} (z)] dz
\notag \\
&= b^{- (n - \ell_j)} (1 + O(1/\log n)) 
\frac{4}{\sqrt{\pi}} y_j e^{- 2 \sqrt{\log b}~y_j}
\notag \\
&~~~~\times \int_{(\ell_j)^{2/5}}^{\ell_j}
z e^{2 \sqrt{\log b}~z}
\sqrt{\frac{\sqrt{t} + n}{\sqrt{t} + \frac{a_n (t)}{n} (n - \ell_j) + y_j - z}} 
~~P_{\rho} [B_{y_j, \ell_j}^{n, t} (z)]~dz,
\end{align}
where 
we have set
\begin{equation*}
A_{y_j, \ell_j}^{n, t} (z) := \left\{
\max_{u \in T_{\ell_j}^{v^*}} \sqrt{\widetilde{L}_{\widetilde{\tau}^{\downarrow}
(s_{y_j, \ell_j}^{n, t} (z))}^{\downarrow} (u)}
> \sqrt{t} + a_n (t) + y_j \right \},
\end{equation*}
\begin{equation*}
B_{y_j, \ell_j}^{n, t} (z) := \left \{
\frac{\max_{u \in T_{\ell_j}} 
L_{\tau (s_{y_j, \ell_j}^{n, t} (z))}^{\ell_j} (u) - s_{y_j, \ell_j}^{n, t} (z)}{2 \sqrt{s_{y_j, \ell_j}^{n, t} (z)}}
> \sqrt{\log b}~\ell_j + z + \Delta_{y_j, \ell_j}^{n, t} (z) \right \}.
\end{equation*}
Here,
\begin{equation*}
s_{y_j, \ell_j}^{n, t} (z) := \left(\sqrt{t} + \frac{a_n (t)}{n} (n - \ell_j) + y_j - z \right)^2
\end{equation*}
and $\Delta_{y_j, \ell_j}^{n, t} (z)$ is the remainder term so that
$\widetilde{P}_{v^{*}} [A_{y_j, \ell_j}^{n, t} (z)] = P_{\rho} [B_{y_j, \ell_j}^{n, t} (z)]$.
(Note that by the symmetry of the $b$-ary tree, 
the law of the Brownian motion on
$\widetilde{T}_{\leq \ell_j}^{v^*}$ starting at $v^*$ is the same as that of the Brownian motion on
$\widetilde{T}_{\leq \ell_j}$ starting at $\rho$.)
One can show that for each $j \geq 1$,
\begin{equation} \label{eq:prop-limit-tail-remainder-term}
\lim_{n \to \infty} \Delta_{y_j, \ell_j}^{n, t} (z) = 0,~~~\text{uniformly in}~
z \in [(\ell_j)^{2/5},~\ell_j]~\text{and}~y_j, t
~\text{satisfying}~(\ref{eq:prop-limit-tail-assumption-uniformity}).
\end{equation}
By Theorem \ref{thm:iso-thm},
we have for any fixed $m \in \mathbb{N}$
\begin{equation*} 
\left(\frac{L_{\tau (s)}^m (v) - s}{2 \sqrt{s}} \right)_{v \in T_{m}}
\longrightarrow
(h_v)_{v \in T_m}
~~~\text{in law as}~s \to \infty,
\end{equation*}
where $(h_v)_{v \in T}$ is a BRW.
By this together with
the definition of $\widetilde{\gamma}_*$ in (\ref{eq:const-gamma-candidate})
and (\ref{eq:prop-limit-tail-remainder-term}),
for all $\varepsilon > 0$, there exists $j_0 \in \mathbb{N}$ such that the following holds:
for each $j \geq j_0$, there exists $n_0 (j) = n_0 (y_j^{-}, y_j^{+}, \ell_j) \in \mathbb{N}$
such that for all $n \geq n_0 (j)$,
\begin{equation} \label{eq:prop-limit-tail-5}
\left|\int_{(\ell_j)^{2/5}}^{\ell_j}
P_{\rho} \left[B_{y_j, \ell_j}^{n, t} (z) \right]
\cdot z e^{2 \sqrt{\log b}~z} dz 
~- ~\widetilde{\gamma}_*~ \right| < \varepsilon,
\end{equation}
uniformly in $y_j$ and $t$ satisfying (\ref{eq:prop-limit-tail-assumption-uniformity}).
Thus, by (\ref{eq:prop-limit-tail-1})-(\ref{eq:prop-limit-tail-3}), 
(\ref{eq:prop-limit-tail-5}), and the definition (\ref{eq:const-beta}) of $\beta_*$,
for all $\varepsilon > 0$, there exists $j_0 \in \mathbb{N}$ such that the following holds:
for each $j \geq j_0$, there exists 
$n_0 (j) = n_0 (y_j^{-}, y_j^{+}, \ell_j) \in \mathbb{N}$ such that for all $n \geq n_0 (j)$,
\begin{equation} \label{eq:prop-limit-tail-6}
\left|\widetilde{E}_{\rho} \left[\Lambda_{y_j, \ell_j, J_{\ell_j}}^n (t) \right]
(y_j)^{-1} e^{2 \sqrt{\log b}~y_j}
- \frac{4}{\sqrt{\pi}} \beta_* \widetilde{\gamma}_* \right| < \varepsilon,
\end{equation}
uniformly in $y_j$ and $t$ satisfying (\ref{eq:prop-limit-tail-assumption-uniformity}).
By Lemma \ref{lem:limit-tail-3}, \ref{lem:limit-tail-4}, and (\ref{eq:prop-limit-tail-6}),
for all $\varepsilon > 0$, there exists $j_0 \in \mathbb{N}$ such that the following holds:
for each $j \geq j_0$, there exists 
$n_0 (j) = n_0(y_j^{-}, y_j^{+}, \ell_j) \in \mathbb{N}$ such that for all $n \geq n_0 (j)$,
(\ref{eq:prop-limit-tail-statement}) holds by replacing $\gamma_*$ with $\widetilde{\gamma}_*$
uniformly in $y_j$ and $t$ satisfying (\ref{eq:prop-limit-tail-assumption-uniformity})
(we take $c_* > 0$ large enough).
Let $\widehat{\gamma}_*$ be the limit of any convergent subsequence
of (\ref{eq:sequence-gamma}).
By taking a sub-subsequence, if necessary, and repeating the above argument,
we have (\ref{eq:prop-limit-tail-statement}) if we replace $\gamma_*$ with
$\widehat{\gamma}_*$.
Thus, the full sequence (\ref{eq:sequence-gamma}) converges
to a finite positive constant and we write $\gamma_*$ to denote the limit.
Therefore, we have (\ref{eq:prop-limit-tail-statement}).~~~$\Box$

\section{Proof of Theorem \ref{thm:convergence-point-process} 
and Corollary \ref{cor:convergence-in-law}}
\label{sec:pf-convergence-point-process}
In this section, we prove Theorem \ref{thm:convergence-point-process} and Corollary 
\ref{cor:convergence-in-law}.
We begin with preliminary lemmas.
Let $(h_v)_{v \in T}$ be a BRW on $T$ defined in Section \ref{sec:intro}.
For each $n \in \mathbb{N}$, we set
\begin{equation*} 
D_n^{(2)} := \sum_{v \in T_n} \left(\sqrt{\log b}~n - h_v \right)^2
e^{- 4 \sqrt{\log b}~\left(\sqrt{\log b}~ n - h_v \right)}.
\end{equation*} 
Then the following holds:
\begin{lem} \label{lem:control-squared-dm}
For all $\varepsilon > 0$,
\begin{equation} \label{eq:control-squared-dm-statement}
\lim_{n \to \infty} \mathbb{P} \left(D_n^{(2)} \geq \varepsilon \right) = 0.
\end{equation}
\end{lem}
{\it Proof.}
Set $m_n := \sqrt{\log b}~n - \frac{3}{4 \sqrt{\log b}} \log n$.
Fix any $y > 0$.
By the equality
$1_{\{\max_{v \in T_n} h_v \leq m_n + y \}} + 1_{\{\max_{v \in T_n} h_v > m_n + y \}} = 1$,
we bound the probability $\mathbb{P} (D_n^{(2)} \geq \varepsilon)$
from above by
\begin{align} \label{eq:control-squared-dm-1}
&\mathbb{P} \left(
\sum_{v \in T_n} \left(\sqrt{\log b}~n - h_v \right)^2
e^{- 4 \sqrt{\log b}~\left(\sqrt{\log b}~ n - h_v \right)}
1_{\left\{h_v \leq m_n + y \right \}}
\geq \varepsilon \right) \notag \\
&+ \mathbb{P} \left(\max_{v \in T_n} h_v > m_n + y \right).
\end{align}

Since $h_v$ is a Gaussian random variable with mean $0$ and variance $n/2$
for each $v \in T_n$,
a simple calculation implies that
\begin{equation*} 
\lim_{n \to \infty} \mathbb{E} \left[
\sum_{v \in T_n} \left(\sqrt{\log b}~n - h_v \right)^2
e^{- 4 \sqrt{\log b}~\left(\sqrt{\log b}~ n - h_v \right)}
1_{\left\{h_v \leq m_n + y \right \}}
\right] = 0.
\end{equation*}
By this,
Lemma \ref{lem:tail-max-brw-1}(i),
and (\ref{eq:control-squared-dm-1}),
we have
\begin{equation*}
\limsup_{n \to \infty} \mathbb{P} \left(D_n^{(2)} \geq \varepsilon \right)
\leq c_1 (1 + y) e^{- 2 \sqrt{\log b}~y}.
\end{equation*}
Since we can take arbitrary $y \in (0,~\infty)$, this implies (\ref{eq:control-squared-dm-statement}).
~~~$\Box$
\\ \\
Recall the definition of $\sigma (\cdot)$ from (\ref{eq:location-point}).
For $a \in \mathbb{R}$, an interval $I \subset \mathbb{R}$, $t > 0$, and $r, n \in \mathbb{N}$
with $r < n$, we set
\begin{equation} \label{eq:event-thin-region-extrema}
\mathcal{B}_{r, t, n}^{\sigma} (I; a)
:= \left\{\sqrt{L_{\tau (t)}^n (v)} - \sqrt{t} - a_n (t) \notin I,~
\forall v \in T_n
~\text{with}~\sigma (v_r) \in \left[a - b^{- r},~a \right] \right \}.
\end{equation}
Then we have the following:
\begin{lem} \label{lem:control-thin-region-extrema}
There exists $c_1 > 0$ such that the following holds:
for any finite interval $I \subset \mathbb{R}$, there exist $c_2 (I) > 0$
and $r_0 (I) \in \mathbb{N}$ such that
for all $t > 0$, $a \in \mathbb{R}$, $r \geq r_0 (I)$, and $n > r$,
\begin{equation} \label{eq:control-thin-region-extrema-statement}
P_{\rho} \left[\left(\mathcal{B}_{r, t, n}^{\sigma} (I; a) \right)^{\text{c}} \right]
\leq c_2 (I) r^3 e^{- c_1 r}.
\end{equation}
\end{lem}
{\it Proof.}
Fix any finite interval $I \subset \mathbb{R}$.
Let $\{\overline{I}, \underline{I} \}$ be the boundary of $I$ with $\underline{I} < \overline{I}$.
Fix any $t > 0$, $a \in \mathbb{R}$, $r \geq r_0(I)$, and $n > r$,
where we take $r_0 (I) > |\overline{I}|$ large enough. 
Recall the event $G_r^n (t)$ from (\ref{eq:event-broken-barrier}).
$P_{\rho} [\left(\mathcal{B}_{r, t, n}^{\sigma} (I; a) \right)^{\text{c}}]$
is bounded from above by
\begin{equation} \label{eq:control-thin-region-extrema-1}
\widetilde{P}_{\rho} \left[
\begin{minipage}{230pt}
$\exists v \in T_n~\text{with}~\sigma (v_r) \in \left[a - b^{-r},~a \right]~\text{s.t.}~\forall s \in [0,~ n]$, \\
$\sqrt{\widetilde{L}_{\widetilde{\tau} (t)}^n (v_s)} \leq \sqrt{t} + \frac{a_n (t)}{n} s
+ \kappa (\log (s \wedge (n-s) ))_{+} + r + 1$, \\
$\sqrt{\widetilde{L}_{\widetilde{\tau} (t)}^n (v)} - \sqrt{t} - a_n (t) \in I$
\end{minipage}
\right] 
+ \widetilde{P}_{\rho} \left(G_r^n (t) \right).
\end{equation}
Next, we estimate the number of $v \in T_n$ satisfying 
$\sigma \left(v_r \right) \in \left[a - b^{-r},~a \right]$.
We may assume $a \in [0, 1]$.
We have a sequence $(x_i)_{i \geq 1}$ with $x_i \in \{0, \dotsc, b - 1 \}$
such that
$a = \sum_{i = 1}^{\infty} \frac{x_i}{b^i}$.
In particular, we have
$\sum_{i = 1}^r \frac{\overline{v}_i - x_i}{b^i} \in [- b^{-r},~b^{-r}]$
for all $v \in T_n$ with $\sigma \left(v_r \right) \in \left[a - b^{-r},~a \right]$
and the label $(\overline{v}_1, \dotsc, \overline{v}_n)$.
By this and a simple observation, one can see that
$\{v \in T_n : \sigma (v_r) \in \left[a - b^{- r},~a \right] \}$ is a subset of
\begin{align*}
&\left\{v \in T_n : \overline{v}_i = x_i,~~1 \leq \forall i \leq r \right\} \notag \\
&\cup \bigcup_{i=1}^r
\left\{v \in T_n : \overline{v}_j = x_j,~~1 \leq \forall j \leq i-1,
~\overline{v}_i - x_i > 0, ~\overline{v}_{j^{\prime}} = x_{j^{\prime}} - (b-1),
~i < \forall j^{\prime} \leq r \right\} \notag \\
&\cup \bigcup_{i=1}^r
\left\{v \in T_n : \overline{v}_j = x_j,~~1 \leq \forall j \leq i-1,
~\overline{v}_i - x_i < 0,~\overline{v}_{j^{\prime}} = x_{j^{\prime}} + (b-1),
~i < \forall j^{\prime} \leq r \right\}.
\end{align*}
This implies 
\begin{equation} \label{eq:control-thin-region-extrema-2}
\left|\left\{v \in T_n : \sigma (v_r) \in \left[a - b^{- r},~a \right] \right \} \right|
\leq c_1 r b^{n-r}.
\end{equation}
By (\ref{eq:control-thin-region-extrema-2}) and
similar arguments to those in (\ref{eq:prop-upper-tail-2}) and (\ref{eq:prop-upper-tail-2-2}), 
the first term on the right-hand side of (\ref{eq:control-thin-region-extrema-1})
is bounded from above by
$c_2 (I) b^{-r} r^3$.
By this and
Lemma \ref{lem:broken-barrier},
we have (\ref{eq:control-thin-region-extrema-statement}).
~~~$\Box$
\\ \\
{\it Proof of Theorem \ref{thm:convergence-point-process}.}
Recall definitions of 
$Z_{\infty}$, $\Xi_{n, t}^{(m)}$, $\beta_*$, and $\gamma_*$ from 
(\ref{eq:critical-mandelbrot-cascade}), (\ref{eq:point-process}),
(\ref{eq:const-beta}), and (\ref{eq:const-gamma}).
For each interval $I \subset \mathbb{R}$, we will write $\partial I$ to denote its boundary. 
Fix any sequence of positive integers $(r_n)_{n \geq 1}$ with
$\lim_{n \to \infty} r_n = \infty$ and $\lim_{n \to \infty} r_n/n \in [0, 1)$.
Fix any $(t_n)_{n \geq 1}$ with
$\lim_{n \to \infty} \sqrt{t_n}/n = \theta \in [0, \infty]$
and $t_n \geq c_* n \log n$,
where $c_*$ is a sufficiently large positive constant.
By, for example, \cite[Proposition 11.1.VIII]{DJ},
in order to show the convergence of
the point process $\Xi_{n, t_n}^{(r_n)}$ to the Cox process
(\ref{eq:limit-cox-process}) as $n \to \infty$,
it is enough to prove the following:
for all finite disjoint intervals $A_i := (\underline{A}_i,~\overline{A}_i] \subset [0, 1],~1 \leq i \leq m$
with $\mathbb{P} \left(Z_{\infty} \left(\partial A_i \right) = 0 \right) = 1$
for all $1 \leq i \leq m$,
finite intervals $B_i := (\underline{B}_i,~\overline{B}_i],~1 \leq i \leq m$,
and positive values $a_1, \dotsc, a_m$,
\begin{align} \label{eq:pf-convergence-point-process-1}
&~~~~\lim_{n \to \infty} E_{\rho} \left[
\exp \left\{- \sum_{i=1}^m a_i \Xi_{n, t_n}^{(r_n)} \left(A_i \times B_i \right) \right \} \right] \notag \\
&= \mathbb{E} \left[
\exp\left\{- \frac{4}{\sqrt{\pi}} \beta_* \gamma_*
\sum_{i=1}^m \left(1 - e^{- a_i} \right) Z_{\infty} (A_i)
\left(e^{- 2\sqrt{\log b}~\underline{B}_i} - e^{- 2 \sqrt{\log b}~\overline{B}_i} \right) \right \} \right].
\end{align}
Let $q < n$ be a positive integer.
To show (\ref{eq:pf-convergence-point-process-1}),
we first prove convergence of $\Xi_{n, t_n}^{(n - q)}$
as $n \to \infty$, $q \to \infty$.
Recall the events (\ref{eq:event-thin-region-extrema}).
We have
\begin{align} \label{eq:pf-convergence-point-process-2}
&~~~~E_{\rho} \left[
\exp \left\{- \sum_{i=1}^m a_i \Xi_{n, t_n}^{(n-q)} \left(A_i \times B_i \right) \right \} \right] \notag \\
&= E_{\rho} \left[
\exp \left\{- \sum_{i=1}^m a_i \Xi_{n, t_n}^{(n-q)} \left(A_i \times B_i \right) \right \}
 1_{\bigcap_{i=1}^m \mathcal{B}_{q, t_n, n}^{\sigma} \left(B_i; \underline{A}_i \right) 
\cap \mathcal{B}_{q, t_n, n}^{\sigma} \left(B_i; \overline{A}_i \right) } \right] \notag \\
&~~~~+ E_{\rho} \left[
\exp \left\{- \sum_{i=1}^m a_i \Xi_{n, t_n}^{(n-q)} \left(A_i \times B_i \right) \right \}
 1_{\bigcup_{i=1}^m \left(\mathcal{B}_{q, t_n, n}^{\sigma} \left(B_i; \underline{A}_i \right) \right)^{\text{c}}
\cup \left(\mathcal{B}_{q, t_n, n}^{\sigma} \left(B_i; \overline{A}_i \right) \right)^{\text{c}}
 } \right] \notag \\
&=: J_1 + J_2.
\end{align}

We estimate $J_1$ in (\ref{eq:pf-convergence-point-process-2}).
Under the event 
$\bigcap_{i=1}^m \mathcal{B}_{q, t_n, n}^{\sigma} \left(B_i; \overline{A}_i \right)
\cap \mathcal{B}_{q, t_n, n}^{\sigma} \left(B_i; \underline{A}_i \right)$, we have for all $v \in T_q$
and $1 \leq i \leq m$,
\begin{align} \label{eq:pf-convergence-point-process-4}
&~~~~1_{\left \{\sigma \left(\text{arg} \max_{v} L_{\tau (t_n)}^n \right)
 \in A_i,~~\max_{u \in T_{n-q}^v} \sqrt{L_{\tau (t_n)}^n (u)} - \sqrt{t_n} - a_n (t_n) \in B_i \right \}}
\notag \\
&= 1_{\left \{\sigma (v)
 \in A_i,~~\max_{u \in T_{n-q}^v} \sqrt{L_{\tau (t_n)}^n (u)} - \sqrt{t_n} - a_n (t_n) \in B_i \right \}}.
\end{align}
Thus, we have
\begin{align} \label{eq:pf-convergence-point-process-5}
J_1 &= E_{\rho} \left[
\exp \left\{- \sum_{v \in T_{q}} \sum_{i=1}^m a_i 
 1_{\left \{\sigma (v) \in A_i,~~
\max_{u \in T_{n-q}^v} \sqrt{L_{\tau (t_n)}^n (u)} - \sqrt{t_n} - a_n (t_n) \in B_i
 \right \}} \right \} \right] \notag \\
&~~~- E_{\rho} \left[
\begin{minipage}{250pt}
$\exp \left\{- \displaystyle \sum_{v \in T_{q}} \sum_{i=1}^m a_i 
 1_{\left \{\sigma (v) \in A_i,~~
\max_{u \in T_{n-q}^v} \sqrt{L_{\tau (t_n)}^n (u)} - \sqrt{t_n} - a_n (t_n) \in B_i \right \}} \right \}$ \\
$1_{\bigcup_{i=1}^m \left(\mathcal{B}_{q, t_n, n}^{\sigma} \left(B_i; \underline{A}_i \right) \right)^{\text{c}}
\cup \left(\mathcal{B}_{q, t_n, n}^{\sigma} \left(B_i; \overline{A}_i \right) \right)^{\text{c}}}$
\end{minipage}
\right] \notag \\
&=: J_{1, 1} - J_{1, 2}.
\end{align}
By Lemma \ref{lem:control-thin-region-extrema}, we have
\begin{equation} \label{eq:pf-convergence-point-process-6}
\max \{J_2, J_{1, 2} \} \leq \sum_{i=1}^m c_1 (B_i) q^3 e^{- c_2 q}.
\end{equation}

We estimate $J_{1, 1}$.
By Theorem \ref{thm:iso-thm}, on the same probability space
(we will write $P$ to denote the probability measure),
we can construct a local time $(L_{\tau (t_n)}^n (v))_{v \in T_{\leq n}}$
and two BRWs $(h_v)_{v \in T_{\leq n}}$, $(h_v^{\prime})_{v \in T_{\leq n}}$
satisfying (\ref{eq:iso-thm-1}) and (\ref{eq:iso-thm-2}).
Fix $\delta \in (0, 1/3)$. We set
\begin{equation} \label{eq:pf-convergence-point-process-7}
C_{q} := \left\{|h_v|, |h_v^{\prime}|
\leq \sqrt{\log b}~q - \frac{3}{4 \sqrt{\log b}} (1 - \delta) \log q,
~\forall v \in T_q \right \}.
\end{equation}
For each $v \in T_{q}$, let $\widetilde{L}^{\downarrow}$ be a local time of a Brownian motion
on $\widetilde{T}_{\leq n-q}^v$ and set
$\widetilde{\tau}^{\downarrow} (p) := \inf \{s \geq 0 : \widetilde{L}_s^{\downarrow} (v) > p \}$.
We omit the subscript $v$ in $\widetilde{L}^{\downarrow}$ and $\widetilde{\tau}^{\downarrow}$.
By Lemma \ref{lem:markov-local-time}, we have
\begin{equation} \label{eq:pf-convergence-point-process-8}
J_{1, 1} 
= E \left[1_{C_{q}} \prod_{v \in T_{q}} K_v \right] 
+ E \left[1_{\left(C_{q}\right)^{\text{c}}} \prod_{v \in T_{q}} K_v \right] 
=: J_{1, 1, 1} + J_{1, 1, 2},
\end{equation}
where for each $v \in T_{q}$, we have set
\begin{equation} \label{eq:pf-convergence-point-process-9}
K_v := \widetilde{E}_v \left[
\exp \left\{- \sum_{i=1}^m a_i 
 1_{\left \{\sigma (v) \in A_i,~
\max_{u \in T_{n-q}^v} \sqrt{\widetilde{L}_{\widetilde{\tau}^{\downarrow} 
\left(L_{\tau (t_n)}^n (v) \right)}^{\downarrow} (u)}
- \sqrt{t_n} - a_n (t_n) \in B_i \right \}} \right \} \right].
\end{equation}
By Lemma \ref{lem:tail-max-brw-1}, we have
\begin{equation} \label{eq:pf-convergence-point-process-10}
J_{1, 1, 2} \leq c_3 (\log q) \cdot (q)^{- \frac{3}{2} \delta}.
\end{equation} 

We estimate $J_{1, 1, 1}$.
Fix $v \in T_q$.
For each $1 \leq i \leq m$, we set
$$E_v^{(i)} (n, t_n, q) :=
\left \{
\max_{u \in T_{n-q}^v} \sqrt{\widetilde{L}_{\widetilde{\tau}^{\downarrow} 
\left(L_{\tau (t_n)}^n (v) \right)}^{\downarrow} (u)}
- \sqrt{t_n} - a_n (t_n) \in B_i \right \}.$$
Since $A_1, \dotsc, A_m$ are disjoint, $K_v$ is equal to
\begin{align} \label{eq:pf-convergence-point-process-11}
&~~~\widetilde{E}_v \left[\prod_{i=1}^m \left\{1 - (1 - e^{- a_i}) 1_{\{\sigma (v) \in A_i \}}
1_{E_v^{(i)} (n, t_n, q)} \right \} \right] \notag \\
&= \widetilde{E}_v \left[1 - \sum_{i=1}^m (1 - e^{- a_i}) 1_{\{\sigma (v) \in A_i \}} 1_{E_v^{(i)} (n, t_n, q)} \right] 
\notag \\
&=\exp \left\{ \log \left(
1 - \sum_{i = 1}^m \left(1 - e^{- a_i} \right) 1_{\left\{\sigma (v) \in A_i \right\}}
\widetilde{P}_v [E_v^{(i)} (n, t_n, q)] \right) \right \}.
\end{align}
On the event $C_{q}$, by (\ref{eq:iso-thm-2}), we have for all $v \in T_q$
\begin{equation} \label{eq:pf-convergence-point-process-12}
\sqrt{t_n} + a_n (t_n)
= \sqrt{L_{\tau (t_n)}^n (v)} + a_{n-q} \left(L_{\tau (t_n)}^n (v) \right)
+ \sqrt{\log b}~q - h_v^{\prime} + \delta_{v, q}^n,
\end{equation}
where
\begin{equation} \label{eq:pf-convergence-point-process-13}
\delta_{v, q}^n := \frac{3\log \left(1- q/n \right)}{4 \sqrt{\log b}} 
+ \frac{\log \left(1 + O \left(q/\sqrt{t_n} \right) \right)}{4 \sqrt{\log b}}
+ O \left(q^2/\sqrt{t_n} \right) .
\end{equation}
We take any sufficiently small $\varepsilon > 0$
and sufficiently large $q_0 \in \mathbb{N}$ which depends on $B_1, \dotsc, B_m$ and $\varepsilon$.
We assume that $q \geq q_0$ and $n \geq n_0$, 
where we take sufficiently large $n_0 = n_0(q, \varepsilon) \in \mathbb{N}$. 
By Proposition \ref{prop:limit-tail} and 
(\ref{eq:pf-convergence-point-process-12})-(\ref{eq:pf-convergence-point-process-13}),
under the event $C_{q}$, we have for all $1 \leq i \leq m$ and $v \in T_{q}$,
$\widetilde{P}_v [E_v^{(i)} (n, t_n, q)]$ is bounded from below by
\begin{align} \label{eq:pf-convergence-point-process-14}
&\widetilde{P}_v \left[
\begin{minipage}{245pt}
$\max_{u \in T_{n-q}^v} \sqrt{\widetilde{L}_{\widetilde{\tau}^{\downarrow} 
\left(L_{\tau (t_n)}^n (v) \right)}^{\downarrow} (u)}
- \sqrt{L_{\tau (t_n)}^n (v)} - a_{n - q} \left(L_{\tau (t_n)}^n (v) \right)$ \\
$\geq \sqrt{\log b}~q - h_v^{\prime} + \underline{B}_i 
+ \delta_{v, q}^n$
\end{minipage}
\right] \notag \\
&- \widetilde{P}_v \left[
\begin{minipage}{245pt}
$\max_{u \in T_{n-q}^v} \sqrt{\widetilde{L}_{\widetilde{\tau}^{\downarrow} 
\left(L_{\tau (t_n)}^n (v) \right)}^{\downarrow} (u)}
- \sqrt{L_{\tau (t_n)}^n (v)} - a_{n - q} \left(L_{\tau (t_n)}^n (v) \right)$ \\
$\geq \sqrt{\log b}~q - h_v^{\prime} + \overline{B}_i 
+ \delta_{v, q}^n$
\end{minipage}
\right] \notag \\
&\geq (1 - c_4 \varepsilon) \frac{4}{\sqrt{\pi}} \beta_* \gamma_*
\left(\sqrt{\log b}~q - h_v^{\prime} + \underline{B}_i \right)
e^{- 2 \sqrt{\log b}~\left(\sqrt{\log b}~q - h_v^{\prime} + \underline{B}_i \right)}
\notag \\
&~~~- (1 + c_4 \varepsilon) \frac{4}{\sqrt{\pi}} \beta_* \gamma_*
\left(\sqrt{\log b}~q - h_v^{\prime} + \overline{B}_i \right)
e^{- 2 \sqrt{\log b}~\left(\sqrt{\log b}~q - h_v^{\prime} + \overline{B}_i \right)}.
\end{align}
We set
\begin{equation*}
D_{q}^{(2)} := \sum_{v \in T_{q}} \left(\sqrt{\log b}~q - h_v^{\prime} \right)^2
e^{- 4 \sqrt{\log b}~\left(\sqrt{\log b}~q - h_v^{\prime} \right)},
\end{equation*}
\begin{equation*}
W_{q} := \sum_{v \in T_{q}} 
e^{- 2 \sqrt{\log b}~\left(\sqrt{\log b}~q - h_v^{\prime} \right)}.
\end{equation*}
Recall the random measure $Z_q$ from (\ref{eq:critical-cascade-measure}).
By (\ref{eq:pf-convergence-point-process-11}), (\ref{eq:pf-convergence-point-process-14}),
and Taylor's expansion of the function $x \mapsto \log (1-x)$,
under the event $C_{q} \cap \left\{D_{q}^{(2)} < \varepsilon \right \}
\cap \left \{W_{q} < \varepsilon \right \}$, 
$\prod_{v \in T_{q}} K_v$ is bounded from above by 
\begin{align} \label{eq:pf-convergence-point-process-15}
&e^{c_5 (B_1, \dotsc, B_m) \varepsilon} 
\cdot \exp \left \{- \frac{4}{\sqrt{\pi}} \beta_* \gamma_* \sum_{i=1}^m \left(1-e^{-a_i} \right)
Z_q (A_i)
\cdot \left(e^{- 2 \sqrt{\log b}~\underline{B}_i} - e^{- 2 \sqrt{\log b}~\overline{B}_i} \right) \right \} \notag \\
&\times \exp \left \{c_4 \varepsilon \frac{4}{\sqrt{\pi}} \beta_* \gamma_* \sum_{i=1}^m \left(1-e^{-a_i}
 \right)
Z_q (A_i) 
\left(e^{- 2 \sqrt{\log b}~\underline{B}_i} + e^{- 2 \sqrt{\log b}~\overline{B}_i} \right) \right \}.
\end{align}
We can obtain a similar lower bound of $\prod_{v \in T_{q}} K_v$.
By Lemma \ref{lem:tail-max-brw-1}, Lemma \ref{lem:control-squared-dm}, and the fact
that $\lim_{q \to \infty} W_{q} = 0$ almost surely (see \cite{Ly}), we have
\begin{equation} \label{eq:pf-convergence-point-process-16}
\lim_{q \to \infty}
P\left(\left(C_{q} \right)^{\text{c}} \cup \left\{D_{q}^{(2)} \geq \varepsilon \right \}
\cup \left\{W_{q} \geq \varepsilon \right \} \right) = 0.
\end{equation}
Thus, by the above estimates,
taking $n \to \infty$, then $q \to \infty$,
and finally $\varepsilon \to 0$,
we have 
\begin{equation} \label{eq:pf-convergence-point-process-17}
\left| E_{\rho} \left[
e^{- \sum_{i=1}^m a_i \Xi_{n, t_n}^{(n-q)} \left(A_i \times B_i \right) } \right] 
- \mathbb{E} \left[
e^{- \frac{4}{\sqrt{\pi}} \beta_* \gamma_*
\sum_{i=1}^m \left(1 - e^{- a_i} \right) Z_{\infty} (A_i)
\left(e^{- 2\sqrt{\log b}~\underline{B}_i} - e^{- 2 \sqrt{\log b}~\overline{B}_i} \right) } \right]
\right| \to 0.
\end{equation}

Next, by using (\ref{eq:pf-convergence-point-process-17}), we will prove
(\ref{eq:pf-convergence-point-process-1}).
Let $z_*$ be a real number with $z_* < \min_{1 \leq i \leq m} \underline{B}_i$.
Take $q_0 = q_0 (z_*) \in \mathbb{N}$ large enough
and fix any $q \geq q_0$.
Take $n \in \mathbb{N}$
large enough so that
$q < n-r_n < n - q$ and $q < n/4$.
We set
\begin{align*}
&U_{z_*, q}^n (t_n) :=
\left \{
\begin{minipage}{180pt}
$\exists v, u \in T_n
~\text{with}~q \leq |v \wedge u| \leq n-q~~\text{s.t.}$ \\
$\sqrt{L_{\tau (t_n)}^n (v)},
~\sqrt{L_{\tau (t_n)}^n (u)} \geq \sqrt{t_n} + a_n (t_n) + z_*$
\end{minipage}
\right \}.
\end{align*}
Under the event $\left(U_{z_*, q}^n (t_n) \right)^{\text{c}}$, 
we have
\begin{align} \label{eq:pf-convergence-point-process-19}
&~~~\left\{\text{arg} \max_v L_{\tau(t_n)}^n : 
v \in T_{n-r_n},~\max_{u \in T_{r_n}^v} \sqrt{L_{\tau (t_n)}^n (u)}
\geq \sqrt{t_n} + a_n (t_n) + z_* \right \} \notag \\
&= \left\{\text{arg} \max_v L_{\tau(t_n)}^n : 
v \in T_q,~\max_{u \in T_{n-q}^v} \sqrt{L_{\tau (t_n)}^n (u)}
\geq \sqrt{t_n} + a_n (t_n) + z_* \right \}.
\end{align}
By (\ref{eq:pf-convergence-point-process-17}),
(\ref{eq:pf-convergence-point-process-19}), 
and Proposition \ref{prop:geometry}, 
we have
(\ref{eq:pf-convergence-point-process-1}).
$\Box$ \\\\
{\em Proof of Corollary \ref{cor:convergence-in-law}.}
Corollary \ref{cor:convergence-in-law} immediately follows form
Theorem \ref{thm:convergence-point-process} and Proposition \ref{prop:tail}(i).
We omit the details. $\Box$

\appendix \section{Appendix} \label{sec:appendix}
In this section, we give proof of some technical estimates.

\subsection{Gaussian process associated with Brownian motion on metric tree} \label{subsec:pf-cov-gff}
We adopt the notation in Section \ref{sec:metric-tree}.
Let $\widetilde{X} = (\widetilde{X}_t, t \geq 0, \widetilde{P}_{x}, x \in \widetilde{T}_{\leq n})$
be a Brownian motion on $\widetilde{T}_{\leq n}$.
Let $\{\widetilde{h}_x : x \in \widetilde{T}_{\leq n} \}$ be a centered Gaussian process
with $\mathbb{E} (\widetilde{h}_x \widetilde{h}_y) 
= \widetilde{E}_x (\widetilde{L}_{H_{\rho}}^n (y))$ for all $x, y \in \widetilde{T}_{\leq n}$,
where $H_x := \inf \{t \ge 0 : \widetilde{X}_t = x \}$.
We will call $\widetilde{h}$ a Gaussian process associated with the Brownian motion $\widetilde{X}$.
We have an explicit representation of the covariance of $\widetilde{h}$.
\begin{lem} \label{lem:covariance-gff}
For all $x, y \in \widetilde{T}_{\leq n}$,
\begin{equation*}
\mathbb{E} (\widetilde{h}_x \widetilde{h}_y) 
= d(\rho, x) + d(\rho, y) - d(x, y).
\end{equation*}
In particular,
the law of $\{\widetilde{h}_{v_s} : 0 \leq s \leq n  \}$
is the same as that of a standard Brownian motion
$\{W_s : 0 \leq s \leq n \}$ on $\mathbb{R},$
where $v_s$ is the point on the unique path in $\widetilde{T}_{\leq n}$ from $\rho$ to $v$
with $d(\rho, v_s) = s/2$.
\end{lem}
{\em Proof.}
Let $\{z_1, z_2 \}, \{w_1, w_2 \}$ be edges of $T$ with
$x \in I_{\{z_1, z_2 \}}$ and $y \in I_{\{w_1, w_2 \}}$.
Let $\theta_s, s \geq 0$ be shift operators.
We have
\begin{align} \label{eq:covariance-gff-1}
&\mathbb{E} (\widetilde{h}_x \widetilde{h}_y) = \widetilde{E}_x (\widetilde{L}_{H_{\rho}}^n (y)) \notag \\
&= \widetilde{E}_x \left(\widetilde{L}_{H_{z_1} \wedge H_{z_2}}^n (y) 
+ \widetilde{L}_{H_{\rho}}^n (y) \circ \theta_{H_{z_1} \wedge H_{z_2}} \right) \notag \\
&= \widetilde{E}_x (\widetilde{L}_{H_{z_1} \wedge H_{z_2}}^n (y))
+ \widetilde{P}_x (H_{z_1} < H_{z_2}) \widetilde{E}_{z_1} (\widetilde{L}_{H_{\rho}}^n (y))
+ \widetilde{P}_x (H_{z_2} < H_{z_1}) \widetilde{E}_{z_2} (\widetilde{L}_{H_{\rho}}^n (y)).
\end{align}
Note that we have
$\widetilde{E}_{w} (\widetilde{L}_{H_{\rho}}^n (w^{\prime})) = \widetilde{E}_{w^{\prime}} (\widetilde{L}_{H_{\rho}}^n (w))$,
$w, w^{\prime} \in \widetilde{T}_{\leq n}$;
see, for example, \cite[Lemma 5.1]{EKMRS}.
Since $\widetilde{X}$ behaves like a standard Brownian motion on each $I_{\{w, w^{\prime} \}}$
($\{w, w^{\prime}\}$ is an edge of $T$),
we have $\widetilde{P}_x (H_{w} < H_{w^{\prime}}) = d(x, w^{\prime})/d(w, w^{\prime})$, 
$x \in I_{\{w, w^{\prime}\}}$; see, for example, \cite[Chapter 2, Exercise 8.13]{KS}.
By these, the right of (\ref{eq:covariance-gff-1}) is equal to
\begin{equation} \label{eq:covariance-gff-1-2}
\widetilde{E}_x (\widetilde{L}_{H_{z_1} \wedge H_{z_2}}^n (y))
+ \frac{d(x, z_2)}{d(z_1, z_2)} \widetilde{E}_{y} (\widetilde{L}_{H_{\rho}}^n (z_1)) 
+ \frac{d(x, z_1)}{d(z_1, z_2)} \widetilde{E}_{y} (\widetilde{L}_{H_{\rho}}^n (z_2)).
\end{equation}
Similar arguments imply that for each $i \in \{1, 2\}$,
\begin{equation} \label{eq:covariance-gff-2}
\widetilde{E}_y (\widetilde{L}_{H_{\rho}}^n (z_i))
= \frac{d(y, w_2)}{d(w_1, w_2)} \widetilde{E}_{w_1} (\widetilde{L}_{H_{\rho}}^n (z_i)) 
+ \frac{d(y, w_1)}{d(w_1, w_2)} \widetilde{E}_{w_2} (\widetilde{L}_{H_{\rho}}^n (z_i)).
\end{equation}
(In (\ref{eq:covariance-gff-2}), we have used the equality
$\widetilde{E}_y (\widetilde{L}_{H_{w_1} \wedge H_{w_2}}^n (z_i)) = 0$ for each $i \in \{1, 2 \}$.)
Note that by (\ref{eq:conti-disc-local-time}) and \cite[Lemma 2.1]{DLP}, for all $i, j \in \{1, 2 \}$,
\begin{equation} \label{eq:covariance-gff-3}
\widetilde{E}_{w_i} (\widetilde{L}_{H_{\rho}}^n (z_j))
= E_{w_i} (L_{H_{\rho}}^n (z_j))
= d(\rho, w_i) + d(\rho, z_j) - d(w_i, z_j).
\end{equation}
We have 
\begin{equation} \label{eq:covariance-gff-4}
\widetilde{E}_x (\widetilde{L}_{H_{z_1} \wedge H_{z_2}}^n (y)) = 0,
~~~~\text{if}~I_{\{z_1, z_2 \}} \neq I_{\{w_1, w_2 \}}.
\end{equation}
Assume that $I_{\{z_1, z_2 \}} = I_{\{w_1, w_2 \}}$.
Note that a Brownian motion on $\widetilde{T}_{\leq n}$
starting at $x$ killed upon $H_{z_1} \wedge H_{z_2}$
has the same law as
a standard Brownian motion on $(0, \frac{1}{2})$ stating at $d(x, z_1)$ killed upon $H_0 \wedge H_{1/2}$.
Let $\{L_t^{\text{B}} (z) : t \ge 0 \}$ be a local time at $z$
of a standard Brownian motion $\{B_t, t \ge 0, P_w^B, w \in \mathbb{R} \}$ on $\mathbb{R}$.
By this and \cite[Chapter VI, Exercise 2.8]{RY}, we have
\begin{align} \label{eq:covariance-gff-5}
\widetilde{E}_x (\widetilde{L}_{H_{z_1} \wedge H_{z_2}}^n (y)) 
&= E_{d(x, z_1)}^B (L_{H_0 \wedge H_{1/2}}^B (d(y, z_1))) \notag \\
&= \frac{d(x, z_2)}{d(z_1, z_2)} d(y, z_1) + \frac{d(x, z_1)}{d(z_1, z_2)} d(y, z_2) - d(x, y),
~~\text{if}~I_{\{z_1, z_2 \}} = I_{\{w_1, w_2 \}}.
\end{align}
Thus, by (\ref{eq:covariance-gff-1})-(\ref{eq:covariance-gff-5}),
we have
$$\mathbb{E} (\widetilde{h}_x \widetilde{h}_y)
= d(\rho, x) + d(\rho, y) - d(x, y).~~~~~\Box$$ 

\subsection{Proof of Lemma \ref{lem:markov-local-time}} \label{subsec:pf-markov-local-time}
In this section, we use notation in Section \ref{sec:metric-tree}.
Let $\widetilde{X} = (\widetilde{X}_t, t \geq 0, \widetilde{P}_{x}, x \in \widetilde{T}_{\leq n})$
be a Brownian motion on $\widetilde{T}_{\leq n}$.
Let $\{\widetilde{h}_x : x \in \widetilde{T}_{\leq n} \}$ be a centered Gaussian process
associated with $\widetilde{X}$ (see Section \ref{subsec:pf-cov-gff} for the definition).
In the proof of Lemma \ref{lem:markov-local-time}, we will use
a variation of Theorem \ref{thm:iso-thm}:
\begin{thm} (\cite[Theorem 1.1]{EKMRS}) \label{thm:variant-iso-thm}
Fix $n \in \mathbb{N}$.
For any $c \in \mathbb{R}$ and $t > 0$,
\begin{align*}
&\text{the law of}~\left\{\widetilde{L}_{\widetilde{\tau} (t)}^n (x) 
+ \frac{1}{2}(\widetilde{h}_x + c)^2 : x \in \widetilde{T}_{\leq n} \right \}
~\text{under}~\widetilde{P}_{\rho} \times \mathbb{P}
~\text{is the same as} \\
&\text{that of}~\left\{\frac{1}{2} (\widetilde{h}_x + \sqrt{2t + c^2})^2 : x \in \widetilde{T}_{\leq n} \right \}
~\text{under}~\mathbb{P}.
\end{align*}
\end{thm}
{\it Proof of Lemma \ref{lem:markov-local-time}.}
Fix $n \in \mathbb{N}$, $t > 0$, and $a \in T_{\leq n} \backslash T_n$.
Let $\widetilde{X}^{\downarrow} 
= (\widetilde{X}_s^{\downarrow}, s \geq 0, \widetilde{P}_a^{\downarrow})$
be a Brownian motion on $\widetilde{T}_{\leq n - |a|}^a$ starting at $a$.
Let $\{\widetilde{L}_s^{\downarrow} (x) : (s, x) \in [0, \infty) \times \widetilde{T}_{\leq n - |a|}^a \}$
be a local time of $\widetilde{X}^{\downarrow}$.
Set $\widetilde{\tau}^{\downarrow} (s) := \inf \{r \geq 0 : \widetilde{L}_r^{\downarrow} (a) > s \}$.
By Lemma \ref{lem:covariance-gff}, we have
\begin{align} \label{eq:markov-gff}
&\left\{\widetilde{h}_x : x \in \widetilde{T}_{\leq n} \right\} 
~\text{under}~\mathbb{P} \notag \\
&\stackrel{d}{=} 
\left\{\widetilde{h}_x^{\downarrow} + \widetilde{h}_a :  
x \in \widetilde{T}_{\leq n - |a|}^a \right\}
\cup \left\{\widetilde{h}_x : x \in \widetilde{T}_{\leq n} \backslash \widetilde{T}_{\leq n - |a|}^a \right \}
~\text{under}~\mathbb{P}^{\downarrow} \times \mathbb{P},
\end{align}
where 
$\{\widetilde{h}_x^{\downarrow} : x \in \widetilde{T}_{\leq n - |a|}^a \}$
is a Gaussian process associated with $\widetilde{X}^{\downarrow}$ 
under $\mathbb{P}^{\downarrow}$.

Fix $m, r \in \mathbb{N}$, $\lambda_1, \cdots, \lambda_m \in [0, \infty)$, 
$\mu_1, \cdots, \mu_r \in [0, \infty)$,
$x_1, \cdots, x_m \in \widetilde{T}_{\leq n - |a|}^a \backslash \{a \}$,
and $y_1, \cdots, y_r \in (\widetilde{T}_{\leq n} \backslash \widetilde{T}_{\leq n - |a|}^a) \cup \{a \}$.
Let $E$, $\widetilde{E}$, $E^{\downarrow}$ be expectations of the product measures
$\widetilde{P}_{\rho} \times \widetilde{P}_a^{\downarrow} \times \mathbb{P}$, 
$\widetilde{P}_{\rho} \times \mathbb{P}$,
$\widetilde{P}_a^{\downarrow} \times \mathbb{P}^{\downarrow}$, respectively.
By (\ref{eq:markov-gff}), we have
\begin{align} \label{eq:lem-markov-local-time-1}
&E \left[e^{- \sum_{i=1}^m \lambda_i 
\left(\widetilde{L}_{\widetilde{\tau}^{\downarrow} \left(\widetilde{L}_{\widetilde{\tau} (t)}^n (a) \right)}
^{\downarrow} (x_i) 
+ \frac{1}{2}(\widetilde{h}_{x_i})^2 \right)}
\cdot e^{- \sum_{i=1}^r \mu_i \left(\widetilde{L}_{\widetilde{\tau} (t)}^n (y_i) + 
\frac{1}{2}(\widetilde{h}_{y_i})^2 \right)} \right] \notag \\
&= \widetilde{E} \left[e^{- \sum_{i=1}^r \mu_i \left(\widetilde{L}_{\widetilde{\tau} (t)}^n (y_i) + 
\frac{1}{2}(\widetilde{h}_{y_i})^2 \right)}
E^{\downarrow} \left[e^{- \sum_{i=1}^m \lambda_i 
\left(\widetilde{L}_{\widetilde{\tau}^{\downarrow} \left(\widetilde{L}_{\widetilde{\tau} (t)}^n (a) \right)}
^{\downarrow} (x_i) 
+ \frac{1}{2}(\widetilde{h}^{\downarrow}_{x_i} + \widetilde{h}_a)^2 \right)} \right] \right].
\end{align}
Applying Theorem \ref{thm:variant-iso-thm} to the local time and the associated Gaussian process
on $\widetilde{T}_{\leq n - |a|}^a$, the right of (\ref{eq:lem-markov-local-time-1})
is equal to
\begin{equation} \label{eq:lem-markov-local-time-2}
\widetilde{E} \left[e^{- \sum_{i=1}^r \mu_i \left(\widetilde{L}_{\widetilde{\tau} (t)}^n (y_i) + 
\frac{1}{2}(\widetilde{h}_{y_i})^2 \right)}
\mathbb{E}^{\downarrow} \left[e^{- \sum_{i=1}^m \lambda_i 
\left(\frac{1}{2}\left(\widetilde{h}_{x_i}^{\downarrow}
+ \sqrt{2\widetilde{L}_{\widetilde{\tau} (t)}^n (a) + (\widetilde{h}_a)^2} \right)^2 \right)} \right] \right].
\end{equation}
Again, by Theorem \ref{thm:variant-iso-thm},
the expectation of (\ref{eq:lem-markov-local-time-2}) is equal to
\begin{equation} \label{eq:lem-markov-local-time-3}
\mathbb{E} \left[e^{- \sum_{i=1}^r \mu_i \left(\frac{1}{2}(\widetilde{h}_{y_i} + \sqrt{2t})^2 \right)}
\mathbb{E}^{\downarrow} \left[e^{- \sum_{i=1}^m \lambda_i 
\left(\frac{1}{2}\left(\widetilde{h}_{x_i}^{\downarrow}
+ \sqrt{(\widetilde{h}_a + \sqrt{2t})^2} \right)^2 \right)} \right] \right].
\end{equation}
By the symmetry of $\widetilde{h}^{\downarrow}$,
conditioned on $\widetilde{h}_a$,
the law of
$(\widetilde{h}^{\downarrow} + |\widetilde{h}_a + \sqrt{2t}|)^2$
is the same as that of
$(\widetilde{h}^{\downarrow} + \widetilde{h}_a + \sqrt{2t})^2$.
By this and (\ref{eq:markov-gff}), the expectation in (\ref{eq:lem-markov-local-time-3}) is equal to
\begin{equation} \label{eq:lem-markov-local-time-4}
\mathbb{E} \left[e^{- \sum_{i=1}^r \mu_i \left(\frac{1}{2}(\widetilde{h}_{y_i} + \sqrt{2t})^2 \right)}
\cdot e^{- \sum_{i=1}^m \lambda_i 
\left(\frac{1}{2}\left(\widetilde{h}_{x_i} + \sqrt{2t} \right)^2 \right)} \right].
\end{equation}
By Theorem \ref{thm:variant-iso-thm}, the expectation in (\ref{eq:lem-markov-local-time-4})
is equal to
\begin{equation} \label{eq:lem-markov-local-time-5}
\widetilde{E} \left[e^{- \sum_{i=1}^r \mu_i 
\left(\widetilde{L}_{\widetilde{\tau} (t)}^n (y_i) + \frac{1}{2}(\widetilde{h}_{y_i})^2 \right)}
\cdot e^{- \sum_{i=1}^m \lambda_i 
\left(\widetilde{L}_{\widetilde{\tau} (t)}^n (x_i) + \frac{1}{2}\left(\widetilde{h}_{x_i}\right)^2 \right)} \right].
\end{equation}
The above argument implies that 
\begin{align*} 
&\left\{\widetilde{L}_{\widetilde{\tau}^{\downarrow} \left(\widetilde{L}_{\widetilde{\tau} (t)}^n (a) \right)}
^{\downarrow} (x) : x \in \widetilde{T}_{\leq n - |a|}^a \right\}
\cup
\left\{\widetilde{L}_{\widetilde{\tau}(t)}^n (x) : x \in \widetilde{T}_{\leq n} \backslash 
\widetilde{T}_{\leq n - |a|}^a \right \}
~\text{under}~\widetilde{P}_{\rho} \times \widetilde{P}_a^{\downarrow} \notag \\
&\stackrel{d}{=}
 \left\{\widetilde{L}_{\widetilde{\tau} (t)}^n (x) : x \in \widetilde{T}_{\leq n} \right\}
~\text{under}~\widetilde{P}_{\rho}.
\end{align*}
This yields the statement of Lemma \ref{lem:markov-local-time}. $\Box$

\subsection{Tail of maximum of local time revisited} \label{subsec:tail-max-revisit}
In the proof of Proposition \ref{prop:geometry}, we need a version of
Proposition \ref{prop:tail}(ii):
\begin{prop} \label{prop:tail-max-revisit}
There exist $c_1 > 0$ and $t_0 > 0$ such that 
for all $n \in \mathbb{N}$, $t \geq t_0$, and $y \in [0,~2 \sqrt{n}]$,
\begin{equation} \label{eq:prop-tail-max-revisit-statement}
P_{\rho} \left(\max_{v \in T_n} \sqrt{L_{\tau (t)}^n (v)} \geq \sqrt{t} + a_n (t) + y \right)
\geq c_1 e^{- 2 \sqrt{\log b}~y}.
\end{equation}
\end{prop}
\begin{rem}
The assumption of $t$ in Proposition \ref{prop:tail-max-revisit}
is weaker than that of Proposition \ref{prop:tail}(ii).
This is the main requirement in the proof of Proposition \ref{prop:geometry}.
\end{rem}
Fix $\varepsilon \in (0,~1/4)$.
For $n \in \mathbb{N}$, $t > 0$, $y > 0$, $0 < r < n$, and $v \in T_n$,
set
\begin{equation*}
A_{v, r}^n (t) := 
\left \{ 
\begin{minipage}{250pt}
$\sqrt{\widetilde{L}_{\widetilde{\tau} (t)}^n (v_s)} 
\in \sqrt{t} + \ell_{y, t, n} (s)  + [- g_n (s), ~- f_n(s)),
~ \forall s \in [r,~ n-r]$, \\
$\sqrt{\widetilde{L}_{\widetilde{\tau} (t)}^n (v)}
\in [\sqrt{t} + a_n (t) + y,~\sqrt{t} + a_n (t) + y + 1)$, \\
$\sqrt{\widetilde{L}_{\widetilde{\tau} (t)}^n (v_{s^{\prime}})}
\geq \sqrt{t} + \ell_{y, t, n} (s^{\prime}) - r^{1/2 + 2 \varepsilon},
~\forall s^{\prime} \in [0,~r] \cup [n-r,~n]$
\end{minipage}
\right \},
\end{equation*}
where 
\begin{equation*}
\ell_{y, t, n} (s) := \left(\frac{a_n (t)}{n} + \frac{y}{n} \right)s,~~~~~0 \leq s \leq n,
\end{equation*}
\begin{equation*}
f_n (s) := \min \left\{s^{1/2 - \varepsilon},~(n-s)^{1/2 - \varepsilon} \right\},
~~~g_n (s) := \min \left\{s^{1/2 + \varepsilon},~(n-s)^{1/2 + \varepsilon} \right\}.
\end{equation*}
To prove Proposition \ref{prop:tail-max-revisit},
we will apply the second moment method
to $\sum_{v \in T_n} 1_{A_{v, r}^n (t)}$.
We first need the following:
\begin{lem} \label{lem:tail-max-revisit-1}
There exist $r_0 \in \mathbb{N}$ and $c_1 > 0$ such that
for all $r_0 \leq r \leq n/4$, $t > 4 r^{1 + 4 \varepsilon}$,
$y \in [0, ~2 \sqrt{n}]$, and $v \in T_n$,
\begin{equation} \label{eq:lem-tail-max-revisit-1-statement}
\widetilde{P}_{\rho} (A_{v, r}^n (t)) \geq c_1 b^{-n} e^{- 2 \sqrt{\log b}~y}.
\end{equation}
\end{lem}
{\it Proof.}
Fix $r_0 \leq r \leq n/4$, $t > 4r^{1+4\varepsilon}$, and $y \in [0, 2 \sqrt{n}]$,
where we take $r_0 \in \mathbb{N}$ large enough.
By Lemma \ref{lem:local-time-bessel} and
(\ref{eq:0-1-dim-bessel}),
$\widetilde{P}_{\rho} (A_{v, r}^n (t))$ is bounded from below by
\begin{equation} \label{eq:lem-tail-max-revisit-1-1}
c_1 \sqrt{\frac{\sqrt{t}}{\sqrt{t} + a_n (t) + y + 1}}
P_0^B \left[
\begin{minipage}{190pt}
$- g_n (s) \leq B_s - \ell_{y, t, n} (s) < - f_n (s),
~\forall s \in [r,~ n-r]$, \\
$B_n \in [a_n (t) + y,~a_n (t) + y + 1)$, \\
$B_{s^{\prime}} \geq \ell_{y, t, n} (s^{\prime}) - r^{1/2 + 2 \varepsilon},
~\forall s^{\prime} \in [0,~r] \cup [n-r,~n]$
\end{minipage}
\right],
\end{equation}
where we have used the following:
under the event that
$\sqrt{X_s/2} \geq \sqrt{t} + \ell_{y, t, n} (s) - g_n (s)$ for all $s \in [r,~n-r]$
and
$\sqrt{X_{s^{\prime}}/2} \geq \sqrt{t} + \ell_{y, t, n} (s^{\prime}) - r^{1/2 + 2 \varepsilon}$
for all $s^{\prime} \in [0,~r] \cup [n-r,~n]$,
we have
$\exp \left(- \frac{3}{8} \int_0^n \frac{ds}{X_s} \right) \geq c_1$.
Since $\left\{B_s - B_n \frac{s}{n} : 0 \leq s \leq n \right\}$ is independent of $B_n$
and has the same law as that of a Brownian bridge from $0$ to $0$ on $[0, ~n]$,
the right-hand side of (\ref{eq:lem-tail-max-revisit-1-1}) is bounded from below by
\begin{align} \label{eq:lem-tail-max-revisit-1-2}
&c_1 \sqrt{\frac{\sqrt{t}}{\sqrt{t} + a_n (t) + y + 1}}
P_{0 \to 0}^n
\left[
\begin{minipage}{175pt}
$- g_n (s) \leq X_s < - f_n (s) - 1,~\forall s \in [r,~ n-r]$, \\
$X_{s^{\prime}} \geq - r^{1/2 + 2 \varepsilon},
~\forall s^{\prime} \in [0,~r] \cup [n-r,~n]$
\end{minipage}
\right] \notag \\
&
\times P_0^B \left(B_n \in [a_n (t) + y,~a_n (t) + y + 1] \right),
\end{align} 
where for $T > 0$ and $p, q \in \mathbb{R}$, $P_{p \to q}^T$ is a probability law
on $\left(C[0, T], \mathcal{B}(C[0, T]) \right)$
($C[0, T]$ is the space of all continuous functions on $[0, T]$
and $\mathcal{B}(C[0, T])$ is the $\sigma$-field generated by cylinder sets in $C[0, T]$)
under which the coordinate process $\{X_s : 0 \leq s \leq T \}$ 
is a Brownian bridge with variance $1/2$ from $p$ to $q$ on $[0,~T]$. 
We have
\begin{align} \label{eq:lem-tail-max-revisit-1-4}
&~~~P_{0 \to 0}^n
\left[
\begin{minipage}{175pt}
$- g_n (s) \leq X_s < - f_n (s) - 1,~\forall s \in [r,~ n-r]$, \\
$X_{s^{\prime}} \geq - r^{1/2 + 2 \varepsilon},
~\forall s^{\prime} \in [0,~r] \cup [n-r,~n]$
\end{minipage}
\right] \notag \\
&= E_{0 \to 0}^n \left[
\begin{minipage}{253pt}
$1_{\left\{- g_n (r) \leq X_r < - f_n (r) - 1,
~ - g_n (n-r) \leq X_{n-r} < - f_n (n-r) - 1\right \}}$ \\
$\times P_{X_r \to X_{n-r}}^{n-2r}
\left(- g_n (s + r) \leq X_s < - f_n (s + r) - 1, \forall s \in [0,~ n-2r] \right)$ \notag \\ 
$\times P_{0 \to X_r}^r \left(X_s \geq - r^{1/2 + 2 \varepsilon},
~\forall s \in [0,~ r] \right)$ \\
$\times P_{X_{n-r} \to 0}^r \left(X_s \geq - r^{1/2 + 2 \varepsilon},
~\forall s \in [0, ~r] \right)$
\end{minipage}
\right] \notag \\
&\geq c_2 P_{0 \to 0}^n
\left(- g_n (s) \leq X_s < - f_n (s) - 1,~\forall s \in [r,~ n-r] \right),
\end{align}
where we have used \cite[Lemma 2.2(a)]{Br2} in the second inequality
which implies that
the last two probabilities in the expectation in the second display of 
(\ref{eq:lem-tail-max-revisit-1-4})
are bounded from below by some constants.
By the symmetry of a Brownian bridge,
\cite[Lemma 2.7 and Proposition 6.1]{Br2},
we have for $r_0 \leq r < n/4$ with $r_0$ large enough
\begin{equation} \label{eq:lem-tail-max-revisit-1-5}
P_{0 \to 0}^n
\left(- g_n (s) \leq X_s < - f_n (s) - 1,~\forall s \in [r,~ n-r] \right)
\geq c_3 P_{0 \to 0}^n
\left(X_s \geq 0,~\forall s \in [r,~ n-r] \right).
\end{equation}
Using \cite[Lemma 2.2(a)]{Br2},
one can show that the right of (\ref{eq:lem-tail-max-revisit-1-5})
is bounded from below by $c_4/n$ for some $c_4 > 0$.
Thus, by (\ref{eq:lem-tail-max-revisit-1-1})-(\ref{eq:lem-tail-max-revisit-1-5}), 
we have (\ref{eq:lem-tail-max-revisit-1-statement}).
~~~$\Box$ \\
Next, we need the following:
\begin{lem} \label{lem:tail-max-revisit-2}
(i) There exist $r_0 \in \mathbb{N}$ and $c_1 > 0$ such that
for all $t > 0$, $r_0 \leq r \leq n/4$, $y \in [0, ~2\sqrt{n}]$,
$k \in \{r, r+1, \dotsc, n-r-1 \}$, $v \in T_n$, and 
$\left(\sqrt{t} + \ell_{y, t, n} (k) - g_n (k) \right)^2 \leq q 
\leq \left(\sqrt{t} + \ell_{y, t, n} (k) - f_n (k) \right)^2$,
\begin{align} \label{eq:lem-tail-max-revisit-2-statement}
&\widetilde{P}_{v_{k}} \left[
\begin{minipage}{210pt}
$\sqrt{\widetilde{L}_{\widetilde{\tau}^{\downarrow} (q)}^{\downarrow} (v_s)}
< \sqrt{t} + \ell_{y, t, n} (s) - f_n (s),~\forall s \in [k,~ n-r]$, \\
$\sqrt{\widetilde{L}_{\widetilde{\tau}^{\downarrow} (q)}^{\downarrow} (v)}
\in [\sqrt{t} + a_n (t) + y,~\sqrt{t} + a_n (t) + y + 1)$ 
\end{minipage}
\right] \notag \\
&\leq c_1 r^{1/2 + \varepsilon} \sqrt{\frac{\sqrt{t} + \ell_{y, t, n} (k) - f_n (k)}
{\sqrt{t} + a_n (t) + y}}
\frac{g_n (k)}{(n-k-r) \sqrt{n-k}}
b^{-(n-k)} \notag \\
&~~~~~~~~~~~~~~~~~~\times e^{\frac{3 \log n}{2n} (n-k)}
e^{\frac{\log \left(\frac{\sqrt{t} + n}{\sqrt{t}} \right)}{2n} (n-k)}
e^{- \frac{2 a_n (t)}{n} f_n (k)},
\end{align}
where $\{\widetilde{L}_s^{\downarrow} (x) : 
(s, x) \in [0,~\infty) \times \widetilde{T}_{\leq n-k}^{v_k} \}$
is a local time of a Brownian motion on $\widetilde{T}_{\leq n-k}^{v_k}$
and $\widetilde{\tau}^{\downarrow} (p) := \inf \{s \geq 0 : \widetilde{L}_s^{\downarrow} (v_k) > p \}$.
\\
(ii) There exist $r_0 \in \mathbb{N}$ and $c_2 > 0$ such that
for all $t > 0$, $r_0 \leq r \leq n/4$, $y \in [0, 2 \sqrt{n}]$, and $v \in T_n$,
\begin{equation*}
\widetilde{P}_{\rho} \left(A_{v, r}^n (t) \right) \leq c_2 r b^{-n} e^{- 2 \sqrt{\log b} y}.
\end{equation*} 
\end{lem}
{\it Proof.}
Recall the probability measure $P_{0 \to 0}^T$ defined 
in the proof of Lemma \ref{lem:tail-max-revisit-1}.
By Lemma \ref{lem:local-time-bessel} and
(\ref{eq:0-1-dim-bessel}), the left of (\ref{eq:lem-tail-max-revisit-2-statement})
is bounded from above by
\begin{equation} \label{eq:lem-tail-max-revisit-2-1}
\sqrt{\frac{\sqrt{q}}{\sqrt{t} + a_n (t) + y}}
P_0^B \left[
\begin{minipage}{230pt}
$B_s < \sqrt{t} - \sqrt{q} + \ell_{y, t, n} (k+s) - f_n (k+s),
~\forall s \in [0,~ n - k - r]$,  \\
$B_{n-k} \in [\sqrt{t} - \sqrt{q} + a_n (t) + y,~\sqrt{t} - \sqrt{q} + a_n (t) + y + 1)$
\end{minipage}
\right].
\end{equation}
Since the process 
$\left\{B_s - \frac{s}{n-k} \cdot B_{n-k} : 0 \leq s \leq n-k \right\}$
is independent of $B_{n-k}$ and has the same law as that of a Brownian bridge with variance $1/2$
from $0$ to $0$ on $[0,~n-k]$,
the probability in (\ref{eq:lem-tail-max-revisit-2-1}) is bounded from above by
\begin{align} \label{eq:lem-tail-max-revisit-2-1-2}
&P_{0 \to 0}^{n-k} \left[
X_s < \sqrt{t} - \sqrt{q} + \ell_{y, t, n} (k+s)
-\frac{s}{n-k} (\sqrt{t} - \sqrt{q} + a_n (t) + y), 
~\forall s \leq n - k - r \right] \notag \\
&\times P_0^B \left[
B_{n-k} \in [\sqrt{t} - \sqrt{q} + a_n (t) + y,~\sqrt{t} - \sqrt{q} + a_n (t) + y + 1) \right].
\end{align}
We estimate two probabilities in (\ref{eq:lem-tail-max-revisit-2-1-2}).
By the assumption of $q$, we have
\begin{align} \label{eq:lem-tail-max-revisit-2-2}
&~~~P_0^B \left(
B_{n-k} \in [\sqrt{t} - \sqrt{q} + a_n (t) + y,~\sqrt{t} - \sqrt{q} + a_n (t) + y + 1) \right)
\notag \\
&\leq \frac{c_1}{\sqrt{n-k}} b^{-(n-k)}
e^{\frac{3 \log n}{2n} (n-k)}
e^{\frac{\log \left(\frac{\sqrt{t} + n}{\sqrt{t}} \right)}{2n} (n-k)}
e^{- \frac{2 a_n (t)}{n} f_n (k)}.
\end{align}
To estimate the other probability, we use \cite[Lemma 3.4]{ABK1}:
for any $x_1, x_2 \in [0, \infty)$, $r_1, r_2 \in [0, \infty)$, and $T > r_1 + r_2$,
\begin{align} \label{eq:version-linear-barrier} 
&~~~P_{0 \to 0}^T \left(
X_s \leq \left(1 - \frac{s}{T} \right) x_1 + \frac{s}{T} x_2,
~r_1 \leq s \leq T - r_2 \right) \notag \\
&\leq \frac{2}{T - r_1 - r_2}
\left \{\left(1 - \frac{r_1}{T} \right) x_1 + \frac{r_1}{T} x_2 + \sqrt{r_1} \right \}
\left \{\frac{r_2}{T} x_1 + \left(1 - \frac{r_2}{T} \right) x_2 + \sqrt{r_2} \right \}.
\end{align}
By (\ref{eq:version-linear-barrier}), we have
\begin{align} \label{eq:lem-tail-max-revisit-2-3}
&P_{0 \to 0}^{n-k} \left(
X_s < \sqrt{t} - \sqrt{q} + \ell_{y, t, n} (k+s)
-\frac{s}{n-k} (\sqrt{t} - \sqrt{q} + a_n (t) + y),~
\forall s \leq n - k - r \right) \notag \\
&\leq c_2 \frac{r^{1/2 + \varepsilon}}{n-k-r} g_n (k).
\end{align}
Thus, by (\ref{eq:lem-tail-max-revisit-2-1})-(\ref{eq:lem-tail-max-revisit-2-3}),
we have (\ref{eq:lem-tail-max-revisit-2-statement}).
By repeating a similar argument,
we can also prove (ii). We omit the detail.~~~$\Box$ \\\\
{\it Proof of Proposition \ref{prop:tail-max-revisit}.}
Fix a sufficiently large positive constant $r$.
Set $Z := \sum_{v \in T_n} 1_{A_{v, r}^n (t)}$.
For each $k \in \{0, \cdots, n-1 \}$, set
$S_k := \sum_{\begin{subarray}{c} v, u \in T_n, \\ |v \wedge u| = k \end{subarray}}
\widetilde{P}_{\rho} \left[A_{v, r}^n (t) \cap A_{u, r}^n (t) \right]$.
We have
\begin{equation} \label{eq:prop-tail-max-revisit-1}
\widetilde{E}_{\rho} [Z^2]
= \widetilde{E}_{\rho} [Z]
+ \sum_{k = 0}^{r-1} S_k 
+ \sum_{k=r}^{n - r - 1} S_k
+ \sum_{k = n-r}^{n-1} S_k.
\end{equation}
For each $0 \leq k \leq n - 1$, set 
$$I_{y, t, n} (k) := [(\sqrt{t} + \ell_{y, t, n} (k) - g_n (k))^2, (\sqrt{t} + \ell_{y, t, n} (k) - f_n (k))^2].$$
For each $q \in I_{y, t, n} (k)$ and $v \in T_n$, let $P_{v_k} (q)$ be the probability in 
(\ref{eq:lem-tail-max-revisit-2-statement}).
In the estimate of $\sum_{k=r}^{n - r - 1} S_k$,
we use the following: by Lemma \ref{lem:markov-local-time},
for any $r \leq k \leq n - r - 1$ and $v, u \in T_n$
with $|v \wedge u| = k$,
$\widetilde{P}_{\rho} \left[A_{v, r}^n (t) \cap A_{u, r}^n (t) \right]$
is bounded from above by 
\begin{equation} \label{eq:prop-tail-max-revisit-2}
\widetilde{P}_{\rho} [A_{u, r}^n (t)] \cdot
\sup_{q \in I_{y, t, n} (k)} P_{v_k} (q).
\end{equation} 
By (\ref{eq:prop-tail-max-revisit-2}) and Lemma \ref{lem:tail-max-revisit-2},
dividing the sum
$\sum_{k=r}^{n - r - 1} S_k$
into $\sum_{k = r}^{\lfloor n/2 \rfloor} S_k + \sum_{\lfloor n/2 \rfloor + 1}^{\lfloor n - n/\log n \rfloor} S_k
+ \sum_{k = \lfloor n - n/\log n \rfloor + 1}^{n - r - 1} S_k$, we have
$\sum_{k=r}^{n - r - 1} S_k \leq c_1 e^{- 2\sqrt{\log b}~y}$.
Similarly, by Lemma \ref{lem:markov-local-time}, for any $0 \leq k \leq r - 1$
and $v, u \in T_n$ with $|v \wedge u| = k$,
$\widetilde{P}_{\rho} \left[A_{v, r}^n (t) \cap A_{u, r}^n (t) \right]$
is bounded from above by 
$\widetilde{P}_{\rho} [A_{u, r}^n (t)] \cdot \sup_{q \in I_{y, t, n} (r)} P_{v_r} (q)$.
By this and Lemma \ref{lem:tail-max-revisit-2},
$\sum_{k = 0}^{r-1} S_k \leq c_2 e^{- 2 \sqrt{\log b}~y}$.
For any $n - r \leq k \leq n - 1$ and $v, u \in T_n$ with $|v \wedge u| = k$,
we bound $\widetilde{P}_{\rho} \left[A_{v, r}^n (t) \cap A_{u, r}^n (t) \right]$
just by $\widetilde{P}_{\rho} \left[A_{u, r}^n (t) \right]$.
Then, by Lemma \ref{lem:tail-max-revisit-2}(ii), 
we have $\sum_{k = n-r}^{n-1} S_k \leq c_3 e^{- 2 \sqrt{\log b}~y}$.

By these estimates, Lemma \ref{lem:tail-max-revisit-1},
and an argument similar to the proof of Proposition \ref{prop:tail}(ii),
we obtain the desired result (\ref{eq:prop-tail-max-revisit-statement}).
$\Box$

\subsection{Tail of maximum of BRW} \label{subsec:tail-max-brw}
In the proof of Proposition \ref{prop:limit-tail}, we use tail estimates of the maximum
of the BRW on $T$.
Let $(h_v)_{v \in T}$ be a BRW on $T$ defined in Section \ref{sec:intro}.
\begin{lem} \label{lem:tail-max-brw-1}
(i) There exist $c_1, c_2 \in (0, ~\infty)$ such that for all $y > 0$ and $n \in \mathbb{N}$,
\begin{equation*}
\mathbb{P} \left(\max_{v \in T_n} h_v
> \sqrt{\log b} ~n - \frac{3}{4 \sqrt{\log b}} \log n + y \right)
\leq c_1 (1 + y) e^{- 2 \sqrt{\log b}~y} e^{- c_2 \frac{y^2}{n}}.
\end{equation*}
(ii) There exist $c_3 > 0$ and $n_0 \in \mathbb{N}$ such that for all $n \geq n_0$
and $y \in [1, ~\sqrt{n}]$,
\begin{equation*}
\mathbb{P} \left(\max_{v \in T_n} h_v
> \sqrt{\log b}~ n - \frac{3}{4 \sqrt{\log b}} \log n + y \right)
\geq c_3 y e^{- 2 \sqrt{\log b}~y}.
\end{equation*}
\end{lem}
Lemma \ref{lem:tail-max-brw-1}(ii) is a special version of Lemma 2.7 of \cite{BDZ2}.
One can easily modify the proof of Lemma 3.8 in \cite{BDZ}
(this is basically a tail estimate of the maximum of a BRW on a $4$-ary tree)
to prove
Lemma \ref{lem:tail-max-brw-1}(i). We omit the details.
\\\\
{\bf Acknowledgments.} \\
The author would like to thank Professor Kumagai and Dr. Kajino for valuable comments
and encouragement,
and Professor Biskup for stimulating discussion and pointing out
a connection between my work and the theory of random multiplicative cascade measures.
This work is partially supported by JSPS KAKENHI 13J01411
and University Grants for student exchange between universities in partnership
under Top Global University Project of Kyoto University.

\end{document}